\documentclass[12pt]{amsart}
\usepackage{latexsym, amsmath, amssymb, amsthm}
\usepackage{epsfig}
\usepackage{amscd}
\usepackage{latexsym}
\usepackage{pstricks,pst-node,cite}

\makeatletter
\newtheorem*{rep@theorem}{\rep@title}
\newcommand{\newreptheorem}[2]{%
\newenvironment{rep#1}[1]{%
 \def\rep@title{#2 \ref{##1}}%
 \begin{rep@theorem}}%
 {\end{rep@theorem}}}
\makeatother

\newtheorem{theorem}{Theorem}[section]
\newreptheorem{theorem}{Theorem}

\newtheorem{cor}[theorem]{Corollary}
\newreptheorem{cor}{Corollary}
\newtheorem{exa}[theorem]{Example}
\newtheorem{que}[theorem]{Question}
\newtheorem{defi}[theorem]{Definition}
\newtheorem{remark}[theorem]{Remark}
\newtheorem{notation}[theorem]{Notation}

\newcommand{\A}{\mathcal{A}}

\newcommand{\C}{\mathcal{C}}
\newcommand{\BS}{\mathbb{S}}

\newcounter{fignum}
\psset{linecolor=black}\newgray{lightgray}{.9}\newgray{darkgray}{.75}\newgray{pup}{.4}
\psset{linewidth=.6pt,dash=3pt 3pt,doublesep=.05,dotsize=1pt 5}
\SpecialCoor

\begin{document}

\title[An alternating labeling on a spanning tree of Seifert graphs]{An alternating labeling on a spanning tree of Seifert graphs
and applications in knot theory}

\author{Dongseok Kim}
\address{Department of Mathematics \\Kyonggi University
\\ Suwon, 443-760 Korea}
\thanks{This work was supported by Kyonggi University Research Grant 2011.}
\email{dongseok@kgu.ac.kr}

\begin{abstract}
In this article, we prove the existence of a co-tree edge alternating
spanning tree of the Seifert graphs of canonical Seifert surfaces.
As an application, we show the existence of basket, flat plumbing and flat plumbing basket
surfaces of a link from its Seifert surface using the Seifert graph of the canonical Seifert surface.
This generalizes the existence of such surfaces only from the braid presentation of the link.
We define the basket number, flat plumbing number and flat
plumbing basket number of a link. Then we provide several upper bounds
for these plumbing numbers, illustrate our upper bounds are sharper than the
previous bounds from braid presentations and study the
relation between these plumbing numbers and the genera of links.
\end{abstract}


\maketitle

\section{Introduction}

Let $\Gamma$ be a finite simple graph with vertex set $V(\Gamma)$ and
edge set $E(\Gamma)$. One classical problem in graph theory is
to find the \emph{complexity} of $\Gamma$, $\kappa(\Gamma)$,
the number spanning trees in a graph $\Gamma$~\cite{Chow, HK}.
The celebrated Kirchhoff's matrix tree theorem finds that
$\kappa(\Gamma)$ is any cofactor of the admittance matrix
(or Laplacian matrix) of $\Gamma$ which is a generalization of
Cayley's formula which provides
$\kappa(K_n)$ of the complete graph $K_n$ on $n$ vertices.
The spanning trees of $\Gamma$ have many wonderful applications not only in graph
theory but also in several mathematical areas, computer
science and engineering\cite{Cha, Chow, CT, DFKLS, GV, GX, HK, MS, Th}.
On the other hand, ever since the exploratory paper by Dirac \cite{Dirac},
the chromatic number has been in the center of graph theory
research. The chromatic number
$\chi(\Gamma)$ of a graph $\Gamma$ is the smallest number of colors needed
to color the vertices of G so that no two adjacent vertices share
the same color. Its rich history can be found in several articles
\cite{HM, theoremas}.

A few graphs can be found in knot theory; a $4$
valent graph can be obtained from a
knot diagram by making the crossings to double
points~\cite{Jones, Kauffman:graph, Kim:transitive},
Tait checkerboard graphs of link diagrams are used in
knot Floer homology~\cite{BL, CK, Greene, We}
and Seifert graphs obtained from Seifert surfaces will
be used in the present article~\cite{HMV, Manchon, Murasugi:knot}.
Seifert graphs are signed graphs and originally
assumed to be planar~\cite{Murasugi:knot}
but here they may not be planar since we keep \emph{the
rotation scheme}, cyclic orders of edges adjacent to each vertex.

Rudolph~\cite{Rudolph:plumbing} first introduced several plumbing surfaces.
These plumbing surfaces are related with the geometry of knot
complements~\cite{Gabai:genera, Stallings:const}.
The existence of a flat plumbing surface of an arbitrary link was
first found by Harashi and Wada~\cite{HW:plumbing} and
the existence of a flat plumbing basket surface was proven
by Furihata, Hirasawa and Kobayashi~\cite{FHK:openbook}.
Both proofs were based on the Alexander theorem; every link is a closed braid.

The author's first preprint about these plumbing surfaces from a canonical
Seifert surface had a critical mistake pointed out by Kobayashi.
By weakening some conditions of plumbings, the author, Kwon and Lee proved
the existence of banded surfaces and flat banded surfaces~\cite{KKL:string}.
The author also proved that every link $L$ is the boundary of an oriented surface
which is obtained from a graph embedding of a dipole graph,
this surface is also known as a braidzel surface~\cite{Nakamura:braidzel},
and a complete bipartite graph $K_{2,n}$,
where all voltage assignments on the edges of dipole graph and $K_{2,n}$ are $0$~\cite{Kim:dipole}.

The definitions of these plumbing surfaces~\cite{Rudolph:plumbing} are very technical
and so it is difficult to handle but the work in~\cite{FHK:openbook} provided
a tangible equivalent definition of a flat plumbing basket surface using
an open book decomposition. Using this definition and results in~\cite{KKL:string},
Choi, Do and the author are working on a new knot tabulation with respect
to the flat plumbing basket number~\cite{CDK}.
The present work is the beginning of this series of results presenting
links as a boundary of the surface obtained in a
embedding of certain graphs as described in~\cite{GT1}. One might consider
these plumbing surfaces as special embeddings of the bouquets of circles.

For a correct construction of a flat plumbing basket surface from a canonical Seifert surface,
we have to consider the Seifert graph of a canonical Seifert
surface and will find a spanning tree
of the Seifert graph with a special labeling on the spanning tree. A spanning tree $T$ of
a graph $\Gamma$ is called a \emph{co-tree  edge alternating spanning tree} if
for any co-tree  edge $e \in E(\Gamma)\setminus E(T)$,
the unique path $P_e$ in $T$ joining both end vertices of the edge $e$
has an alternating signs with respect to a labeling $\mu$ on $T$.

We provide an algorithm to find a co-tree  edge alternating spanning tree
in the following theorem.

\begin{reptheorem}{maintheorem1}
For a connected bipartite graph $\Gamma$ and a vertex $v$,
there exists an algorithm to determine a co-tree  edge alternating spanning tree $T$
with respect to the depth labeling $\mu_v : E(T) \rightarrow \{ +, - \}$.
\end{reptheorem}

As applications of Theorem~\ref{maintheorem1}, we not only
construct basket surfaces, flat plumbing surfaces
and flat plumbing basket surfaces
but also define the \emph{basket number, flat plumbing number and flat
plumbing basket number} of a link which are the minimum number
of annuli plumbings required to have each of these plumbing surfaces.
Consequently, we provide some upper bounds for these
plumbing numbers from braid presentations of
links and canonical Seifert surfaces of links.
We also compare the upper bounds for those
which were obtained from a braid presentation of a link.
In a recent article by Hirose and Nakashima~\cite{HN},
the flat plumbing basket numbers of knots of $9$ crossings or less
were studied.

The outline of this paper is as follows.
In section~\ref{main},
we first review some preliminary definitions in knot theory,
then we provide the definition of the Seifert graph
$\Gamma(S)$ obtained from a canonical Seifert surface $S$.
In section~\ref{premain} we prove Theorem~\ref{maintheorem1}.
In section~\ref{secbasket}, we first
review the definitions of these plumbing numbers.
In subsection~\ref{baskets} we find a new
upper bound for basket number as follows and
demonstrate this bound is sharp in Example~\ref{basketexa}.

\begin{reptheorem}{baskettheorem1}
Let $L$ be a link which is the closure of a braid $\beta \in B_n$
where the length of the braid $\beta$ is $m$.
Then the basket number of $L$ is less than or equal to $m-n+1$, $i. e.$,
$$bk(L)\le m-n+1.$$
\end{reptheorem}

In subsection~\ref{flatbaskets}, we first find an upper bound
for flat plumbing basket number from the braid
presentative of a link as follows and demonstrate this bound
is better than the previous one in Example~\ref{fpbsexa}.

\begin{reptheorem}{fpbktheoremext}
Let $L$ be an oriented link which is a closed $n$-braid with a braid
word $\beta$ whose length is $m$ and
let $ps(\sigma_i^{\pm 1})$ be the power sum
of $\sigma_i^{\pm 1}$ in $\beta$ for all $i=1, 2,
\ldots, n-1$.
Let $\gamma$ be the cardinality of the set
$$\Omega= \{ i | 1
\le i \le n-1, \sigma_i ~{\rm{and}}~ \sigma_i^{-1}~ {\rm{both}~\rm{appear}~ \rm{in}}~ \beta\}.$$
Let
$$\epsilon_i = \begin{cases} 1 ~~~& {\rm{if}}~ 1 \le
ps(\sigma_i^{1}) \le ps(\sigma_i^{- 1}) ~{\rm{or}}~ps(\sigma_i^{- 1})=0, \\
-1 ~~~&{\rm{if}}~ 1 \le ps(\sigma_i^{-1}) \le ps(\sigma_i^{1})
~{\rm{or}}~ps(\sigma_i^{1})=0. \end{cases}
$$
Then the flat plumbing basket number of $L$ is bounded by
$m + n-1 -4\gamma +2\sum_{i=1}^{n-1} ps(\sigma_i^{\epsilon_i})$, $i. e.,$
$$fpbk(L)\le m + n-1 -4\gamma +2\sum_{i=1}^{n-1} ps(\sigma_i^{\epsilon_i}).$$
\end{reptheorem}

Next we find an upper bound for flat plumbing basket number by
constructing flat plumbing basket surfaces from canonical Seifert surfaces by choosing a
spanning tree and alternating labeling on it provided in
Theorem~\ref{maintheorem1} as follows.
We also demonstrate that our upper bound is sharper than
the previous one from braid presentation in Example~\ref{52flat}
and Example~\ref{fpbsexa2}.

\begin{reptheorem}{flatbktheorem3}
Let $\Gamma$ be an Seifert graph of canonical Seifert surface $S$ of a link $L$ with
$|V(\Gamma)|=n$, $|E(\Gamma)|=m$ and the sign labeling $\phi$.
Let $G(\Gamma)$ be the induced graph of $\Gamma$.
Let $T$ be a co-tree  edge alternating spanning tree of $\Gamma$ and $\mu$ a
labeling on $T$ chosen by Theorem~\ref{maintheorem1}.
Let $\delta(T)$ be the cardinality of the set
$$\Psi(T) =\{e\in E(T)|~\mu(e)\neq\phi(\overline{e}) ~{\rm{for~all}}~ \overline{e} \in \Gamma(e) \},$$
and let $\zeta(T)$ be the cardinality of the set
$$\Upsilon(T) =\{\overline{e}\in E(\Gamma(T))~|~\mu(e)=\phi(\overline{e}),~\overline{e}\in \Gamma(e),~e \in E(T)-\Omega(T) \}.$$
and let $\eta(T)$ be the cardinality of the set
$$\Phi(T) = \{ \overline{e}\in E(\Gamma)-E(\Gamma(T)) ~|~\mu(\overline{e}) = \nu(e) \}$$
where $\nu(e) =+$($-$, resp.) if there is one extra positive(negative, respectively) sign
in the path $P_e$ joining end vertices of the edge $e$ in $T$.
Then the
flat plumbing basket number of $L$ is bounded by $m - 3(n-1) + 2 ( 2\delta(T)+\zeta(T) + \eta(T))$, $i. e.$,
$$ fpbk(L)\le m - 3(n-1) + 2 ( 2\delta(T)+\zeta(T) + \eta(T)).$$
\end{reptheorem}

We also obtain Corollary~\ref{bktheorem3} for an upper bound for the
basket number from a canonical Seifert surface using Theorem~\ref{flatbktheorem3}.
In subsection~\ref{flats}, we find a new
upper bounds for flat plumbing number in Theorem~\ref{fptheorem1} from
braid presentative of $L$ and Theorem~\ref{flattheorem2}
from its canonical Seifert surface.
At last, we study the relations between flat plumbing numbers
and the genera of a link in subsection~\ref{relation}.

\section{Seifert graphs of Seifert surfaces}
\label{main}

We first review some preliminary definitions in knot theory
in subsection~\ref{pregraph}. For more terms in knot theory,
we refer the readers
to~\cite{BZ:knot, Murasugi:knot}.
For general terminology in graph theory, we refer the reader
to~\cite{BM:GTWA, GT1}.
Then we provide the definition of the Seifert graph $\Gamma$
obtained from a Seifert surface $S$ in subsection~\ref{preSeifert}.

\subsection{Preliminaries in knot theory} \label{pregraph}

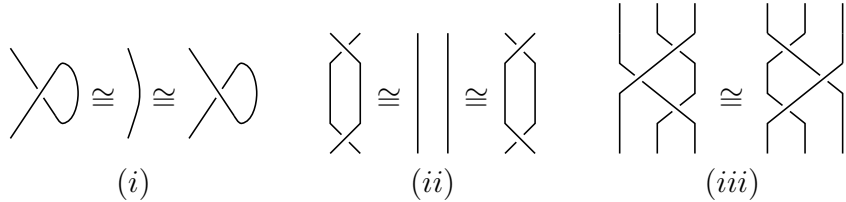
\begin{figure}
$$
\begin{pspicture}[shift=-1.2](-.45,-1.4)(.55,1)
\psline(.2,.3)(-.4,-.6)
\psline(-.4,.6)(-.04,.06) \psline(.04,-.06)(.2,-.3)
\pscurve(.2,-.3)(.3,-.4)(.5,0)(.3,.4)(.2,.3)
\end{pspicture}
\cong
\begin{pspicture}[shift=-1.2](-.25,-1.4)(0,1)
\pscurve(-.2,.6)(-.1,.3)(-.05,.1)(-.05,-.1)(-.1,-.3)(-.2,-.6)
\rput[t](-.125,-1){$(i)$}
\end{pspicture}
\cong
\begin{pspicture}[shift=-1.2](-.45,-1.4)(.55,1)
\psline(-.4,.6)(.2,-.3)
\psline(.2,.3)(.04,.06) \psline(-.04,-.06)(-.4,-.6)
\pscurve(.2,-.3)(.3,-.4)(.5,0)(.3,.4)(.2,.3)
\end{pspicture}
\quad\quad
\begin{pspicture}[shift=-1.2](-.3,-1.4)(.3,1)
\psline(-.2,.8)(.2,.4)(.2,-.4)(-.2,-.8)
\psline(.2,.8)(.05,.65)
\psline(-.05,.55)(-.2,.4)(-.2,-.4)(-.05,-.55)
\psline(.05,-.65)(.2,-.8)
\end{pspicture}
\cong
\begin{pspicture}[shift=-1.2](-.3,-1.4)(.3,1)
\psline(-.2,-.8)(-.2,.8)
\psline(.2,-.8)(.2,.8)
\rput[t](0,-1){$(ii)$}
\end{pspicture}
\cong
\begin{pspicture}[shift=-1.2](-.3,-1.4)(.3,1)
\psline(.2,.8)(-.2,.4)(-.2,-.4)(.2,-.8)
\psline(-.2,.8)(-.05,.65)
\psline(.05,.55)(.2,.4)(.2,-.4)(.05,-.55)
\psline(-.05,-.65)(-.2,-.8)
\end{pspicture}
\quad\quad
\begin{pspicture}[shift=-1](-.7,-1.4)(.7,1)
\psline(-.5,-1)(-.5,-.2)(.5,.6)
\pcline(.5,.6)(.5,1)
\psline(0,-1)(0,-.6)(.2,-.44)
\psline(.3,-.36)(.5,-.2)(.5,.2)(.3,.36)
\psline(.2,.44)(0,.6)
\pcline(0,1)(0,.6)
\psline(.5,-1)(.5,-.6)(-.2,-.04)
\psline(-.3,.04)(-.5,.2)(-.5,.6)
\pcline(-.5,.6)(-.5,1)
\rput[t](1,-1.2){$(iii)$}
\end{pspicture}
\cong
\begin{pspicture}[shift=-1](-.7,-1.4)(.7,1.1)
\psline(-.5,-1)(-.5,-.6)(.5,.2)
\psline(.5,.2)(.5,.6)
\pcline(.5,.6)(.5,1)
\psline(0,-1)(0,-.6)(-.2,-.44)
\psline(-.3,-.36)(-.5,-.2)(-.5,.2)(-.3,.36)
\psline(-.2,.44)(0,.6)
\pcline(0,1)(0,.6)
\psline(.5,-1)(.5,-.2)(.3,-.04)
\psline(.2,.04)(-.5,.6)
\pcline(-.5,.6)(-.5,1)
\end{pspicture}
$$
\caption{$(i)$ Reidemeister move I,
$(ii)$ Reidemeister move II and
$(iii)$ Reidemeister move III.} \label{Reide}
\end{figure}

A \emph{link} $L$ is an embedding of $n$ copies of $\mathbb{S}^1$ in
$\mathbb{S}^3$. The number of components of $L$ is denoted by $\ell(L)$.
In the case $\ell(L)=1$, a link is called a \emph{knot}.
Throughout the article, we will assume all links are \emph{tame}
which means all links can be in a form of a finite union of line segments,
in the language in graph theory, a knot can be considered as an
embedding of the cycle graph $C_k$ of $k$ vertices into $\mathbb{S}^3$,
which is called a \emph{spatial graph} of $C_k$. Two links are \emph{equivalent}
if there is an isotopy between them.
In the case of prime knots, this equivalence is the same as
the existence of an orientation preserving homeomorphism on $\mathbb{S}^3$,
which sends a knot to the other knot. Although the equivalent
class of a link $L$ is called a \emph{link type},
throughout the article, a link really means the equivalent class of link $L$.

A useful way to visualise and manipulate knots is to project the knot onto
a plane. A small change in the direction of projection will ensure that
it is one-to-one except at the double points, called \emph{crossings},
where the shadow of the knot crosses itself once transversely~\cite{Kauffman}.
At each crossing, to be able to recreate the original knot, the over-strand
must be distinguished from the under-strand. This is often done by creating
a break in the strand going underneath. The resulting diagram is an immersed
plane curve with the additional data of which strand is over and which is
under at each crossing. These diagrams are called \emph{knot diagrams} when
they represent a knot and \emph{link diagrams} when they represent a link.
Reidemeister proved that two link diagrams belonging to the same link type
can be related a finite sequence of planar isotopies and three Reidemeister moves
in Figure~\ref{Reide}~\cite{Reidemeister}.

\begin{figure}
$$
\begin{pspicture}[shift=-2](-.2,-.7)(4.2,3.7)
\psline(0,.5)(0,3) \psline[arrowscale=1.5]{->}(0,1.74)(0,1.76)
\pccurve[angleA=90,angleB=180](0,3)(.5,3.5)
\psline(.5,3.5)(1,3.5)
\pccurve[angleA=0,angleB=135](1,3.5)(1.2,3.3)
\pccurve[angleA=-45,angleB=180](1.3,3.2)(1.5,3)
\psline(1.5,3)(2,3) \psline[arrowscale=1.5]{->}(1.74,3)(1.76,3)
\pccurve[angleA=0,angleB=-135](2,3)(2.25,3.25)
\pccurve[angleA=45,angleB=180](2.25,3.25)(2.5,3.5)
\psline(2.5,3.5)(3.5,3.5) \psline[arrowscale=1.5]{->}(2.99,3.5)(3.01,3.5)
\pccurve[angleA=0,angleB=90](3.5,3.5)(4,3)
\psline(4,3)(4,.5) \psline[arrowscale=1.5]{->}(4,1.76)(4,1.74)
\pccurve[angleA=-90,angleB=0](4,.5)(3.75,0)
\pccurve[angleA=-90,angleB=180](3.5,.5)(3.75,0)
\pccurve[angleA=90,angleB=-45](3.5,.5)(3.28,.72)
\pccurve[angleA=135,angleB=-90](3.22,.78)(3,1)
\psline(3,1)(3,1.5) \psline[arrowscale=1.5]{->}(3,1.24)(3,1.26)
\pccurve[angleA=90,angleB=-135](3,1.5)(3.25,1.75)
\pccurve[angleA=45,angleB=-90](3.25,1.75)(3.5,2)
\psline(3.5,2)(3.5,2.75) \psline[arrowscale=1.5]{->}(3.5,2.37)(3.5,2.38)
\pccurve[angleA=90,angleB=0](3.5,2.75)(3.25,3)
\psline(3.25,3)(3,3)
\pccurve[angleA=180,angleB=45](3,3)(2.78,2.78)
\pccurve[angleA=-135,angleB=0](2.72,2.72)(2.5,2.5)
\pccurve[angleA=180,angleB=90](2.5,2.5)(2,2)
\pccurve[angleA=-90,angleB=45](2,2)(1.75,1.75)
\pccurve[angleA=-135,angleB=90](1.75,1.75)(1.5,1.5)
\psline(1.5,1.5)(1.5,1) \psline[arrowscale=1.5]{->}(1.5,1.26)(1.5,1.24)
\pccurve[angleA=-90,angleB=135](1.5,1)(1.72,.78)
\pccurve[angleA=-45,angleB=90](1.78,.72)(2,.5)
\psline(2,.5)(2,.25)
\pccurve[angleA=-90,angleB=180](2,.25)(2.25,0)
\psline(2.25,0)(2.75,0) \psline[arrowscale=1.5]{->}(2.49,0)(2.51,0)
\pccurve[angleA=0,angleB=-90](2.75,0)(3,.25)
\psline(3,.25)(3,.5) \psline[arrowscale=1.5]{->}(3,.37)(3,.38)
\pccurve[angleA=90,angleB=-135](3,.5)(3.25,.75)
\pccurve[angleA=45,angleB=-90](3.25,.75)(3.5,1)
\psline(3.5,1)(3.5,1.5) \psline[arrowscale=1.5]{->}(3.5,1.24)(3.5,1.26)
\pccurve[angleA=90,angleB=-45](3.5,1.5)(3.28,1.71)
\pccurve[angleA=135,angleB=-90](3.22,1.78)(3,2)
\psline(3,2)(3,2.5) \psline[arrowscale=1.5]{->}(3,2.24)(3,2.26)
\pccurve[angleA=90,angleB=-45](3,2.5)(2.75,2.75)
\pccurve[angleA=135,angleB=-45](2.75,2.75)(2.5,3)
\pccurve[angleA=135,angleB=-45](2.5,3)(2.3,3.2)
\pccurve[angleA=135,angleB=0](2.2,3.3)(2,3.5)
\psline(2,3.5)(1.5,3.5) \psline[arrowscale=1.5]{->}(1.76,3.5)(1.74,3.5)
\pccurve[angleA=180,angleB=45](1.5,3.5)(1.25,3.25)
\pccurve[angleA=-135,angleB=0](1.25,3.25)(1,3)
\pccurve[angleA=180,angleB=90](1,3)(.5,2.75)
\psline(.5,2.75)(.5,2.5)
\pccurve[angleA=-90,angleB=135](.5,2.5)(1.5,2)
\pccurve[angleA=-45,angleB=135](1.5,2)(1.72,1.78)
\pccurve[angleA=-45,angleB=90](1.78,1.72)(2,1.5)
\psline(2,1.5)(2,1) \psline[arrowscale=1.5]{->}(2,1.26)(2,1.24)
\pccurve[angleA=-90,angleB=45](2,1)(1.75,.75)
\pccurve[angleA=-135,angleB=90](1.75,.75)(1.5,.5)
\psline(1.5,.5)(1.5,.25) \psline[arrowscale=1.5]{->}(1.5,.38)(1.5,.37)
\pccurve[angleA=-90,angleB=0](1.5,.25)(1.25,0)
\psline(1.25,0)(.5,0) \psline[arrowscale=1.5]{->}(.88,0)(.87,0)
\pccurve[angleA=180,angleB=-90](.5,0)(0,.5)
\rput[t](2,-.5){$(i)$} \rput[t](-.4,2){$7_5$}
\end{pspicture}
\quad\quad  \begin{pspicture}[shift=-2](-.2,-.7)(4.2,3.7)
\psline(0,.5)(1,3)\psline[arrowscale=1.5]{->}(1,3)(1.1,3.25)
\psline(1,.5)(.53,1.68)
\psline(.47,1.82)(0,3)\psline[arrowscale=1.5]{->}(0,3)(-.1,3.25)
\psline(4,.5)(3,3)\psline[arrowscale=1.5]{->}(3,3)(2.9,3.25)
\psline(3,.5)(3.47,1.68)
\psline(3.53,1.82)(4,3)\psline[arrowscale=1.5]{->}(4,3)(4.1,3.25)
\rput(1.2,1.75){$\Rightarrow$}\rput(2.8,1.75){$\Leftarrow$}
\pccurve[angleA=-80,angleB=90](1.5,3.25)(1.85,1.75)
\pccurve[angleA=-90,angleB=80](1.85,1.75)(1.5,.5)
\pccurve[angleA=-100,angleB=90](2.5,3.25)(2.15,1.75)
\pccurve[angleA=-90,angleB=100](2.15,1.75)(2.5,.5)
\psline[arrowscale=1.5]{->}(1.85,1.7699)(1.85,1.77)
\psline[arrowscale=1.5]{->}(2.15,1.7699)(2.15,1.77)
\rput[t](2,-.5){$(ii)$}
\end{pspicture}
$$

$$\begin{pspicture}[shift=-2](-.2,-.7)(4.2,3.7)
\psline(0,.5)(0,3) \psline[arrowscale=1.5]{->}(0,1.74)(0,1.76)
\pccurve[angleA=90,angleB=180](0,3)(.5,3.5)
\psline(.5,3.5)(1,3.5)
\pccurve[angleA=0,angleB=0](1,3.5)(1,3)
\pccurve[angleA=180,angleB=180](1.5,3.5)(1.5,3)
\psline(1.5,3)(2,3) \psline[arrowscale=1.5]{->}(1.74,3)(1.76,3)
\pccurve[angleA=0,angleB=0](2,3)(2,3.5)
\pccurve[angleA=180,angleB=180](2.5,3.5)(2.5,3.5)
\psline(2.5,3.5)(3.5,3.5) \psline[arrowscale=1.5]{->}(2.99,3.5)(3.01,3.5)
\pccurve[angleA=0,angleB=90](3.5,3.5)(4,3)
\psline(4,3)(4,.5) \psline[arrowscale=1.5]{->}(4,1.76)(4,1.74)
\pccurve[angleA=-90,angleB=0](4,.5)(3.75,0)
\pccurve[angleA=-90,angleB=180](3.5,.5)(3.75,0)
\pccurve[angleA=90,angleB=-90](3.5,.5)(3.35,.75)
\pccurve[angleA=90,angleB=-90](3.35,.75)(3.5,1)
\psline(3,1)(3,1.5) \psline[arrowscale=1.5]{->}(3,1.24)(3,1.26)
\pccurve[angleA=90,angleB=-90](3,1.5)(3.15,1.75)
\pccurve[angleA=90,angleB=-90](3.15,1.75)(3,2)
\psline(3.5,2)(3.5,2.75) \psline[arrowscale=1.5]{->}(3.5,2.37)(3.5,2.38)
\pccurve[angleA=90,angleB=0](3.5,2.75)(3.25,3)
\psline(3.25,3)(3,3)
\pccurve[angleA=180,angleB=0](3,3)(2.75,2.85)
\pccurve[angleA=180,angleB=0](2.75,2.85)(2.5,3)
\pccurve[angleA=180,angleB=90](2.5,2.5)(2,2)
\pccurve[angleA=-90,angleB=90](2,2)(1.85,1.75)
\pccurve[angleA=-90,angleB=90](1.85,1.75)(2,1.5)
\psline(1.5,1.5)(1.5,1) \psline[arrowscale=1.5]{->}(1.5,1.26)(1.5,1.24)
\pccurve[angleA=-90,angleB=90](1.5,1)(1.65,.75)
\pccurve[angleA=-90,angleB=90](1.65,.75)(1.5,.5)
\psline(2,.5)(2,.25)
\pccurve[angleA=-90,angleB=180](2,.25)(2.25,0)
\psline(2.25,0)(2.75,0) \psline[arrowscale=1.5]{->}(2.49,0)(2.51,0)
\pccurve[angleA=0,angleB=-90](2.75,0)(3,.25)
\psline(3,.25)(3,.5) \psline[arrowscale=1.5]{->}(3,.37)(3,.38)
\pccurve[angleA=90,angleB=-90](3,.5)(3.15,.75)
\pccurve[angleA=90,angleB=-90](3.15,.75)(3,1)
\psline(3.5,1)(3.5,1.5) \psline[arrowscale=1.5]{->}(3.5,1.24)(3.5,1.26)
\pccurve[angleA=90,angleB=-90](3.5,1.5)(3.35,1.75)
\pccurve[angleA=90,angleB=-90](3.35,1.75)(3.5,2)
\psline(3,2)(3,2.5) \psline[arrowscale=1.5]{->}(3,2.24)(3,2.26)
\pccurve[angleA=90,angleB=0](3,2.5)(2.75,2.65)
\pccurve[angleA=180,angleB=0](2.75,2.65)(2.5,2.5)
\pccurve[angleA=180,angleB=180](2.5,3)(2.5,3.5)
\psline(2,3.5)(1.5,3.5) \psline[arrowscale=1.5]{->}(1.76,3.5)(1.74,3.5)
\pccurve[angleA=180,angleB=90](1,3)(.5,2.75)
\psline(.5,2.75)(.5,2.5)
\pccurve[angleA=-90,angleB=90](.5,2.5)(1.5,2)
\pccurve[angleA=-90,angleB=90](1.5,2)(1.65,1.75)
\pccurve[angleA=-90,angleB=90](1.65,1.75)(1.5,1.5)
\psline(2,1.5)(2,1) \psline[arrowscale=1.5]{->}(2,1.26)(2,1.24)
\pccurve[angleA=-90,angleB=90](2,1)(1.85,.75)
\pccurve[angleA=-90,angleB=90](1.85,.75)(2,.5)
\psline(1.5,.5)(1.5,.25) \psline[arrowscale=1.5]{->}(1.5,.38)(1.5,.37)
\pccurve[angleA=-90,angleB=0](1.5,.25)(1.25,0)
\psline(1.25,0)(.5,0) \psline[arrowscale=1.5]{->}(.88,0)(.87,0)
\pccurve[angleA=180,angleB=-90](.5,0)(0,.5)
\rput[t](2,-.5){$(iii)$}
\end{pspicture}
\quad\quad  \begin{pspicture}[shift=-2](-.2,-.7)(4.2,3.7)
\psline(0,.5)(0,3) \psline[arrowscale=1.5]{->}(0,1.74)(0,1.76)
\pccurve[angleA=90,angleB=180](0,3)(.5,3.5)
\psline(.5,3.5)(1,3.5)
\pccurve[angleA=0,angleB=135](1,3.5)(1.2,3.3)
\pccurve[angleA=-45,angleB=180](1.3,3.2)(1.5,3)
\psline(1.5,3)(2,3) \psline[arrowscale=1.5]{->}(1.74,3)(1.76,3)
\pccurve[angleA=0,angleB=-135](2,3)(2.25,3.25)
\pccurve[angleA=45,angleB=180](2.25,3.25)(2.5,3.5)
\psline(2.5,3.5)(3.5,3.5) \psline[arrowscale=1.5]{->}(2.99,3.5)(3.01,3.5)
\pccurve[angleA=0,angleB=90](3.5,3.5)(4,3)
\psline(4,3)(4,.5) \psline[arrowscale=1.5]{->}(4,1.76)(4,1.74)
\pccurve[angleA=-90,angleB=0](4,.5)(3.75,0)
\pccurve[angleA=-90,angleB=180](3.5,.5)(3.75,0)
\pccurve[angleA=90,angleB=-45](3.5,.5)(3.28,.72)
\pccurve[angleA=135,angleB=-90](3.22,.78)(3,1)
\psline(3,1)(3,1.5) \psline[arrowscale=1.5]{->}(3,1.24)(3,1.26)
\pccurve[angleA=90,angleB=-135](3,1.5)(3.25,1.75)
\pccurve[angleA=45,angleB=-90](3.25,1.75)(3.5,2)
\psline(3.5,2)(3.5,2.75) \psline[arrowscale=1.5]{->}(3.5,2.37)(3.5,2.38)
\pccurve[angleA=90,angleB=0](3.5,2.75)(3.25,3)
\psline(3.25,3)(3,3)
\pccurve[angleA=180,angleB=45](3,3)(2.78,2.78)
\pccurve[angleA=-135,angleB=0](2.72,2.72)(2.5,2.5)
\pccurve[angleA=180,angleB=90](2.5,2.5)(2,2)
\pccurve[angleA=-90,angleB=45](2,2)(1.75,1.75)
\pccurve[angleA=-135,angleB=90](1.75,1.75)(1.5,1.5)
\psline(1.5,1.5)(1.5,1) \psline[arrowscale=1.5]{->}(1.5,1.26)(1.5,1.24)
\pccurve[angleA=-90,angleB=135](1.5,1)(1.72,.78)
\pccurve[angleA=-45,angleB=90](1.78,.72)(2,.5)
\psline(2,.5)(2,.25)
\pccurve[angleA=-90,angleB=180](2,.25)(2.25,0)
\psline(2.25,0)(2.75,0) \psline[arrowscale=1.5]{->}(2.49,0)(2.51,0)
\pccurve[angleA=0,angleB=-90](2.75,0)(3,.25)
\psline(3,.25)(3,.5) \psline[arrowscale=1.5]{->}(3,.37)(3,.38)
\pccurve[angleA=90,angleB=-135](3,.5)(3.25,.75)
\pccurve[angleA=45,angleB=-90](3.25,.75)(3.5,1)
\psline(3.5,1)(3.5,1.5) \psline[arrowscale=1.5]{->}(3.5,1.24)(3.5,1.26)
\pccurve[angleA=90,angleB=-45](3.5,1.5)(3.28,1.71)
\pccurve[angleA=135,angleB=-90](3.22,1.78)(3,2)
\psline(3,2)(3,2.5) \psline[arrowscale=1.5]{->}(3,2.24)(3,2.26)
\pccurve[angleA=90,angleB=-45](3,2.5)(2.75,2.75)
\pccurve[angleA=135,angleB=-45](2.75,2.75)(2.5,3)
\pccurve[angleA=135,angleB=-45](2.5,3)(2.3,3.2)
\pccurve[angleA=135,angleB=0](2.2,3.3)(2,3.5)
\psline(2,3.5)(1.5,3.5) \psline[arrowscale=1.5]{->}(1.76,3.5)(1.74,3.5)
\pccurve[angleA=180,angleB=45](1.5,3.5)(1.25,3.25)
\pccurve[angleA=-135,angleB=0](1.25,3.25)(1,3)
\pccurve[angleA=180,angleB=90](1,3)(.5,2.75)
\psline(.5,2.75)(.5,2.5)
\pccurve[angleA=-90,angleB=135](.5,2.5)(1.5,2)
\pccurve[angleA=-45,angleB=135](1.5,2)(1.72,1.78)
\pccurve[angleA=-45,angleB=90](1.78,1.72)(2,1.5)
\psline(2,1.5)(2,1) \psline[arrowscale=1.5]{->}(2,1.26)(2,1.24)
\pccurve[angleA=-90,angleB=45](2,1)(1.75,.75)
\pccurve[angleA=-135,angleB=90](1.75,.75)(1.5,.5)
\psline(1.5,.5)(1.5,.25) \psline[arrowscale=1.5]{->}(1.5,.38)(1.5,.37)
\pccurve[angleA=-90,angleB=0](1.5,.25)(1.25,0)
\psline(1.25,0)(.5,0) \psline[arrowscale=1.5]{->}(.88,0)(.87,0)
\pccurve[angleA=180,angleB=-90](.5,0)(0,.5)
\rput(.75,.75){$a$}\rput(2.5,1.25){$b$}
\rput(3.75,1.75){$c$}\rput(1.75,3.25){$d$}
\rput[t](2,-.5){$(iv)$} \rput[t](.75,1.75){$S(7_5)$}
\end{pspicture}
$$
\caption{The Seifert algorithm to produce a canonical Seifert surface of the knot $7_5$.} \label{SeifertS}
\end{figure}
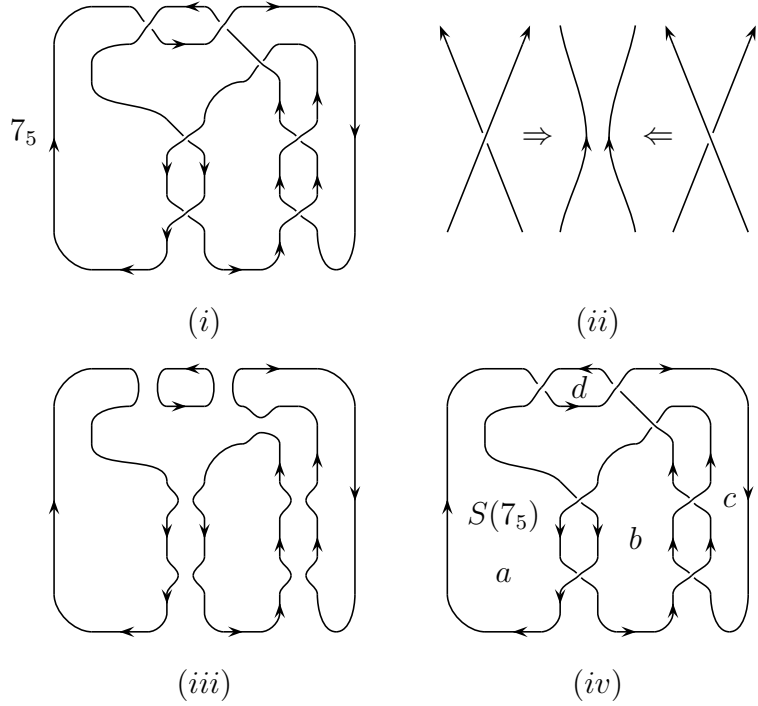

A compact orientable surface
$S$ is a \emph{Seifert surface} of $L$ if the boundary of $S$ is isotopic to
$L$. The existence of such a surface was first proven by Seifert
using an algorithm on a diagram of $L$: first, we oriented each components of the link,
resolving each crossings by the rule illustrated in Figure~\ref{SeifertS} $(ii)$,
the resulting simple closed curves are called \emph{Seifert circles}
and a Seifert surface
is obtained by connecting discs, bounded by Seifert circles,
by half twisted bands as same as the original crossings
as illustrated in Figure~\ref{SeifertS} $(i), (iii)$ and $(iv)$.
This algorithm was named after him as \emph{Seifert's algorithm}~\cite{Seifert:def}.
A Seifert surface of a link $L$ obtained by applying Seifert's
algorithm for a diagram of $L$ is called
a \emph{canonical Seifert surface}, denoted by $S(L)$. However, not all Seifert
surface is canonical~\cite{Brittenham:free}.

Some Seifert surfaces feature extra structures. Seifert
surfaces obtained by annuli plumbings are the main subjects of this
article. Even though higher dimensional plumbings can be defined
here we will only concentrate on \emph{annuli plumbings}.
This is often called a \emph{Murasugi sum} and it
has been studied extensively for the fibreness of links and surfaces
\cite{Gabai:genera, Stallings:const}. To show the existence of these
plumbing surfaces of a link, it is common to present the link as the closure of
a braid in a classical Artin group \cite{FHK:openbook, HW:plumbing}.
Furthermore, a few different ways to find braid presentations of a link
have been found by Alexander~\cite{Birman:Braids}, Morton
\cite{Morton:Threading}, Vogel~\cite{vogel:cmh} and
Yamada~\cite{yamada}. In particular, the work of Yamada is closely
related with Seifert's algorithm and has been generalized to find
another beautiful presentation of the braid groups~\cite{BKL}. Several authors have shown
the existence of basket surfaces, flat plumbing surfaces and flat plumbing basket
surfaces using a braid presentative of a link, where the link is presented by
a closure of a braid $\sigma_1 \sigma_2 \ldots \sigma_n W(\sigma_1, \sigma_2, \ldots ,\sigma_n)$
 and the disc is chosen by the union of Seifert discs connected by
 half twist bands corresponding to the first $n$ braid word
 $\sigma_1 \sigma_2 \ldots \sigma_n$~\cite{FHK:openbook, HW:plumbing}.

\subsection{Seifert graphs of Seifert surfaces} \label{preSeifert}

A canonical Seifert surface $S$ gives rise to a natural signed graph,
which is called the \emph{Seifert graph $\Gamma(S)$}
by shrinking each disc to a point and at the same time the width of
the half twisted band is shrunk to a signed edge as illustrated in
Figure~\ref{graph} $(ii)$~\cite{Murasugi:knot}.
In the process, we kept the cyclic order of the adjacent edges at each vertex,
thus, one may consider Seifert graphs as maps without the direction of edges.
These processes can also be performed on arbitrary Seifert surfaces.
A simple graph
obtained from $\Gamma$ by identifying edges in the same parallel class is called
\emph{the induced graph} of  $\Gamma$, denoted by $G(\Gamma)$
as illustrated in Figure~\ref{Figure8} $(ii)$.
Since a link $L$ is tame and its Seifert surface $S(L)$ is compact,
the Seifert graph $\Gamma(S(L))$ is finite.
By separating discs by local orientation as indicated on each
vertices in Figure~\ref{graph} $(ii)$, it can be considered as a bipartite graph.
Although, the example in Figure~\ref{graph} $(ii)$ is planar,
Seifert graph is not planar in general as illustrated in Figure~\ref{Figure8} $(ii)$.
If the Seifert surface is connected, then its Seifert graph is also connected.

\begin{figure}
$$
\begin{pspicture}[shift=-2](-.2,-.7)(4.2,3.7)
\psline(0,.5)(0,3) \psline[arrowscale=1.5]{->}(0,1.74)(0,1.76)
\pccurve[angleA=90,angleB=180](0,3)(.5,3.5)
\psline(.5,3.5)(1,3.5)
\pccurve[angleA=0,angleB=135](1,3.5)(1.2,3.3)
\pccurve[angleA=-45,angleB=180](1.3,3.2)(1.5,3)
\psline(1.5,3)(2,3) \psline[arrowscale=1.5]{->}(1.74,3)(1.76,3)
\pccurve[angleA=0,angleB=-135](2,3)(2.25,3.25)
\pccurve[angleA=45,angleB=180](2.25,3.25)(2.5,3.5)
\psline(2.5,3.5)(3.5,3.5) \psline[arrowscale=1.5]{->}(2.99,3.5)(3.01,3.5)
\pccurve[angleA=0,angleB=90](3.5,3.5)(4,3)
\psline(4,3)(4,.5) \psline[arrowscale=1.5]{->}(4,1.76)(4,1.74)
\pccurve[angleA=-90,angleB=0](4,.5)(3.75,0)
\pccurve[angleA=-90,angleB=180](3.5,.5)(3.75,0)
\pccurve[angleA=90,angleB=-45](3.5,.5)(3.28,.72)
\pccurve[angleA=135,angleB=-90](3.22,.78)(3,1)
\psline(3,1)(3,1.5) \psline[arrowscale=1.5]{->}(3,1.24)(3,1.26)
\pccurve[angleA=90,angleB=-135](3,1.5)(3.25,1.75)
\pccurve[angleA=45,angleB=-90](3.25,1.75)(3.5,2)
\psline(3.5,2)(3.5,2.75) \psline[arrowscale=1.5]{->}(3.5,2.37)(3.5,2.38)
\pccurve[angleA=90,angleB=0](3.5,2.75)(3.25,3)
\psline(3.25,3)(3,3)
\pccurve[angleA=180,angleB=45](3,3)(2.78,2.78)
\pccurve[angleA=-135,angleB=0](2.72,2.72)(2.5,2.5)
\pccurve[angleA=180,angleB=90](2.5,2.5)(2,2)
\pccurve[angleA=-90,angleB=45](2,2)(1.75,1.75)
\pccurve[angleA=-135,angleB=90](1.75,1.75)(1.5,1.5)
\psline(1.5,1.5)(1.5,1) \psline[arrowscale=1.5]{->}(1.5,1.26)(1.5,1.24)
\pccurve[angleA=-90,angleB=135](1.5,1)(1.72,.78)
\pccurve[angleA=-45,angleB=90](1.78,.72)(2,.5)
\psline(2,.5)(2,.25)
\pccurve[angleA=-90,angleB=180](2,.25)(2.25,0)
\psline(2.25,0)(2.75,0) \psline[arrowscale=1.5]{->}(2.49,0)(2.51,0)
\pccurve[angleA=0,angleB=-90](2.75,0)(3,.25)
\psline(3,.25)(3,.5) \psline[arrowscale=1.5]{->}(3,.37)(3,.38)
\pccurve[angleA=90,angleB=-135](3,.5)(3.25,.75)
\pccurve[angleA=45,angleB=-90](3.25,.75)(3.5,1)
\psline(3.5,1)(3.5,1.5) \psline[arrowscale=1.5]{->}(3.5,1.24)(3.5,1.26)
\pccurve[angleA=90,angleB=-45](3.5,1.5)(3.28,1.71)
\pccurve[angleA=135,angleB=-90](3.22,1.78)(3,2)
\psline(3,2)(3,2.5) \psline[arrowscale=1.5]{->}(3,2.24)(3,2.26)
\pccurve[angleA=90,angleB=-45](3,2.5)(2.75,2.75)
\pccurve[angleA=135,angleB=-45](2.75,2.75)(2.5,3)
\pccurve[angleA=135,angleB=-45](2.5,3)(2.3,3.2)
\pccurve[angleA=135,angleB=0](2.2,3.3)(2,3.5)
\psline(2,3.5)(1.5,3.5) \psline[arrowscale=1.5]{->}(1.76,3.5)(1.74,3.5)
\pccurve[angleA=180,angleB=45](1.5,3.5)(1.25,3.25)
\pccurve[angleA=-135,angleB=0](1.25,3.25)(1,3)
\pccurve[angleA=180,angleB=90](1,3)(.5,2.75)
\psline(.5,2.75)(.5,2.5)
\pccurve[angleA=-90,angleB=135](.5,2.5)(1.5,2)
\pccurve[angleA=-45,angleB=135](1.5,2)(1.72,1.78)
\pccurve[angleA=-45,angleB=90](1.78,1.72)(2,1.5)
\psline(2,1.5)(2,1) \psline[arrowscale=1.5]{->}(2,1.26)(2,1.24)
\pccurve[angleA=-90,angleB=45](2,1)(1.75,.75)
\pccurve[angleA=-135,angleB=90](1.75,.75)(1.5,.5)
\psline(1.5,.5)(1.5,.25) \psline[arrowscale=1.5]{->}(1.5,.38)(1.5,.37)
\pccurve[angleA=-90,angleB=0](1.5,.25)(1.25,0)
\psline(1.25,0)(.5,0) \psline[arrowscale=1.5]{->}(.88,0)(.87,0)
\pccurve[angleA=180,angleB=-90](.5,0)(0,.5)
\rput(.75,1.25){$a$}\rput(2.5,1.25){$b$}
\rput(3.75,1.75){$c$}\rput(1.75,3.25){$d$}
\rput[t](2,-.5){$(i)$} \rput[t](-.8,2){$S(7_5)$}
\end{pspicture}
\quad\quad  \begin{pspicture}[shift=-2](-.3,-.7)(2.3,3.7) \psline(0,1)(2,1)
\psline(0,1)(1,3) \psline(1,3)(2,1)
\psarc(.5,1){.5}{180}{0}
\psarc(1.5,1){.5}{180}{0}
\pccurve[angleA=70,angleB=110](1,1)(2,1)
\pscircle[fillstyle=solid,fillcolor=darkgray](0,1){.25}
\pscircle[fillstyle=solid,fillcolor=darkgray](1,1){.25}
\pscircle[fillstyle=solid,fillcolor=darkgray](2,1){.25}
\pscircle[fillstyle=solid,fillcolor=darkgray](1,3){.25}
\rput(0,1){$+$} \rput(1,1){$-$} \rput(2,1){$+$} \rput(1,3){$-$}
\rput[t](.5,1.35){$+$}\rput[t](.5,.45){$+$}
\rput[t](1.5,.45){$+$}\rput[t](1.5,.9){$+$}\rput[t](1.5,1.55){$+$}
\rput[t](.4,2.5){$+$}\rput[t](1.6,2.5){$+$}
\rput[t](-.3,1.45){$a$}\rput[t](1,1.65){$b$}
\rput[t](2.2,1.45){$c$}\rput[t](1,3.65){$d$}
\rput[t](1,-.5){$(ii)$} \rput[t](0,3.5){$\Gamma(S(7_5))$}
\end{pspicture}
\quad\quad  \begin{pspicture}[shift=-2](-.3,-.7)(2.3,3.7) \psline(0,1)(1,1)
\psline(0,1)(1,3) \psline(1,3)(2,1)
\pscircle[fillstyle=solid,fillcolor=darkgray](0,1){.1}
\pscircle[fillstyle=solid,fillcolor=darkgray](1,1){.1}
\pscircle[fillstyle=solid,fillcolor=darkgray](2,1){.1}
\pscircle[fillstyle=solid,fillcolor=darkgray](1,3){.1}
\rput[t](.5,1.4){$+$}
\rput[t](.4,2.5){$+$}\rput[t](1.6,2.5){$+$}
\rput[t](-.3,1.25){$a$}\rput[t](1,1.45){$b$}
\rput[t](2.2,1.25){$c$}\rput[t](1,3.45){$d$}
\rput[t](1,-.5){$(iii)$} \rput[t](0,3.5){$T$}
\end{pspicture}
$$
\caption{$(i)$ A Seifert surface $S(7_5)$, $(ii)$ its corresponding signed Seifert graph $\Gamma(S(7_5))$ and
$(iii)$ a spanning tree $T$ of $\Gamma(S(7_5))$.} \label{graph}
\end{figure}
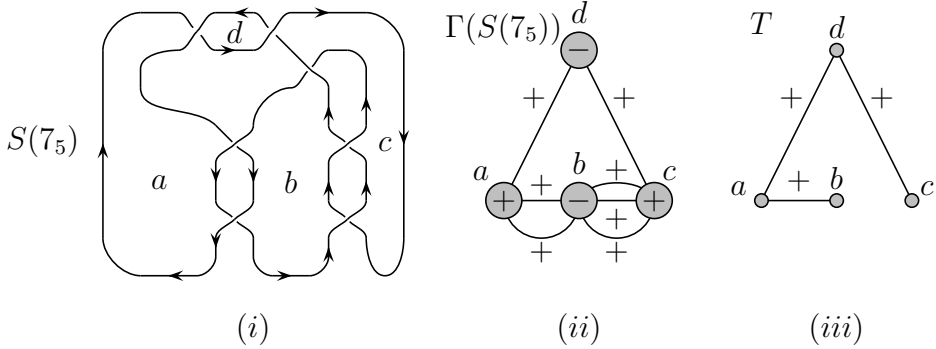

\begin{figure}
$$
\begin{pspicture}[shift=.2](-2.2,-2.7)(2.2,2.2)
\psarc(0,0){2}{15}{75}
\pccurve[angleA=165,angleB=20](2;75)(1.5;110)
\psarc(0,0){1.5}{110}{160}
\pccurve[angleA=-110,angleB=135](1.5;160)(1.3;177)
\pccurve[angleA=-45,angleB=120](1.2;183)(1;210)
\psarc(0,0){1}{210}{330}
\pccurve[angleA=60,angleB=-70](1;330)(1.5;20)
\psarc(0,0){1.5}{20}{70}
\pccurve[angleA=160,angleB=-45](1.5;70)(1.7;87)
\pccurve[angleA=135,angleB=15](1.8;93)(2;105)
\psarc(0,0){2}{105}{255}
\pccurve[angleA=-15,angleB=-160](2;255)(1.5;290)
\psarc(0,0){1.5}{290}{-20}
\pccurve[angleA=70,angleB=-45](1.5;-20)(1.3;-3)
\pccurve[angleA=135,angleB=-60](1.2;3)(1;30)
\psarc(0,0){1}{30}{150}
\pccurve[angleA=-120,angleB=110](1;150)(1.5;200)
\psarc(0,0){1.5}{200}{250}
\pccurve[angleA=-15,angleB=135](1.5;250)(1.65;267)
\pccurve[angleA=-45,angleB=195](1.8;273)(2;285)
\psarc(0,0){2}{285}{20}
\rput(0,-2.5){$(i)$} \rput(0,0){$S(4_1)$}
\end{pspicture} \quad\quad
\begin{pspicture}[shift=.2](-1.2,-2.7)(1.8,2.2)
\pccurve[angleA=0,angleB=0](0,0)(0,1.5)
\pccurve[angleA=180,angleB=180](0,0)(0,1.5)
\pccurve[angleA=-90,angleB=90](0,0)(0,-1.5)
\pccurve[angleA=90,angleB=0](0,0)(-.2,.2)
\pccurve[angleA=180,angleB=180](-.35,.2)(-.35,-1.7)
\pccurve[angleA=0,angleB=-90](-.35,-1.7)(0,-1.5)
\pscircle[fillstyle=solid,fillcolor=darkgray](0,1.5){.1}
\pscircle[fillstyle=solid,fillcolor=darkgray](0,0){.1}
\pscircle[fillstyle=solid,fillcolor=darkgray](0,-1.5){.1}
\rput(0,-2.5){$(ii)$} \rput(1,-.5){$\Gamma(S(4_1))$}
\end{pspicture}
\quad\quad
\begin{pspicture}[shift=.2](-.2,-2.7)(2,2.2)
\psline(0,1.5)(0,-1.5)
\pscircle[fillstyle=solid,fillcolor=darkgray](0,1.5){.1}
\pscircle[fillstyle=solid,fillcolor=darkgray](0,0){.1}
\pscircle[fillstyle=solid,fillcolor=darkgray](0,-1.5){.1}
\rput(0,-2.5){$(ii)$} \rput(1.5,0){$G(\Gamma(S(4_1)))$}
\end{pspicture}
$$
\caption{$(i)$ A canonical Seifert surface
obtained from a closed braid diagram of the figure eight
knot $4_1$, $(ii)$ its Seifert graph $\Gamma(S(4_1))$ and $(iii)$
its induced graph $G(\Gamma(S(4_1)))$.}
\label{Figure8}
\end{figure}
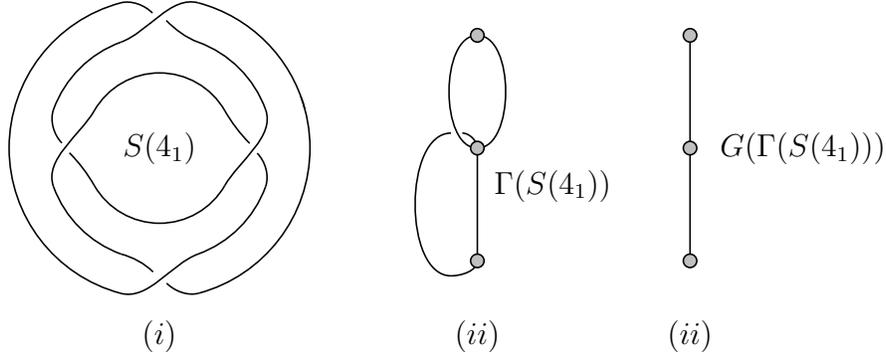

\begin{notation}\label{seifertgraphdef}
Let $S(L)$ be a canonical Seifert surface of a link $L$.
The number of Seifert circles in $S(L)$ is denoted by $s(S(L))$ and
the number of the half twisted bands in $S(L)$ is denoted by $c(S(L))$.
\end{notation}

It is fairly easy to see that $s(S(L))=|V(\Gamma(S(L)))|$, the
cardinality of the vertex set, and $c(S(L))= |E(\Gamma(S(L)))|$,
the cardinality of the edge set. If the surface $S$ has a genus $g$ and
the link $L=\partial S$ has $\ell(L)$ components, by Euler characteristics formula, we have
$$s(S(L)) - c(S(L)) + \ell(L) =2- 2g$$
because $\ell(L)$ is the number of faces in the $2$-cell embedding of
$\Gamma(S)$ into the surface $S$ of genus $g$.

A spanning tree $T$ of $\Gamma(S(L))$ is depicted in Figure~\ref{graph} $(iii)$.
Recall that the number of edges of a spanning tree of a connected
graph with $n$ vertices is $n-1$.
Since it is bipartite by the orientation of the surface $S(L)$ and all edge $e$'s in $E(\Gamma(S(L)))$
are half twisted bands in $S(L)$, one can see that the length of the path in
the tree $T$ joining both end vertices $u$ and $v$ of any co-tree edge $e=\{u,v\} \in E(\Gamma(S(L)))\setminus E(T)$ is odd.

Spanning trees of the Seifert graph play a key role in the research of plumbing surfaces.
If we use a braid presentation of a link, its Seifert graph of the canonical Seifert surface is
a path with multi-edges. Thus, there is no ambiguity for the choice
of a spanning tree $T$ and the existence of a desired labeling on $T$.
Let us provide a proper definition for labelings as follows.

\begin{defi}
A graph $\Gamma$ is \emph{signed} if there is a function $\phi : E(\Gamma) \rightarrow \{ + , - \}$.
A \emph{labeling} on $\Gamma$ means to be an edge $2$-coloring which
is a function $\mu : E(\Gamma) \rightarrow \{+, - \}$ unless stated differently.
For a spanning tree $T$ of a connected graph $\Gamma$, a labeling
$\mu: E(T) \rightarrow \{ + , - \}$ is
\emph{alternating} on a spanning tree $T$ of a graph $\Gamma$ if
for any co-tree edge $e \in E(\Gamma) \setminus E(T)$, the unique path
$v_0, v_1, v_2, \ldots, v_{\ell}$ in $T$ joining both end vertices
of the edge $e$ satisfies that $\mu(v_i v_{i+1}) \neq \mu(v_{i+1}
v_{i+2})$ for any $i=0, 1, \ldots, \ell -2$. A spanning tree $T$
is called \emph{co-tree edge alternating spanning tree} if there exists an alternating
labeling on $T$.
\end{defi}

\begin{figure}
$$
\begin{pspicture}[shift=-1.2](-.3,-1.4)(.3,.8)
\psline(-.2,-.8)(-.2,.8) \psline[arrowscale=1.5]{->}(-.2,-.1)(-.2,.1)
\psline(.2,-.8)(.2,.8) \psline[arrowscale=1.5]{->}(.2,-.1)(.2,.1)
\end{pspicture} \Leftrightarrow
\begin{pspicture}[shift=-1.2](-.3,-1.4)(.3,.8)
\psline(-.2,.8)(.2,.4)(.2,-.4)(-.2,-.8)
\psline(.2,.8)(.05,.65) \psline[arrowscale=1.5]{->}(-.2,-.1)(-.2,.1)
\psline(-.05,.55)(-.2,.4)(-.2,-.4)(-.05,-.55)
\psline(.05,-.65)(.2,-.8) \psline[arrowscale=1.5]{->}(.2,-.1)(.2,.1)
\end{pspicture}
~~:~~
\begin{pspicture}[shift=-1.2](-.7,-1.4)(.7,.8)
\psline(-.6,.8)(-.5,0)(-.6,-.8)
\psline(.6,.8)(.5,0)(.6,-.8)
\pscircle[fillstyle=solid,fillcolor=darkgray](-.5,0){.1}
\pscircle[fillstyle=solid,fillcolor=darkgray](.5,0){.1}
\rput[t](-.9,-1){$(i)$}
\end{pspicture}
\Leftrightarrow
\begin{pspicture}[shift=-1.2](-.7,-1.4)(.7,.8)
\psline(-.6,.8)(-.5,0)(-.6,-.8)
\psline(.6,.8)(.5,0)(.6,-.8)
\pccurve[angleA=45,angleB=135](-.5,0)(.5,0)
\pccurve[angleA=-45,angleB=-135](-.5,0)(.5,0)
\pscircle[fillstyle=solid,fillcolor=darkgray](-.5,0){.1}
\pscircle[fillstyle=solid,fillcolor=darkgray](.5,0){.1}
\rput(0,.4){$+$}
\rput(0,-.4){$-$}
\end{pspicture}
\quad\quad
\begin{pspicture}[shift=-1.2](-.3,-1.4)(.3,.8)
\psline(-.2,-.8)(-.2,.8) \psline[arrowscale=1.5]{->}(-.2,-.1)(-.2,.1)
\psline(.2,-.8)(.2,.8) \psline[arrowscale=1.5]{->}(.2,.1)(.2,-.1)
\end{pspicture} \Leftrightarrow
\begin{pspicture}[shift=-1.2](-.3,-1.4)(.3,.8)
\psline(-.2,.8)(.2,.4)(.2,-.4)(-.2,-.8)
\psline(.2,.8)(.05,.65) \psline[arrowscale=1.5]{<-}(-.2,-.1)(-.2,.1)
\psline(-.05,.55)(-.2,.4)(-.2,-.4)(-.05,-.55)
\psline(.05,-.65)(.2,-.8) \psline[arrowscale=1.5]{<-}(.2,.1)(.2,-.1)
\end{pspicture}
~~:~~
\begin{pspicture}[shift=-1.2](-.2,-1.4)(.2,.8)
\psline(0,.8)(0,-.8)
\pscircle[fillstyle=solid,fillcolor=darkgray](0,0){.1}
\rput[t](-.5,-1){$(ii)$}
\end{pspicture}
\Leftrightarrow
\begin{pspicture}[shift=-1.2](-.2,-1.4)(.2,.8)
\psline(0,.8)(0,-.8)
\pscircle[fillstyle=solid,fillcolor=darkgray](0,-.6){.1}
\pscircle[fillstyle=solid,fillcolor=darkgray](0,.6){.1}
\pscircle[fillstyle=solid,fillcolor=darkgray](0,0){.1}
\rput(.3,.3){$-$}
\rput(.3,-.3){$+$}
\end{pspicture}
$$
\caption{The effects on the Seifert graph by Reidemeister move II with $(i)$ the same direction and
$(ii)$ the opposite direction on two parallel segments of a link.} \label{Reide2}
\end{figure}
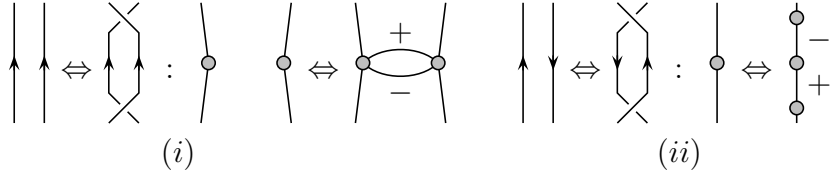

In general, there may not exist a spanning path in a Seifert graph.
However, in the case of $\Gamma(S(7_5))$, there exists a spanning
path as given in Figure~\ref{graph} $(iii)$,
thus, the alternating labeling $\mu$ on the spanning path in in Figure~\ref{graphtree1} $(ii)$ will
satisfy Theorem~\ref{maintheorem1}.
To obtain a flat plumbing basket surface of a link $L$ from a
canonical Seifert surface $S$ of the given link $L$,
we need to find a co-tree edge alternating spanning tree $T$.
The existence of such a co-tree edge alternating spanning tree can be
stated in the language of graph theory as in Theorem~\ref{maintheorem1}.

Since the Reidemeister move II will be frequently used in section~\ref{secbasket},
we demonstrate its effects on a Seifert graph in Figure~\ref{Reide2}.

\section{A main theorem for the algorithm} \label{premain}

Before we provide a proof of Theorem~\ref{maintheorem1}, let us review some definitions
which will be used in the proof.

\begin{defi}
Let $\Gamma$ be a connected finite bipartite graph with a
spanning tree $T$, a vertex $v \in V(\Gamma)$.
A \emph{depth labeling}  $\mu_v : E(T) \rightarrow \{ +, - \}$
of $T$ is defined as follows;
we define the \emph{depth} $d(u)$ of a vertex $u\in V(T)$ to
be the distance between the vertices $u$ and $v$,
the \emph{depth} $d(e)$ of an edge $e=\{v_i,v_j\}\in E(T)$ to
be the maximum of the depth of $v_i$ and $v_j$ and
the \emph{depth labeling} $\mu_v (e)$ of an edge $e=\{u,w\} \in T$
to be the sign of $(-1)^{d(e)}$ as illustrated in Figure~\ref{graph1} $(iii)$.
The vertices adjacent to $v$ whose depths are one more than
the depth of $v$ are called \emph{the children} of $v$.
For a spanning tree $T$ of $\Gamma$ with respect to the depth $d$,
one can see that the depths of the vertices in the path in $T$ joining the end vertices
$u$ and $w$ of a co-tree edge $e=\{u,w\} \in E(\Gamma)\setminus E(T)$ has the minimum
at $v_{\{u,w\}}$ which is the least common ancestor of the vertices $u$ and $w$,
let us call it \emph{the least common ancestor} of the edge $e$, denoted by $v_e$.
\end{defi}

\begin{theorem}
For a connected bipartite graph $\Gamma$ and a vertex $v$,
there exists an algorithm to determine a co-tree edge alternating spanning tree $T$
with respect to the depth labeling $\mu_v$.  \label{maintheorem1}
\begin{proof}

\begin{figure}
$$
\begin{pspicture}[shift=-2](-1.8,-2.8)(1.8,2.3)
\psline(-1,2)(1,0)(-1,1)(1,0)
\psline(1,0)(-1,0)(1,-1)(-1,-1)(1,0)(-1,1)(1,-1)(-1,-2)(1,1)(-1,1)
\rput(-1,2){$\bullet$} \rput(-1,1){$\bullet$}
\rput(-1,0){$\bullet$} \rput(-1,-1){$\bullet$}
\rput(-1,-2){$\bullet$} \rput(1,1){$\bullet$}
\rput(1,0){$\bullet$} \rput(1,-1){$\bullet$}
\rput(1.3,0){$v$} \rput(-1.5,0){$\Gamma$}
\rput[t](0,-2.3){$(i)$}
\end{pspicture} \quad
\begin{pspicture}[shift=-2](-1.8,-2.8)(1.8,2.5)
\psline(-1.5,1)(0,2)(-.5,1)(-.5,0)(.5,-1)
\psline(.5,1)(0,2)(1.5,1)(.5,0)(.5,1)(0,2)(-.5,1)(.5,0)(.5,-1)
\rput(0,2){$\bullet$} \rput(-1.5,1){$\bullet$}
\rput(-.5,1){$\bullet$} \rput(.5,1){$\bullet$}
\rput(1.5,1){$\bullet$} \rput(-.5,0){$\bullet$}
\rput(.5,0){$\bullet$} \rput(.5,-1){$\bullet$}
\rput(0,2.3){$v$} \rput[t](-1.5,2){$\Gamma$}
\rput[t](0,-2.3){$(ii)$}
\end{pspicture}
\quad
\begin{pspicture}[shift=-2](-1.8,-2.8)(1.8,2.5)
\psline(-1.5,1)(0,2)(-.5,1)(-.5,0)
\psline[linestyle=dashed](-.5,0)(.5,-1)
\psline(.5,1)(0,2)(1.5,1)
\psline[linestyle=dashed](1.5,1)(.5,0)
\psline(.5,0)(.5,1)(0,2)(-.5,1)
\psline[linestyle=dashed,linewidth=1.4pt](-.5,1)(.5,0)
\psline(.5,0)(.5,-1)
\rput(0,2){$\bullet$} \rput(-1.5,1){$\bullet$}
\rput(-.5,1){$\bullet$} \rput(.5,1){$\bullet$}
\rput(1.5,1){$\bullet$} \rput(-.5,0){$\bullet$}
\rput(.5,0){$\bullet$} \rput(.5,-1){$\bullet$}
\rput(-1,1.6){$-$} \rput(-.5,1.4){$-$}
\rput(1,1.6){$-$} \rput(.1,1.4){$-$}
\rput(-.7,.5){$+$} \rput(.3,.5){$+$}
\rput(.3,-.5){$-$} \rput(0,2.3){$v=v_{e_1}$}
\rput(-1.7,.7){$v_1$} \rput(-.8,1){$v_2$}
\rput(.8,1){$v_3$} \rput(1.7,.7){$v_4$}
\rput(0,.8){$e_1$} \rput(-.8,0){$v_5$}
\rput(.8,0){$v_6$} \rput(.8,-1){$v_7$}
\rput[t](-1.5,2){$T$} \rput[t](0,-2.3){$(iii)$}
\end{pspicture}
$$
$$
\begin{pspicture}[shift=-2](-1.8,-1.8)(1.8,2.5)
\psline(-1.5,1)(0,2)
\psline[linestyle=dashed](-.5,1)(0,2)
\psline(.5,1)(0,2) \psline(1.5,1)(0,2)
\psline(.5,0)(.5,1) \psline(-.5,0)(-.5,1)
\psline[linestyle=dashed,linewidth=1.4pt](.5,0)(1.5,1)
\psline(.5,0)(-.5,1) \psline(.5,-1)(.5,0)
\psline[linestyle=dashed](.5,-1)(-.5,0)
\rput(0,2){$\bullet$} \rput(-1.5,1){$\bullet$}
\rput(-.5,1){$\bullet$} \rput(.5,1){$\bullet$}
\rput(1.5,1){$\bullet$} \rput(-.5,0){$\bullet$}
\rput(.5,0){$\bullet$} \rput(.5,-1){$\bullet$}
\rput(-1,1.6){$-$} \rput(-.2,.4){$-$}
\rput(1,1.6){$-$} \rput(.1,1.4){$-$}
\rput(-.7,.5){$+$} \rput(.3,.5){$+$}
\rput(.3,-.4){$-$} \rput(0,2.3){$v=v_{e_2}$}
\rput(-1.7,.7){$v_1$} \rput(-.8,1){$v_5$}
\rput(.8,1){$v_2$} \rput(1.7,.7){$v_3$}
\rput(1.2,.4){$e_2$} \rput(-.8,0){$v_7$}
\rput(.8,-.1){$v_4$} \rput(.8,-1){$v_6$}
\rput[t](-1.5,2){$T_1$} \rput[t](0,-1.5){$(iv)$}
\end{pspicture}
\quad
\begin{pspicture}[shift=-2](-1.8,-1.8)(1.8,2.5)
\psline(-1.5,1)(0,2)
\psline[linestyle=dashed](-.5,1)(0,2)
\psline(.5,1)(0,2)
\psline[linestyle=dashed](1.5,1)(0,2)
\psline(.5,0)(.5,1) \psline(-.5,0)(-.5,1)
\psline(.5,0)(1.5,1) \psline(.5,0)(-.5,1)
\psline(.5,-1)(.5,0) \psline[linestyle=dashed,linewidth=1.4pt](.5,-1)(-.5,0)
\rput(0,2){$\bullet$} \rput(-1.5,1){$\bullet$}
\rput(-.5,1){$\bullet$} \rput(.5,1){$\bullet$}
\rput(1.5,1){$\bullet$} \rput(-.5,0){$\bullet$}
\rput(.5,0){$\bullet$} \rput(.5,-1){$\bullet$}
\rput(-1,1.6){$-$} \rput(-.2,.4){$-$}
\rput(1.2,.4){$-$} \rput(.1,1.4){$-$}
\rput(-.7,.5){$+$} \rput(.3,.5){$+$}
\rput(.3,-.4){$-$} \rput(0,2.3){$v$}
\rput(-1.7,.7){$v_1$} \rput(-.8,1){$v_4$}
\rput(.8,1){$v_2$} \rput(1.7,.7){$v_6$}
\rput(-.1,-.7){$e_3$} \rput(-.8,0){$v_7$}
\rput(1.35,-.1){$v_3=v_{e_3}$}
\rput(.8,-1){$v_5$} \rput[t](-1.5,2){$T_2$}
\rput[t](0,-1.5){$(v)$}
\end{pspicture}
\quad
\begin{pspicture}[shift=-2](-1.8,-1.8)(1.8,2.5)
\psline(-1.5,1)(0,2) \psline[linestyle=dashed](-.5,1)(0,2)
\psline(.5,1)(0,2) \psline[linestyle=dashed](1.5,1)(0,2)
\psline(.5,0)(.5,1) \psline(-.5,0)(-.5,1)
\psline(.5,0)(1.5,1) \psline(.5,0)(-.5,1)
\psline[linestyle=dashed](.5,-1)(.5,0)
\psline(.5,-1)(-.5,0) \rput(0,2){$\bullet$}
\rput(-1.5,1){$\bullet$} \rput(-.5,1){$\bullet$}
\rput(.5,1){$\bullet$} \rput(1.5,1){$\bullet$}
\rput(-.5,0){$\bullet$} \rput(.5,0){$\bullet$}
\rput(.5,-1){$\bullet$} \rput(-1,1.6){$-$}
\rput(-.2,.4){$-$} \rput(1.2,.4){$-$}
\rput(.1,1.4){$-$} \rput(-.7,.5){$+$}
\rput(.3,.5){$+$} \rput(-.3,-.6){$-$}
\rput(0,2.3){$v$} \rput(-1.7,.7){$v_1$}
\rput(-.8,1){$v_4$} \rput(.8,1){$v_2$}
\rput(1.7,.7){$v_5$} \rput(-.8,0){$v_6$}
\rput(1.35,-.1){$v_3$} \rput(.8,-1){$v_7$}
\rput[t](-1.5,2){$\overline T$}
\rput[t](0,-1.5){$(vi)$}
\end{pspicture}
$$
\caption{$(i)$ A connected bipartite graph $\Gamma$ with a fixed vertex $v$,
$(ii)$ a drawing of $\Gamma$ with respect to the depth from the root vertex $v$,
$(iii)$ a depth labeling $\mu_v$ of $T$ and
$(iv-vi)$ an algorithmic process to change $(T, \mu_v)$ to the
desired $({\overline T}, {\overline \mu}_v)$.}  \label{graph1}
\end{figure}
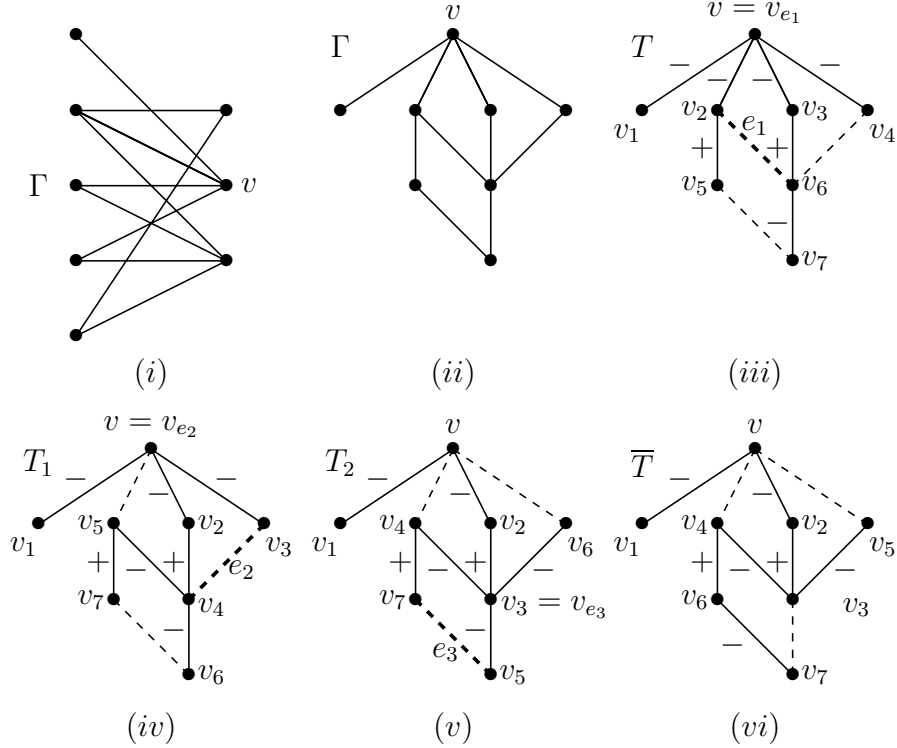

From a connected bipartite graph $\Gamma$ and the fixed vertex $v$,
we redraw $\Gamma$ with respect to the depth from the root
vertex $v$ as illustrated in Figure~\ref{graph1} $(ii)$.
Then, we choose a spanning tree $T$ in $\Gamma$ and consider
the depth labeling $\mu_d$ of $T$.
If $T$ is a co-tree  edge alternating spanning tree with respect
to the depth labeling $\mu_v$,
then we have a desired spanning tree and a labeling on $T$ and
it is easy to see that the least common ancestor of each co-tree
edge $e=\{u,w\} \in E(\Gamma)\setminus E(T)$ must be either one of
vertices $u$ or $w$.

Suppose $T$ is not a co-tree edge alternating spanning tree
with respect to the depth labeling $\mu_v$.
Then, there exists a cotree edge $e \in E(\Gamma)\setminus E(T)$ such that
the unique path $P_e$ in $T$ joining both end vertices of the edge $e$
does not have an alternating signs with respect to the depth labeling $\mu_v$.
This may restate that there exists a co-tree edge $e \in E(\Gamma)\setminus E(T)$ such that
the least common ancestor $v_e$ of the edge $e=\{u,w\}$ is neither $u$ nor $w$.
Then, the path $P_e$ is the union of two paths $P_u$ joining
the vertices $v_e$ and $u$ and $P_w$ joining vertices $v_e$ and $w$.
Each of these two paths $P_u$ and $P_w$ has alternating signs with respect
to the depth labeling $\mu_d$.
But the two children of $v_e$ in $P_u$ and $P_w$ have the same sign.
Since the lengths of two paths $P_u$ and $P_w$ are odd and even,
we can choose the shorter one.
We remove one of children edges of $v_e$ which belong to the
shorter path and add the edge $e$ to get a new spanning tree $\overline{T}$.
Since the new co-tree edge $e'$ produced by the algorithm has the
property that the least common ancestor of the edge $e'$ is either one of
vertices of $e'$, inductively we can remove all co-tree edges
for which the unique path in $T$ joining both end vertices of the edge
does not have an alternating signs with respect to the depth labeling $\mu_v$.

This process does depend on the order of the edges $e$'s, thus we
repeat the process among the edges for which the least common
ancestor of the end vertices of the edge are not an end vertex of the edge
that has the minimal depth at the the least common ancestor and if more than
two edges' least common ancestors have the same depth, we use an
order given by the new labeling of vertices.
Here, let us deal with an example in Figure~\ref{graph1} $(iii)$.
There are three co-tree  edges $\{v_2,v_6\}, \{v_4,v_6\}$ and $\{v_5,v_7\}$
whose the least common ancestors of the edges are the root $v$.
But the vertex $v_2$ has the smallest labeling on $v_i$ in these edges.
So, we pick the edge $e_1 = \{v_2,v_6\}$.
For algorithmic process, we remove the edge $\{v,v_2\}$ and add the edge
$e_1 = \{v_2,v_6\}$ to have a new spanning tree $T_1$ with a new
vertex labeling as in Figure~\ref{graph1} $(iv)$.
Now there are only two co-tree  edges $\{v_3,v_5\}$ and $\{v_6,v_7\}$
whose the least common ancestor of the edge is not the end vertices of the edge.
The least common ancestor of the edge $\{v_3,v_5\}$ has the depth $0$, so pick
the edge $e_2 = \{v_3,v_5\}$. By repeating the process, we get $T_2$ in Figure~\ref{graph1} $(v)$.
By the similar reasoning, we choose the edge $e_2 = \{v_5,v_7\}$ and we finally get the
desired spanning tree and a labeling $({\overline T}, {\overline \mu}_v)$
as in Figure~\ref{graph1} $(vi)$.

For each new spanning tree in the process, we named the vertices in
lexicographical order, the first by the depth of the vertices and
the second from the left to the right as demonstrated in Figure~\ref{graph1} $(iii-vi)$.

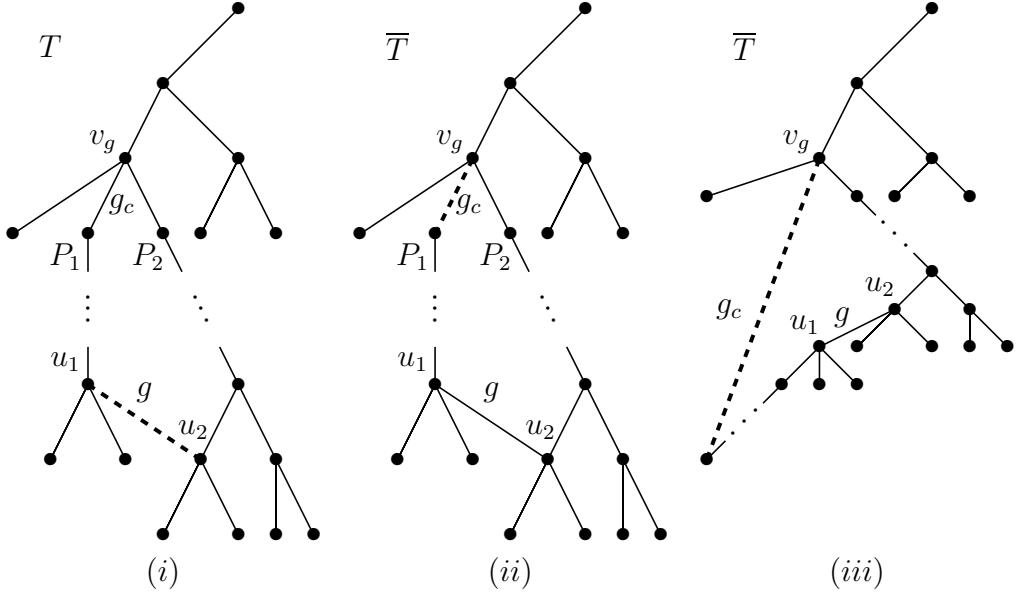
\begin{figure}
$$
\begin{pspicture}[shift=-2](-.1,-4.7)(4.1,3.2)
\psline(3,3)(2,2)(1.5,1)(0,0)
\psline(1.5,1)(1,0)(1,-.5)
\psline(1,-1.5)(1,-2)(.5,-3)(1,-2)(1.5,-3)
\psline(1.5,1)(2,0)(2.25,-.5)
\psline(2.75,-1.5)(3,-2)(2.5,-3)(2,-4)(2.5,-3)(3,-4)
\psline(3,-2)(3.5,-3)(3.5,-4)(3.5,-3)(4,-4)
\psline(2,2)(3,1)(2.5,0)(3,1)(3.5,0)
\psline[linestyle=dashed, linewidth=1.5pt](1,-2)(2.5,-3)
\rput(3,3){$\bullet$} \rput(2,2){$\bullet$}
\rput(1.5,1){$\bullet$} \rput(3,1){$\bullet$}
\rput(0,0){$\bullet$} \rput(1,0){$\bullet$}
\rput(2,0){$\bullet$} \rput(2.5,0){$\bullet$}
\rput(3.5,0){$\bullet$} \rput(1,-2){$\bullet$}
\rput(3,-2){$\bullet$} \rput(.5,-3){$\bullet$}
\rput(1.5,-3){$\bullet$} \rput(2.5,-3){$\bullet$}
\rput(3.5,-3){$\bullet$} \rput(2,-4){$\bullet$}
\rput(3,-4){$\bullet$} \rput(3.5,-4){$\bullet$}
\rput(4,-4){$\bullet$} \rput(1.2,1.23){$v_g$}
\rput(1.45,.4){$g_c$} \rput(.7,-1.7){$u_1$}
\rput(2.4,-2.6){$u_2$} \rput(.7,-.3){$P_1$}
\rput(1.8,-.3){$P_2$} \rput(1.75,-2.1){$g$}
\rput(1,-1){$\cdot$} \rput(1,-.85){$\cdot$}
\rput(1,-1.15){$\cdot$} \rput(2.5,-1){$\cdot$}
\rput(2.425,-.85){$\cdot$} \rput(2.575,-1.15){$\cdot$}
\rput(.5,2.5){$T$} \rput[t](2,-4.3){$(i)$}
\end{pspicture} \quad
\begin{pspicture}[shift=-2](-.1,-4.7)(4.1,3.2)
\psline(3,3)(2,2)(1.5,1)(0,0)
\psline[linestyle=dashed, linewidth=1.5pt](1.5,1)(1,0)
\psline(1,0)(1,-.5)
\psline(1,-1.5)(1,-2)(.5,-3)(1,-2)(1.5,-3)
\psline(1.5,1)(2,0)(2.25,-.5)
\psline(2.75,-1.5)(3,-2)(2.5,-3)(2,-4)(2.5,-3)(3,-4)
\psline(3,-2)(3.5,-3)(3.5,-4)(3.5,-3)(4,-4)
\psline(2,2)(3,1)(2.5,0)(3,1)(3.5,0)
\psline(1,-2)(2.5,-3)
\rput(3,3){$\bullet$} \rput(2,2){$\bullet$}
\rput(1.5,1){$\bullet$} \rput(3,1){$\bullet$}
\rput(0,0){$\bullet$} \rput(1,0){$\bullet$}
\rput(2,0){$\bullet$} \rput(2.5,0){$\bullet$}
\rput(3.5,0){$\bullet$} \rput(1,-2){$\bullet$}
\rput(3,-2){$\bullet$} \rput(.5,-3){$\bullet$}
\rput(1.5,-3){$\bullet$} \rput(2.5,-3){$\bullet$}
\rput(3.5,-3){$\bullet$} \rput(2,-4){$\bullet$}
\rput(3,-4){$\bullet$} \rput(3.5,-4){$\bullet$}
\rput(4,-4){$\bullet$} \rput(1.2,1.23){$v_g$}
\rput(1.45,.4){$g_c$} \rput(.7,-1.7){$u_1$}
\rput(2.4,-2.6){$u_2$} \rput(.7,-.3){$P_1$}
\rput(1.8,-.3){$P_2$} \rput(1.75,-2.1){$g$}
\rput(1,-1){$\cdot$} \rput(1,-.85){$\cdot$}
\rput(1,-1.15){$\cdot$} \rput(2.5,-1){$\cdot$}
\rput(2.425,-.85){$\cdot$} \rput(2.575,-1.15){$\cdot$}
\rput(.5,2.5){$\overline T$} \rput[t](2,-4.3){$(ii)$}
\end{pspicture}
\quad
\begin{pspicture}[shift=-2](-.1,-3.7)(4.1,4.2)
\psline(3,4)(2,3)(1.5,2)(0,1.5)
\psline[linestyle=dashed, linewidth=1.5pt](1.5,2)(0,-2)
\psline(.25,-1.75)(0,-2)
\psline(2.5,0)(1.5,-.5)(1,-1)(.75,-1.25)
\psline(1.5,-1)(1.5,-.5)(2,-1)
\psline(1.5,2)(2,1.5)(2.25,1.25)
\psline(2.75,.75)(3,.5)(2.5,0)(2,-.5)(2.5,0)(3,-.5)
\psline(3,.5)(3.5,0)(3.5,-.5)(3.5,0)(4,-.5)
\psline(2,3)(3,2)(2.5,1.5)(3,2)(3.5,1.5)
\rput(3,4){$\bullet$} \rput(2,3){$\bullet$}
\rput(1.5,2){$\bullet$} \rput(3,2){$\bullet$}
\rput(0,1.5){$\bullet$} \rput(0,-2){$\bullet$}
\rput(2,1.5){$\bullet$} \rput(2.5,1.5){$\bullet$}
\rput(3.5,1.5){$\bullet$} \rput(1.5,-.5){$\bullet$}
\rput(3,.5){$\bullet$} \rput(1,-1){$\bullet$}
\rput(2,-1){$\bullet$} \rput(2.5,0){$\bullet$}
\rput(3.5,0){$\bullet$} \rput(2,-.5){$\bullet$}
\rput(3,-.5){$\bullet$} \rput(3.5,-.5){$\bullet$}
\rput(1.5,-1){$\bullet$} \rput(4,-.5){$\bullet$}
\rput(2.35,1.15){$\cdot$} \rput(2.5,1){$\cdot$}
\rput(2.65,.85){$\cdot$} \rput(.35,-1.65){$\cdot$}
\rput(.5,-1.5){$\cdot$} \rput(.65,-1.35){$\cdot$}
\rput(1.2,2.23){$v_g$}
\rput(.3,0){$g_c$} \rput(1.3,-.2){$u_1$}
\rput(2.3,.3){$u_2$}  \rput(1.8,-.1){$g$}
\rput(.5,3.5){$\overline T$}
\rput[t](2,-3.3){$(iii)$}
\end{pspicture}
$$
\caption{$(i)$ A part of spanning tree $T$ including $g, g_c, P_1, P_2$,
$(ii)$ the resulting spanning tree $\overline T$ by the process described in the proof and
$(iii)$ a new diagram of $\overline T$ with respect to the depth.}  \label{graph3}
\end{figure}

Next, we prove that this algorithm always produce the desired
spanning tree and an alternating labeling on it.
Suppose that our algorithmic process does not work.
Let $\Gamma$ be a counterexample. For any vertex $v$ and any spanning tree $T$ of $\Gamma$,
let $\A(T,v)$ be the set of all co-tree edges
$e \in E(\Gamma)\setminus E(T)$ whose least common ancestor $v_e$
of the end vertices of the edge $e$ is not an end vertex of the edge $e$.
If $\A=\emptyset$, then $T$ is a co-tree  edge alternating spanning tree.
Thus, it contradicts the hypothesis that $(\Gamma, T)$ is a counterexample.

For any vertex $v$ and any spanning tree $T$ of $\Gamma$,
denoted by an ordered triple $(\Gamma, T, v)$, we define $\eta(\Gamma, T, v)$ to be
$$\eta(\Gamma, T, v) = \sum_{ u \in V(\Gamma)} d_v (u)$$
where $d_v (u)$ is the distance between the vertices $u$ and $v$,
$i. e.$ the number of edges in the path joining the vertices $u$ and $v$.
Since $\Gamma$ is a counterexample, there exist a vertex $v$ and
a spanning tree $T$ of $\Gamma$ such that $\eta(\Gamma, T, v)$ is maximal.
By the above claim, $\A\not=\emptyset$. But if so, we will show that there exists
a spanning tree $\overline{T}$ of $\Gamma$ such that
$\eta(\Gamma, T, v) < \eta(\Gamma, \overline{T}, v)$.
This contraction completes the proof of the theorem. In fact, our claim is that
a fixed vertex $v$ and a spanning tree $T$ of the maximum $\eta(\Gamma, T, v)$ must
be a co-tree  edge alternating spanning tree with respect to the depth coloring $\mu_v$.

Suppose if $\A\not=\emptyset$, then there exists a co-tree edge
$f \in E(\Gamma)\setminus E(T)$ whose
least common ancestor $v_f$ of the end vertices of the
edge $f$ is not an end vertex of the edge $f$.
Among all such edges $f \in \A$, we pick an edge $g=\{u_1, u_2\}$
that has the minimal $d(v_g)$.
Let $P_g$ be the path in $T$ joining $u_1$ and $u_2$. Since the edge
$g \in \C$, the path $P_g$ is the union of two paths, $P_1$
joining the vertices $u_1$ and $v_g$ and the path
$P_2$ joining the vertices $u_2$ and $v_g$. We further assume that
the path $P_1$ has the shorter length then $P_2$ as
depicted in Figure~\ref{graph3} $(i)$.
If the edge $g \in E(\Gamma)\setminus E(T)$ replaces the child edge $g_{c}$ of $v_g$
which belongs to the shorter path $P_1$, then we get a new spanning tree $\overline{T}$.
Using the same vertex $v$ as a root, we have a new depth $\overline d$ on $\overline{T}$.
We claim that $\eta(\Gamma, T, v) < \eta(\Gamma, \overline{T}, v)$ as follows.
We denote that $E(\Gamma)\setminus E(T) = \{ g$, $e_1$, $e_2$, $\ldots$, $e_n \}$
and $E(\Gamma)\setminus E({\overline T}) = \{ g_c, e_1, e_2, \ldots, e_n \}$.
For a vertex $u\in V(\Gamma)$, there are two possible cases;
an ancestor of the vertex $u$ is either not in the path $P_1$
or it is. Since the paths $P_1$ and $P_2$ have none zero length,
both cases do occur if we consider all vvertices of $\Gamma$. 
If an ancestor of the vertex $u$ is not in the path $P_1$,
it is easy to see that the depth of the path joining the
vertices $u$ and $v$ has not changed. Thus, we have
$d_v(u)= {\overline d}_v (u)$.
Otherwise, we choose a vertex
$w\in V(\Gamma)$ to have the maximal depth among all such
ancestors of the vertex $u$.
Let $Q$ be the path joining the vertices $u$ and $w$ and
$R$ be the path joining the vertices $f_1$ and $w$.
Then,
$$d_v(u) = d_v(w) + \ell(Q) \le d_v(f_1) + \ell(Q) <
{\overline d}_v(f_1) + \ell(Q) + \ell(R)={\overline d}_v(u).$$
where $\ell(P)$ is the length of the path $P$.

Therefore, by summing all $u \in V(\Gamma)$, we have
$\eta(\Gamma, T, v) < \eta(\Gamma, \overline{T}, v)$. However,
this contradicts our assertion that $(\Gamma, T, v)$ has
the maximal  $\eta(\Gamma, T, v)$ over all spanning trees $T$. This completes
the proof of Theorem~\ref{maintheorem1}.
\end{proof}
\end{theorem}

\section{Applications in knot theory : Plumbing numbers}\label{secbasket}

In this section, we find a few applications of the main theorem in knot theory.
There are several interesting plumbing
surfaces\cite{Rudolph:plumbing} but, in particular,  we would like to discuss three
plumbing surfaces and their plumbing numbers. We will
review the definitions of these plumbing
surfaces and define corresponding plumbing numbers and
prove theorems for the upper bounds of each plumbing
number.

\subsection{Definitions}\label{prelim}

\begin{figure}
$$
\begin{pspicture}[shift=-1.6](-.2,-2.2)(6,1.7)
\psline[linewidth=2pt](0,-1.5)(.8,-.5)
\psline[linecolor=darkgray, linewidth=1.5pt](1.2,0)(4.5,0)
\psline[linestyle=dashed, linewidth=2pt](.8,-.5)(2,1)
\psline[linewidth=2pt](2.4,1.5)(2,1)
\pccurve[angleA=110,angleB=180](.8,-.5)(1.4,1)
\pccurve[linestyle=dashed, angleA=0,angleB=90](1.4,1)(1.6,.5)
\psline(1.4,1)(4.4,1)
\pccurve[angleA=110,angleB=180](4.1,-.5)(4.4,1)
\pccurve[angleA=0,angleB=90](4.4,1)(4.9,.5)
\psline[linewidth=2pt](3.3,-1.5)(5.7,1.5)
\psline(.8,-.5)(4.1,-.5) \psline[linestyle=dashed](1.6,.5)(4.9,.5)
\pccurve[angleA=200,angleB=0](5.7,1.5)(4.6,1.7)
\pccurve[angleA=180,angleB=0](4.6,1.7)(3.5,1.3)
\pccurve[angleA=180,angleB=-20](3.5,1.3)(2.4,1.5)
\pccurve[angleA=200,angleB=0](3.3,-1.5)(2.2,-1.7)
\pccurve[angleA=180,angleB=0](2.2,-1.7)(1.1,-1.3)
\pccurve[angleA=180,angleB=-20](1.1,-1.3)(0,-1.5)
\rput(3,.2){$\alpha$} \rput(2.5,-.28){$C_{\alpha}$}
\rput(2,-1){$S$} \rput(1.25,.6){$B_{\alpha}$}
\rput[t](2.7,-2){$(i)$}
\end{pspicture}\quad
\begin{pspicture}[shift=-1.6](-.2,-2.2)(6,1.7)
\psline[linewidth=2pt](0,-1.5)(.8,-.5)
\psline(.8,-.5)(1.6,.5) \psline[linewidth=2pt](1.6,.5)(2.4,1.5)
\pccurve[linewidth=2pt, angleA=45,angleB=135](1.6,.5)(4.9,.5)
\pccurve[linewidth=2pt, angleA=45,angleB=135](.8,-.5)(4.1,-.5)
\psline[linewidth=2pt](3.3,-1.5)(4.1,-.5)
\psline(4.1,-.5)(4.9,.5) \psline[linewidth=2pt](4.9,.5)(5.7,1.5)
\psline[linestyle=dotted](.8,-.5)(4.1,-.5)
\psline[linestyle=dotted](1.6,.5)(4.9,.5)
\pccurve[angleA=200,angleB=0](5.7,1.5)(4.6,1.7)
\pccurve[angleA=180,angleB=0](4.6,1.7)(3.5,1.3)
\pccurve[angleA=180,angleB=-20](3.5,1.3)(2.4,1.5)
\pccurve[angleA=200,angleB=0](3.3,-1.5)(2.2,-1.7)
\pccurve[angleA=180,angleB=0](2.2,-1.7)(1.1,-1.3)
\pccurve[angleA=180,angleB=-20](1.1,-1.3)(0,-1.5)
\psline(4.2,.2)(5.2,.2) \psline[arrowscale=1.5]{->}(5.18,.2)(5.2,.2)
\rput(2.5,-.28){$C_{\alpha}$} \rput(5.5,.2){$A_{0}$}
\rput(2,-1){$\overline{S}$}
\rput[t](2.7,-2){$(ii)$}
\end{pspicture}
$$
\caption{$(i)$ A geometric shape of $\alpha, B_{\alpha}$
and $C_{\alpha}$ on a Seifert surface $S$
and $(ii)$ a new Seifert surface $\overline{S}$ obtained
from $S$ by a top $A_0$ plumbing along the path
$\alpha$.}
\label{topfig}
\end{figure}
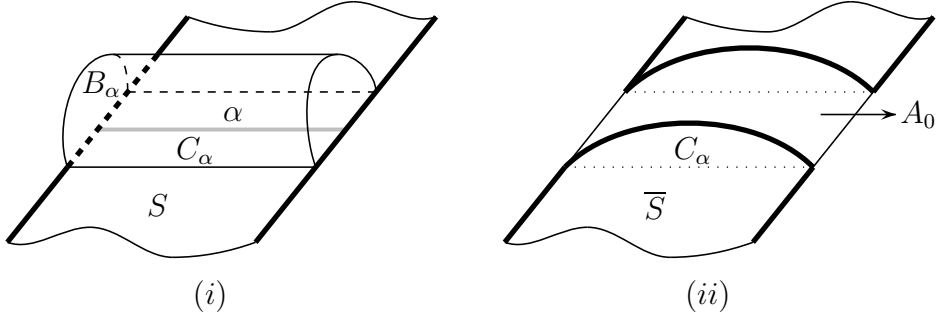

Spaces, maps, etc., are piecewise smooth unless stated differently.
Let $M$ be an oriented manifold. ${-}M$ denotes $M$ with its orientation reversed
and when notation requires it, $+M$ denotes $M$. For a suitable subset
$S \subset M$, $N_M (S)$ denotes a closed regular neighborhood of $S$ in $(M, \partial M)$
where an ordered pair $(S, T)$ stands a condition $T \subset S$ and a
map between ordered pairs $f: (S, T) \rightarrow (U,V)$ is a
map $f: S \rightarrow U$ which requires to preserve
subsets so that $f(T) \subset V$.
For a suitable codimension-$1$ submanifold $S \subset M$ (resp.,
submanifold pair
$(S, \partial S) \subset (M, \partial M))$, a emph{collaring}
is an orientation-preserving embedding
$S \times [0, 1] \rightarrow M$ (resp., $(S, \partial S)
\times[0, 1] \rightarrow (M, \partial M))$ extending $id_S = id_{S \times\{0\}}$;
a \emph{collar} of $S$ in $M$ (resp., of $(S, \partial S)$
in $(M, \partial M)$) is the image col$_M (S)$
(resp., col$_{(M, \partial M)} (S, \partial S)$) of a collaring.
The push-off of $S$ determined by
a collaring of $S$ or $(S, \partial S)$, denoted by $S^+$,
is the image by the collaring of
$S \times \{1\}$ with the orientation of $S$;
let $S^{-} :=$ ${-}$ $S^{+}$ such that $S$ and $S^{-}$ are
oriented submanifolds of the boundary of col$_M (S)$~\cite{Rudolph:plumbing}.

An \emph{arc} is a manifold homeomorphic to the unit interval $[0, 1]$.
An arc $\alpha$ in a $2$-manifold $S$ with a boundary is \emph{proper}
if $\partial \alpha = \alpha \cap \partial S$
Let $S$ be a Seifert surface in $\BS^3$. Let Top$(S)$ be col$_{\BS^3} (S)$ and
let $D^2$ be a $2$-dimensional disc.
Now, we first define the top plumbing as follows.

\begin{defi}(\cite{Rudolph:plumbing}) \label{topdef}
Let $\alpha$ be a proper arc on a Seifert surface $S$.
Let $C_{\alpha}$ be \emph{col}$_{(S, \partial S)} (\alpha, \partial \alpha)$ which is called
\emph{the gluing region}.
Let $B_{\alpha}$ be \emph{col}$_{(S, \partial S)} (\alpha, \partial \alpha)$ (so $B_{\alpha}$
is a $3$-cell in \emph{top}$(S)$, that is the positive normal to $S$
along $C_{\alpha}=S \cap B_{\alpha} \subset \partial B_{\alpha}$
points into $B_{\alpha}$) as depicted in Figure~\ref{topfig}.
Let $A_n\subset B_{\alpha}$ be an $n$-full twisted annulus such that
$A_n \cap \partial B_{\alpha}= C_{\alpha}$. Then \emph{top plumbing} on $S$ along a path $\alpha$
is the new surface $\overline{S}= S \cup A_n$ where $A_n, C_{\alpha}, B_{\alpha}$ satisfy
the previous conditions.
\end{defi}

Although the bottom plumbing was defined too, we only
use top plumbing for the rest of the article and simply call it a \emph{plumbing}.
Rudolph found a few interesting results including every arborescent
Seifert surfaces are baskets~\cite{Rudolph:plumbing}.

\begin{remark} \label{remark-1}
Two consecutive plumbings are non-commutative in general.
To obtain plumbing surfaces from a canonical Seifert surface,
some pairs of plumbings are not commutative.
\end{remark}

Because two plumbings are non-commutative as given in Remark~\ref{remark-1},
the order of plumbings has to be chosen carefully.
However, the exact order of plumbing will not be discuss in this article.
One may find the details in~\cite{CDK}.

Throughout Section~\ref{secbasket}, we will assume
all links are not splittable nd prime. A link $L$ is \emph{splittable} if there exists a $2$-sphere $\BS^2 \subset \BS^3$
such that the intersection of $\BS^2$ and the link is an emptyset and
each of two $3$-balls bounded by $\BS^2$ contains a nonempty subset of the link.
If there do exist a $2$-sphere $\BS^2 \subset \BS^3$
such that the intersection of $\BS^2$ and the link $L$ is a set of two points,
then we have two links $L_i$ which is obtained by the union of a 
path on $\BS^2$ joining two points and $L \cap B_i^3$
where $B_i^3$ is the $3$-ball bounded by $\BS^2$ for $i=1, 2$.
If none of $L_i$ is a trivial knot, the link $L$ is called \emph{composite}, denoted by $L=L_1 \# L_2$. 
A link $L$ is \emph{prime} if it is not composite.
Since the following plumbing numbers of a composite link $L=L_1 \# L_2$ is
the sum of the plumbing numbers of $L_i$ as shown in
Theorem~\ref{connectsum}, we can handle each $L_i$ separately.
A link diagram $D(L)$ of a link $L$ is \emph{reducible}
if the Reidemeister move I or mover II can be
performed to decrease the number of crossings.
If a link diagram is not reducible,
we say it is \emph{reduced}.
For the closed braid $\overline{\beta}$, the reducibility
is equivalent to say that there does not exist,
$\sigma_i\sigma_i^{-1}$, $\sigma_i^{-1}\sigma_i$ nor 
$\sigma_i$ or $\sigma_i^{-1}$ appear just once in $\beta$.

\begin{remark} \label{remark0}
\begin{enumerate}
\item
If a link $L$ is splittable, then we have to apply a Reidemeister move \emph{II} first.
If a closed braid $\overline{\beta}$, $\beta \in B_n$ is not splittable,
then for $i=1, 2, \ldots, n-1$, one of
$\sigma_i$, $\sigma_i^{-1}$, must be in $\beta$.

\item
If a closed braid $\overline{\beta}$, $\beta \in B_n$ is prime and only one of
$\sigma_i$, $\sigma_i^{-1}$ appears in $\beta$ for some $i \in \{1, 2, \ldots, n-1\}$,
then it has to appear more than once.
\end{enumerate}
\end{remark}

Now, we define a basket surface of a link $L$.

\begin{defi}(\cite{Rudolph:plumbing}) \label{basketdef}
Let $A_n \subset \mathbb{S}^3$ denote an $n$-twisted unknotted
annulus. A Seifert surface $S$ is a \emph{basket surface} if it is
$2$-disc $D^2$ or it can be constructed by plumbing $A_n$ to a
basket surface $S_0$ along a proper arc $\alpha \subset D^2
\subset S_0$, denoted by $S_0 *_{\alpha} A_n$.
If a link $L$ is a boundary of a basket surface, then it
is called a \emph{basket surface of $L$}.
The \emph{flat plumbing number} of $L$, denoted by $bk(L)$, is
the minimal number of flat annuli to obtain a basket surface of $L$.
\end{defi}

\begin{defi}(\cite{Rudolph:plumbing}) \label{fpdef}
A Seifert surface $S$ is a \emph{flat plumbing surface} if it is
$2$-disc $D^2$ or it can be constructed by plumbing $A_0$ to a
flat plumbing surface $S_0$ along a proper arc $\alpha \subset
S_0$, denoted by $S = S_0 *_{\alpha} A_0$.
If a link $L$ is a boundary of a flat plumbing surface, then it
is called a \emph{flat plumbing surface of $L$}.
The \emph{flat plumbing number} of $L$, denoted by $fp(L)$, is
the minimal number of flat annuli to obtain a flat plumbing surface of $L$.
\end{defi}

We note that for a flat plumbing surface, the gluing regions
$C_{\alpha}$ in the construction are not
necessarily contained in $D^2$. Hayashi and Wada showed every
oriented link is a boundary of a flat plumbing surface which is obtained by
finitely many flat plumbings~\cite{HW:plumbing}. Thus, the flat plumbing number
of a link $L$ is well defined.

\begin{figure}
$$
\begin{pspicture}[shift=-1.6](-2.3,-1.7)(2.5,1.7)
\psarc(.75,.25){1}{0}{90} \psarc(.75,-.25){1}{270}{0}
\psarc(-.75,.25){1}{90}{180} \psarc(-.75,-.25){1}{180}{270}
\psarc(.75,.25){.5}{0}{180} \psarc(.75,-.25){.5}{180}{0}
\psarc(-.75,.25){.5}{0}{180} \psarc(-.75,-.25){.5}{180}{0}
\psline(-.75,1.25)(.75,1.25) \psline(-.75,-1.25)(.75,-1.25)
\pccurve[angleA=-90,angleB=90](-1.75,.25)(-1.25,-.25)
\pccurve[angleA=-90,angleB=90](-.25,.25)(.25,-.25)
\pccurve[angleA=-90,angleB=90](1.25,.25)(1.75,-.25)
\pccurve[angleA=-90,angleB=45](1.75,.25)(1.55,.05)
\pccurve[angleA=-135,angleB=90](1.45,-.05)(1.25,-.25)
\pccurve[angleA=-90,angleB=45](.25,.25)(.1,.05)
\pccurve[angleA=-135,angleB=90](-.1,-.05)(-.25,-.25)
\pccurve[angleA=-90,angleB=45](-1.25,.25)(-1.45,.05)
\pccurve[angleA=-135,angleB=90](-1.55,-.05)(-1.75,-.25)
\end{pspicture}
\cong
\begin{pspicture}[shift=-1.6](-2.5,-1.7)(2.3,1.7)
\pccurve[doubleline=true, angleA=-70,angleB=70](1.48;90)(1.48;-90)
\pccurve[doubleline=true, angleA=-10,angleB=100](-.9,1.28)(1,-1.28)
\pccurve[doubleline=true, angleA=30,angleB=150](-1.92,0)(1.92,0)
\pccurve[doubleline=true, angleA=90,angleB=180](-1,-1.28)(.9,1.28)
\psellipse[fillstyle=none, linewidth=3pt](0,0)(2,1.5)
\end{pspicture}
$$
\caption{A flat plumbing basket surface of the trefoil knot.}
\label{figure8fpb}
\end{figure}
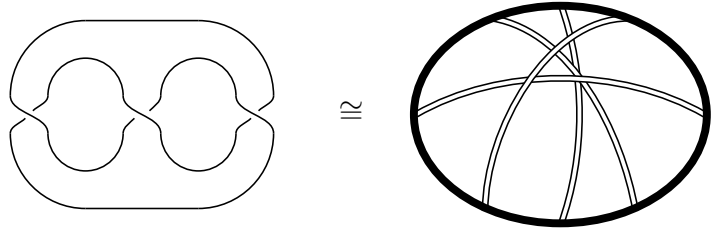

\begin{defi}(\cite{Rudolph:plumbing}) \label{fpbdef}
A Seifert surface $S$ is a \emph{flat plumbing basket surface} if
$S = D^2$ or if $S = S_0
*_{\alpha} A_0$ which can be constructed by plumbing $A_0$ to a
flat plumbing basket surface $S_0$ along a proper arc $\alpha \subset D^2\subset S_0$.
We say that a link $L$ admits a \emph{flat plumbing basket
presentation} if there exists a flat plumbing basket $S$ such that
$\partial S$ is equivalent to $L$.
The \emph{flat plumbing basket number} of $L$, denoted
by $fpbk(L)$, is the minimal number of flat annuli to obtain a
flat plumbing basket surface of $L$.
\end{defi}

\begin{figure}
$$
\begin{pspicture}[shift=-1.2](-.7,-1.8)(4.2,1.2)
\psarc[doubleline=true](2.5,0){1}{-5}{185}
\psarc[doubleline=true](2,0){1}{-5}{185}
\psarc[doubleline=true](1.5,0){1}{-5}{185}
\psarc[doubleline=true](1,0){1}{-5}{185}
\psframe[linecolor=lightgray,fillstyle=solid,fillcolor=lightgray](-.5,-1)(4,0)
\psline(-.03,0)(-.5,0)(-.5,-1)(4,-1)(4,0)(3.53,0)
\psline(.03,0)(.47,0) \psline(.53,0)(.97,0) \psline(1.03,0)(1.47,0)
\psline(1.53,0)(1.97,0) \psline(2.03,0)(2.47,0) \psline(2.53,0)(2.97,0)
\psline(3.03,0)(3.47,0)
\rput(1.75,-.5){{$D^2$}}
\rput(1.75,-1.7){{$(i)$}}
\end{pspicture} \quad
\begin{pspicture}[shift=-1.2](-.7,-1.8)(4.2,1.2)
\psarc[doubleline=true](2,0){1}{-5}{185}
\psarc[doubleline=true](2.5,0){1}{-5}{185}
\psarc[doubleline=true](1.5,0){1}{-5}{185}
\psarc[doubleline=true](1,0){1}{-5}{185}
\psframe[linecolor=lightgray,fillstyle=solid,fillcolor=lightgray](-.5,-1)(4,0)
\psline(-.03,0)(-.5,0)(-.5,-1)(4,-1)(4,0)(3.53,0)
\psline(.03,0)(.47,0) \psline(.53,0)(.97,0) \psline(1.03,0)(1.47,0)
\psline(1.53,0)(1.97,0) \psline(2.03,0)(2.47,0) \psline(2.53,0)(2.97,0)
\psline(3.03,0)(3.47,0)
\rput(1.75,-.5){{$D^2$}}
\rput(1.75,-1.7){{$(ii)$}}
\end{pspicture}
$$
\caption{Flat plumbing basket surfaces of $(i)$ the trefoil knot
and $(ii)$ the figure eight knot in the trivial open book decomposition.} \label{figure84band}
\end{figure}
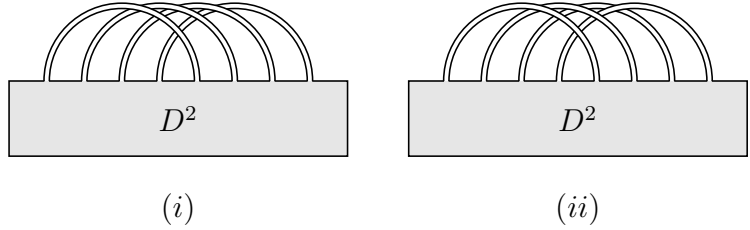

A flat plumbing basket surface of the trefoil knot is given in Figure~\ref{figure8fpb}.
In~\cite{FHK:openbook}, it is
shown that every link admits a flat plumbing basket presentation.
Thus, the flat plumbing basket number
of a link $L$ is well defined. Since every flat plumbing basket
surface is a basket surface, the basket number
of a link $L$ is also well defined.

An alternative definition of the flat plumbing basket
surfaces is given~\cite{FHK:openbook}
and this is very easy to follow. The \emph{trivial
open book decomposition} of $\mathbb{R}^3$
is a decomposition of $\mathbb{R}^3$ into the half planes
in the following form. In cylindrical coordinates, it can be presented
$$ \mathbb{R}^3 = \bigcup_{\theta \in [0, 2\pi)}
\{(r, \theta, z) | r \ge 0, z \in \mathbb{R} \}$$
where $\{(r, \theta, z) | r \ge 0, z \in \mathbb{R} \}$
is called a \emph{page} for $\theta \in [0, 2\pi)$.
Let $\mathcal{O}$ be the \emph{trivial open book decomposition} of the
$3$-sphere $\BS^3$ which is obtained from the trivial
open book decomposition of $\mathbb{R}^3$
by the one point compactification. A Seifert surface
is said to be a flat plumbing basket surface
if it consists of a single page of $\mathcal{O}$ and
finitely many bands which are
embedded in distinct pages~\cite{FHK:openbook}.
Flat plumbing basket surfaces
of $(i)$ the trefoil knot and $(ii)$ the figure
eight knot in the trivial open book decomposition are depicted 
in Figure~\ref{figure84band} where $D^2$ is presented as a shaded
rectangular region and the top horizontal line of
the rectangle is in the $z$-axis and the top
hemi-spherical annuli are contained in different pages.

\begin{remark} \label{remark1}
\begin{enumerate}
\item If one of basket, flat plumbing and flat plumbing
basket number of a link $L$ is zero, then $L$ is the unknot.
Equivalently, every nontrivial knots must have nonzero basket,
flat plumbing and flat plumbing basket numbers.

\item One can see that the boundary of these surfaces
with $n$-plumbing has at most $n+1$ components,
and the number of components is always congruent to $n+1$ modulo $2$.
Therefore, the basket, flat plumbing and flat plumbing
basket numbers of a knot have to be even.
\end{enumerate}
\end{remark}

For the last part of this subsection, let us prove the following theorem
which addresses plumbing numbers of composite links.

\begin{theorem} \label{connectsum}
The plumbing numbers of a composite link $L=L_1 \# L_2$ is
the sum of the plumbing numbers of $L_i$.
\begin{proof}
Let $S_i$ be plumbing surfaces of $L_i$ of the
plumbing number $n_i$ with the $2$-disc $D_i^2$ for $i=1, 2$. Then the band sum
of two surfaces where the connection was chosen to be a band between
the disc $D^2_1$ and $D^2_2$ is a plumbing surface of $L=L_1 \# L_2$ of
the plumbing number $n_1+n_2$. Thus, the plumbing number of $L=L_1 \# L_2$
is less than or equal to the sum of the plumbing numbers of $L_1$ and $L_2$.

For converse, let $S$ be a plumbing surface of $L$ of the plumbing number of $L$.
Since the link $L$ is composite, there exists a $2$-sphere $\BS^2 \subset \BS^3$
such that the intersection of $\BS^2$ and the link $L$ is a set of two points,
where $L_i$ is the link obtained by the union of a path on
$\BS^2$ joining two points and $L \cap B_i^3$
where $B_i^3$ is the $3$-ball bounded by $\BS^2$ for $i=1, 2$.
Since $S$ and $\BS^2$ are compact and piecewise smooth, we may assume
the intersection of these two surfaces are transversal, which implies
that the intersection is the union of $1$-dimensional manifolds.
However, by the assumption of $\BS^2$, the intersection of these two surfaces
are a single arc and a finite (possibly empty) set of circles. If we cap
off these circles by $2$-dimensional discs,
only a connected component $S_i$ in $B_i^3$ which has a boundary is a Seifert surfaces of $L_i$
and all the other components are homeomorphic to $\BS^2$ for $i=1, 2$. Then $S_1$ and $S_2$
are plumbing surfaces of $S_1$ and $S_2$ respectively and the connected sum of $S_1$ and $S_2$
is a plumbing surface of $L$ with the same plumbing number.
Thus, the plumbing number of $L=L_1 \# L_2$
is bigger than or equal to the sum of the plumbing numbers of $L_1$
and $L_2$. It completes the proof of the theorem.
\end{proof}
\end{theorem}

\subsection{Basket number} \label{baskets}

If the link $L$ is presented as a closure of a braid $\beta \in B_n$, we can choose
a braid of the following form $\beta=\sigma_{n-1}\sigma_{n-2}
\ldots \sigma_2\sigma_1 W(\sigma_1^{\pm 1}, \sigma_2^{\pm 1},
\ldots, \sigma_{n-1}^{\pm 1})$ where $W$ is a word in
$\sigma_i^{\pm 1}$~\cite{Birman:Braids}. The \emph{length} of the word $W$ is the number of
letters $\sigma_i^{\pm 1}$ in the word $W$. By a simple
modification of the idea of the main theorems
in~\cite[Theorem 2.4]{FHK:openbook} \cite{HW:plumbing}, we have the
following theorem that addresses an upper bound for the basket number of $L$.

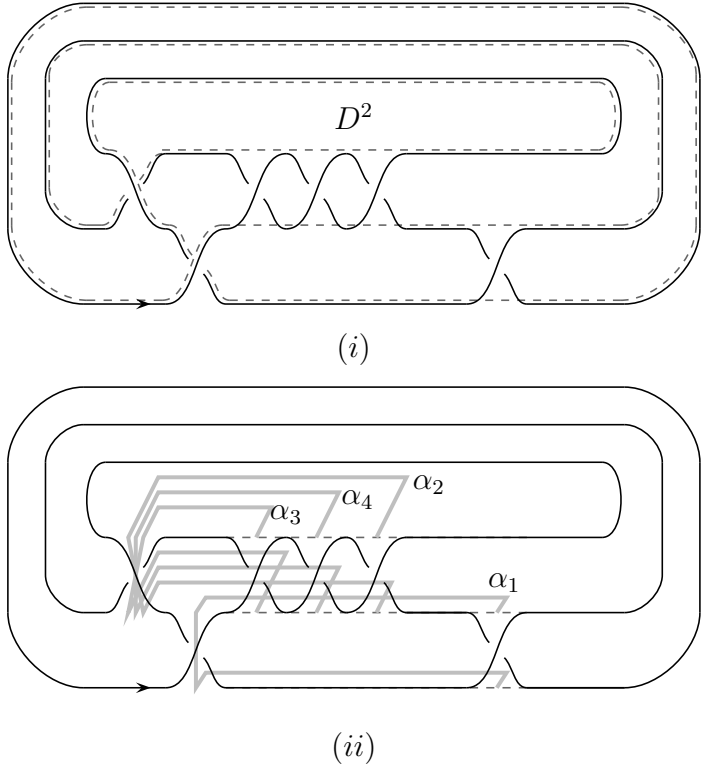
\begin{figure}
$$
\begin{pspicture}[shift=-1.2](0,-1.8)(9.5,3.1)
\psline[linecolor=pup, linestyle=dashed](.05,.05)(.05,1.95)
\pccurve[linecolor=pup, linestyle=dashed, angleA=90,angleB=180](.05,1.95)(1.05,2.95)
\psline[linecolor=pup, linestyle=dashed](1.05,2.95)(8.15,2.95)
\pccurve[linecolor=pup, linestyle=dashed, angleA=0,angleB=90](8.15,2.95)(9.15,2.05)
\psline[linecolor=pup, linestyle=dashed](9.15,2.05)(9.15,.05)
\pccurve[linecolor=pup, linestyle=dashed, angleA=-90,angleB=0](9.15,.05)(8.15,-.95)
\psline[linecolor=pup, linestyle=dashed](8.15,-.95)(2.95,-.95)
\pccurve[linecolor=pup, linestyle=dashed, angleA=180,angleB=-45](2.95,-.95)(2.65,-.55)
\pccurve[linecolor=pup, linestyle=dashed, angleA=135,angleB=0](2.45,-.35)(2.15,.05)
\pccurve[linecolor=pup, linestyle=dashed, angleA=180,angleB=0](2.15,.05)(1.35,1.05)
\pccurve[linecolor=pup, linestyle=dashed, angleA=180,angleB=180](1.35,1.05)(1.35,1.95)
\psline[linecolor=pup, linestyle=dashed](1.35,1.95)(7.85,1.95)
\pccurve[linecolor=pup, linestyle=dashed, angleA=0,angleB=0](7.85,1.95)(7.85,1.05)
\psline[linecolor=pup, linestyle=dashed](7.85,1.05)(2.05,1.05)
\pccurve[linecolor=pup, linestyle=dashed, angleA=180,angleB=45](2.05,1.000005)(1.75,.65)
\pccurve[linecolor=pup, linestyle=dashed, angleA=-135,angleB=0](1.55,.45)(1.25,.05)
\psline[linecolor=pup, linestyle=dashed](1.25,.05)(1.05,.05)
\pccurve[linecolor=pup, linestyle=dashed, angleA=180,angleB=-90](1.05,.05)(.55,.55)
\psline[linecolor=pup, linestyle=dashed](.55,.55)(.55,2.05)
\pccurve[linecolor=pup, linestyle=dashed, angleA=90,angleB=180](.55,2.05)(1.05,2.45)
\psline[linecolor=pup, linestyle=dashed](1.05,2.45)(8.15,2.45)
\pccurve[linecolor=pup, linestyle=dashed, angleA=0,angleB=90](8.15,2.45)(8.65,1.95)
\psline[linecolor=pup, linestyle=dashed](8.65,1.95)(8.65,.45)
\pccurve[linecolor=pup, linestyle=dashed, angleA=-90,angleB=0](8.65,.45)(8.15,.05)
\psline[linecolor=pup, linestyle=dashed](8.15,.05)(2.85,.05)
\pccurve[linecolor=pup, linestyle=dashed, angleA=180,angleB=0](2.85,.05)(2.05,-.95)
\psline[linecolor=pup, linestyle=dashed](2.05,-.95)(1.05,-.95)
\pccurve[linecolor=pup, linestyle=dashed, angleA=180,angleB=-90](1.05,-.95)(.05,.05)
\psline(0,0)(0,2)
\pccurve[angleA=90,angleB=180](0,2)(1,3)
\psline(1,3)(8.2,3)
\pccurve[angleA=0,angleB=90](8.2,3)(9.2,2)
\psline(9.2,2)(9.2,0)
\pccurve[angleA=-90,angleB=0](9.2,0)(8.2,-1)
\pccurve[angleA=180,angleB=-45](5.3,0)(5,.4)
\pccurve[angleA=135,angleB=0](4.8,.6)(4.5,1)
\pccurve[angleA=180,angleB=0](4.5,1)(3.7,0)
\pccurve[angleA=180,angleB=-45](3.7,0)(3.4,.4)
\pccurve[angleA=135,angleB=0](3.2,.6)(2.9,1)
\psline(2.9,1)(2.1,1)
\pccurve[angleA=180,angleB=45](2.1,1)(1.8,.6)
\pccurve[angleA=-135,angleB=0](1.6,.4)(1.3,0)
\psline(1.3,0)(1,0)
\pccurve[angleA=180,angleB=-90](1,0)(.5,.5)
\psline(.5,.5)(.5,2)
\pccurve[angleA=90,angleB=180](.5,2)(1,2.5)
\psline(1,2.5)(8.2,2.5)
\pccurve[angleA=0,angleB=90](8.2,2.5)(8.7,2)
\psline(8.7,2)(8.7,.5)
\pccurve[angleA=-90,angleB=0](8.7,.5)(8.2,0)
\psline(8.2,-1)(6.9,-1)
\pccurve[angleA=180,angleB=-45](6.9,-1)(6.6,-.6)
\pccurve[angleA=135,angleB=0](6.4,-.4)(6.1,0)
\psline(6.1,0)(5.3,0) \psline(8.2,0)(6.9,0)
\pccurve[angleA=180,angleB=0](6.9,0)(6.1,-1)
\psline(6.1,-1)(2.9,-1)
\pccurve[angleA=180,angleB=-45](2.9,-1)(2.6,-.6)
\pccurve[angleA=135,angleB=0](2.4,-.4)(2.1,0)
\pccurve[angleA=180,angleB=0](2.1,0)(1.3,1)
\pccurve[angleA=180,angleB=180](1.3,1)(1.3,2)
\psline(1.3,2)(7.9,2)
\pccurve[angleA=0,angleB=0](7.9,2)(7.9,1)
\psline(7.9,1)(5.3,1)
\pccurve[angleA=180,angleB=0](5.3,1)(4.5,0)
\pccurve[angleA=180,angleB=-45](4.5,0)(4.2,.4)
\pccurve[angleA=135,angleB=0](4,.6)(3.7,1)
\pccurve[angleA=180,angleB=0](3.7,1)(2.9,0)
\pccurve[angleA=180,angleB=0](2.9,0)(2.1,-1)
\psline(2.1,-1)(1,-1)
\pccurve[angleA=180,angleB=-90](1,-1)(0,0)
\psline[arrowscale=1.5]{->}(1.8,-1)(1.9,-1)
\rput(4.6,-1.6){{$(i)$}} \rput(4.6,1.5){{$D^2$}}
\end{pspicture}$$
$$
\begin{pspicture}[shift=-1.2](0,-2)(9.5,3.1)
\psline[linecolor=pup, linestyle=dashed](2.9,1)(6.9,1)
\psline[linecolor=pup, linestyle=dashed](2.9,0)(6.9,0)
\psline[linecolor=pup, linestyle=dashed](2.9,-1)(6.9,-1)
\psline[linewidth=1.5pt, linecolor=darkgray](4.9,0)(5.1,.4)(2,.4)(1.8,0)(1.6,1)(2,1.8)(5.3,1.8)(4.9,1)
\psline[linewidth=1.5pt, linecolor=darkgray](4.1,0)(4.4,.6)(2,.6)(1.7,0)(1.7,1)(2,1.6)(4.4,1.6)(4.1,1)
\psline[linewidth=1.5pt, linecolor=darkgray](3.3,0)(3.7,.8)(2,.8)(1.6,0)(1.8,1)(2,1.4)(3.5,1.4)(3.3,1)
\psline[linewidth=1.5pt, linecolor=darkgray](6.5,-1)(6.63,-.8)(2.63,-.8)(2.5,-1)(2.5,0)(2.63,.2)(6.63,.2)(6.5,0)
\psline(0,0)(0,2) \pccurve[angleA=90,angleB=180](0,2)(1,3)
\psline(1,3)(8.2,3) \pccurve[angleA=0,angleB=90](8.2,3)(9.2,2)
\psline(9.2,2)(9.2,0) \pccurve[angleA=-90,angleB=0](9.2,0)(8.2,-1)
\psline(8.2,-1)(6.9,-1) \psline(6.1,0)(5.3,0)
\pccurve[angleA=180,angleB=-45](5.3,0)(5,.4)
\pccurve[angleA=135,angleB=0](4.8,.6)(4.5,1)
\pccurve[angleA=180,angleB=0](4.5,1)(3.7,0)
\pccurve[angleA=180,angleB=-45](3.7,0)(3.4,.4)
\pccurve[angleA=135,angleB=0](3.2,.6)(2.9,1)
\psline(2.9,1)(2.1,1) \pccurve[angleA=180,angleB=45](2.1,1)(1.8,.6)
\pccurve[angleA=-135,angleB=0](1.6,.4)(1.3,0)
\psline(1.3,0)(1,0) \pccurve[angleA=180,angleB=-90](1,0)(.5,.5)
\psline(.5,.5)(.5,2) \pccurve[angleA=90,angleB=180](.5,2)(1,2.5)
\psline(1,2.5)(8.2,2.5) \pccurve[angleA=0,angleB=90](8.2,2.5)(8.7,2)
\psline(8.7,2)(8.7,.5) \pccurve[angleA=-90,angleB=0](8.7,.5)(8.2,0)
\psline(8.2,-1)(6.9,-1)
\pccurve[angleA=180,angleB=-45](6.9,-1)(6.6,-.6)
\pccurve[angleA=135,angleB=0](6.4,-.4)(6.1,0)
\psline(6.1,0)(5.3,0) \psline(8.2,0)(6.9,0)
\pccurve[angleA=180,angleB=0](6.9,0)(6.1,-1)
\psline(6.1,-1)(2.9,-1)
\pccurve[angleA=180,angleB=-45](2.9,-1)(2.6,-.6)
\pccurve[angleA=135,angleB=0](2.4,-.4)(2.1,0)
\pccurve[angleA=180,angleB=0](2.1,0)(1.3,1)
\pccurve[angleA=180,angleB=180](1.3,1)(1.3,2)
\psline(1.3,2)(7.9,2)
\pccurve[angleA=0,angleB=0](7.9,2)(7.9,1)
\psline(7.9,1)(5.3,1)
\pccurve[angleA=180,angleB=0](5.3,1)(4.5,0)
\pccurve[angleA=180,angleB=-45](4.5,0)(4.2,.4)
\pccurve[angleA=135,angleB=0](4,.6)(3.7,1)
\pccurve[angleA=180,angleB=0](3.7,1)(2.9,0)
\pccurve[angleA=180,angleB=0](2.9,0)(2.1,-1)
\psline(2.1,-1)(1,-1)
\pccurve[angleA=180,angleB=-90](1,-1)(0,0)
\psline[arrowscale=1.5]{->}(1.8,-1)(1.9,-1)
\rput(4.65,1.5){{$\alpha_4$}} \rput(3.7,1.3){{$\alpha_3$}}
\rput(5.6,1.7){{$\alpha_2$}} \rput(6.6,.4){{$\alpha_1$}}
\rput(4.6,-1.8){{$(ii)$}}
\end{pspicture}$$
\caption{$(i)$ The knot $5_2$ as a closed braid and $D^2$ for a basket surface $S(5_2)$
$(ii)$ the paths $\alpha_1$, $\alpha_2$, $\alpha_3$,
$\alpha_4$ on the disc $D^2$.} \label{52complete}
\end{figure}

Before we proceed to prove Theorem~\ref{baskettheorem1}, let us deal with an example.

\begin{exa}\label{bk52}
The basket number of the knot $5_2$ is less than or equal to $4$.
\begin{proof}
The knot $5_2$ is presented by a closed braid $\sigma_2
\sigma_1^{-1} (\sigma_2)^{-3} (\sigma_1)^{-1} \in B_3$ as
depicted in Figure~\ref{52complete} $(i)$ where the word $W=(\sigma_2)^{-3} (\sigma_1)^{-1}$.
To construct a basket surface, we first choose a $2$-disc $D^2$ as the gray dashed line
illustrated in Figure~\ref{52complete} $(i)$ which is the union of three discs
bounded by three Seifert circles and connected by the half
twisted bands presented by two leftmost braid generators $\sigma_2 \sigma_1^{-1}$.
The last letter in word $W$ is $(\sigma_1)^{-1}$
has the same sign compared to the generator $(\sigma_1)^{-1}$ used for
$D^2$. Thus, we perform an $A_{-1}$ plumbing along the arc $\alpha_1$
chosen on the disc $D^2$.
The next last letter in $W$ is $\sigma_2^{-1}$ and
it has the different sign compared to the generator $\sigma_2$ used for
$D^2$. Thus, we perform a $A_0$ plumbing along the arc
$\alpha_2$ in Figure~\ref{52complete} $(ii)$.
We continue $A_0$ plumbings  along the arcs $\alpha_3$ and $\alpha_4$
to obtain its canonical Seifert surface
$S$ as a basket surface of the knot $5_2$.
\end{proof}
\end{exa}

Let us remark that if $\sigma_i^s$ is chosen for the disc $D^2$,
every half twisted bands presented by $\sigma_i^{-s}$ can be obtained
by a flat plumbing while the half twisted bands presented by
$\sigma_i^{s}$ can be obtained by an $A_s$  plumbing. Using this idea,
we obtain the following theorem.

\begin{theorem}
Let $L$ be a link which is the closure of a braid $\beta \in B_n$
where the length of the braid $\beta$ is $m$.
Then the basket number of $L$ is less than or equal to $m-n+1$, $i. e.$,
$$bk(L)\le m-n+1.$$ \label{baskettheorem1}
\end{theorem}
\begin{proof}
Since $L$ is not splittable, for $i=1, 2, \ldots, n-1$,
its braid presentative $\beta\in B_n$ contains
at least one of $\sigma_i$ or $\sigma_i^{-1}$, say
$\sigma_i^{\epsilon_i}$ where $\epsilon_i \in \{1, -1\}$.
As we have seen in an example for the knot $5_2$,
we pick a $2$-disc $D^2$ which is obtained from $n$
disjoint discs bounded by $n$ Seifert circles
by attaching $(n-1)$ twisted bands presented by
$\sigma_{1}^{\epsilon_1}$, $\sigma_{2}^{\epsilon_2}$,
$\ldots$, $\sigma_{n-1}^{\epsilon_{n-1}}$.
For any other letter in $\beta$, it is fairly easy to see that
one can pick a $3$-cell $D_{\alpha_1}$ on
top of $D^2$ along $C_{\alpha_1}$ of $\alpha$ on $D^2$
which satisfies
$A_{\pm 1}(\mathrm{or}~ A_0) \cap \partial D_{\alpha_1}= C_{\alpha_1}$ as
depicted in Figure~\ref{52complete}.
Inductively, the canonical Seifert surface $S$ can
be obtained from the disc $D^2$ by plumbing
exactly $m-n+1$ times.
\end{proof}

Although, we have chosen $\sigma_{n-1}$ $\sigma_{n-2}$
$\ldots $ $\sigma_2$ $\sigma_1$ to construct
$D^2$, we might choose some of $\sigma_i$s to be inverses or a
different order as we have seen an example for the knot $5_2$.
The following example demonstrates that the
inequality in Theorem~\ref{baskettheorem1} is sharp.

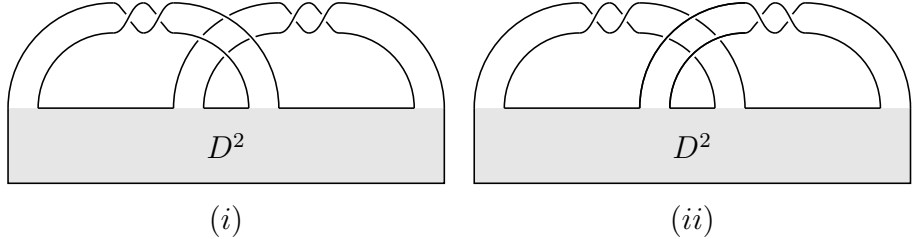
\begin{figure}
$$
\begin{pspicture}[shift=-1.4](-2.4,-1.7)(3.8,2)
\psframe[linecolor=lightgray,fillstyle=solid,fillcolor=lightgray](-2.2,-1)(3.6,0)
\psline(-2.2,0)(-2.2,-1)(3.6,-1)(3.6,0)
\psline(3.2,0)(1.4,0)
\psline(1,0)(.4,0)
\psline(0,0)(-1.8,0)
\pccurve[angleA=0,angleB=180](-.8,1)(-.4,1.4)
\pccurve[angleA=0,angleB=180](-.4,1)(0,1.4)
\pccurve[angleA=0,angleB=135](-.8,1.4)(-.63,1.25)
\pccurve[angleA=-45,angleB=180](-.57,1.15)(-.4,1)
\pccurve[angleA=0,angleB=135](-.4,1.4)(-.23,1.25)
\pccurve[angleA=-45,angleB=180](-.17,1.15)(0,1)
\pccurve[angleA=0,angleB=180](1.4,1)(1.8,1.4)
\pccurve[angleA=0,angleB=180](1.8,1)(2.2,1.4)
\pccurve[angleA=0,angleB=135](1.4,1.4)(1.57,1.25)
\pccurve[angleA=-45,angleB=180](1.63,1.15)(1.8,1)
\pccurve[angleA=0,angleB=135](1.8,1.4)(1.97,1.25)
\pccurve[angleA=-45,angleB=180](2.03,1.15)(2.2,1)
\psarc(0,0){1.4}{0}{90} \psarc(0,0){1}{0}{90}
\psarc(-.8,0){1.4}{90}{180} \psarc(-.8,0){1}{90}{180}
\psarc(2.2,0){1.4}{0}{90} \psarc(2.2,0){1}{0}{90}
\psarc(1.4,0){1.4}{90}{118} \psarc(1.4,0){1.4}{122}{136} \psarc(1.4,0){1.4}{140}{180}
\psarc(1.4,0){1}{90}{108} \psarc(1.4,0){1}{113}{132} \psarc(1.4,0){1}{137}{180}
\rput(.7,-.5){{$D^2$}}
\rput(.7,-1.5){{$(i)$}}
\end{pspicture}
\begin{pspicture}[shift=-1.4](-2.4,-1.7)(3.8,2)
\psframe[linecolor=lightgray,fillstyle=solid,fillcolor=lightgray](-2.2,-1)(3.6,0)
\psline(-2.2,0)(-2.2,-1)(3.6,-1)(3.6,0)
\psline(3.2,0)(1.4,0)
\psline(1,0)(.4,0)
\psline(0,0)(-1.8,0)
\pccurve[angleA=0,angleB=180](-.8,1)(-.4,1.4)
\pccurve[angleA=0,angleB=180](-.4,1)(0,1.4)
\pccurve[angleA=0,angleB=135](-.8,1.4)(-.63,1.25)
\pccurve[angleA=-45,angleB=180](-.57,1.15)(-.4,1)
\pccurve[angleA=0,angleB=135](-.4,1.4)(-.23,1.25)
\pccurve[angleA=-45,angleB=180](-.17,1.15)(0,1)
\pccurve[angleA=0,angleB=180](1.4,1.4)(1.8,1)
\pccurve[angleA=0,angleB=180](1.8,1.4)(2.2,1)
\pccurve[angleA=0,angleB=-135](1.4,1)(1.57,1.15)
\pccurve[angleA=45,angleB=180](1.63,1.25)(1.8,1.4)
\pccurve[angleA=0,angleB=-135](1.8,1)(1.97,1.15)
\pccurve[angleA=45,angleB=180](2.03,1.25)(2.2,1.4)
\psarc(0,0){1.4}{0}{40} \psarc(0,0){1.4}{44}{58} \psarc(0,0){1.4}{62}{90}
\psarc(0,0){1}{0}{43} \psarc(0,0){1}{48}{67} \psarc(0,0){1}{72}{90}
\psarc(-.8,0){1.4}{90}{180} \psarc(-.8,0){1}{90}{180}
\psarc(2.2,0){1.4}{0}{90} \psarc(2.2,0){1}{0}{90}
\psarc(1.4,0){1.4}{90}{180} \psarc(1.4,0){1}{90}{180}
\psarc(1.4,0){1.4}{90}{180} \psarc(1.4,0){1}{90}{180}
\rput(.7,-.5){{$D^2$}}
\rput(.7,-1.5){{$(ii)$}}
\end{pspicture}
$$
\caption{Basket surfaces of $(a)$ the trefoil knot and
$(b)$ the figure eight knot whose basket numbers are $2$~\cite{KKL:string}.} \label{2bandsfig}
\end{figure}

\begin{exa} \label{basketexa}
The basket number of the trefoil knot and figure eight knot is $2$.
\begin{proof}
Let $K_1$ be the trefoil knot which is presented by $\sigma_1 \sigma_1 \sigma_1 \in B_2$ and
let $K_2$ be the figure eight knot which is presented by $\sigma_2 \sigma_1 \sigma_2 \sigma_1 \in B_3$.
Then by applying Theorem~\ref{baskettheorem1}
for $W_1 =\sigma_1 \sigma_1$ for $K_1$ and $W_2 =\sigma_2 \sigma_1$ for $K_2$,
the basket numbers of these two knots are less than or equal to $2$. However,
by Remark~\ref{remark1} $(1)$, nontrivial knots can not have the basket
number zero and by Remark~\ref{remark1} $(2)$, the basket numbers of
these two knots must be even.
Therefore, the basket numbers of the trefoil knot and the figure eight knot are $2$.
Basket surfaces of these two knots are illustrated in Figure~\ref{2bandsfig}
which was found in~\cite{KKL:string}.
\end{proof}
\end{exa}

An upper bound for basket number of a link from its canonical Seifert surface will be
obtained in Corollary~\ref{bktheorem3} in subsection~\ref{flatbaskets}.

\subsection{Flat plumbing basket number} \label{flatbaskets}

Let $L$ be an oriented link which is a closed $n$-braid with a braid
word $\sigma_{n-1}\sigma_{n-2}\ldots\sigma_1 W$ where the length of
$W$ is $m$ and $W$ has $p$ positive letters.
Furihata, Hirasawa and Kobayashi~\cite{FHK:openbook} first found
an upper bound for the flat plumbing basket number of
$L$ in Theorem~\ref{falbktheorem1}.

\begin{theorem} (\cite{FHK:openbook})
Let $L$ be an oriented link which is a closed $n$-braid with a braid
word $\sigma_{n-1}\sigma_{n-2}\ldots\sigma_1 W$ where the length of
$W$ is $m$ and $W$ has $p$ positive letters, then there exists a
flat plumbing basket surface $S$ with $m + 2p$ bands such that
$\partial S$ is isotopic to $L$, $i.e.,$ $fpbk(L)\le m+2p$.
\label{falbktheorem1}
\end{theorem}

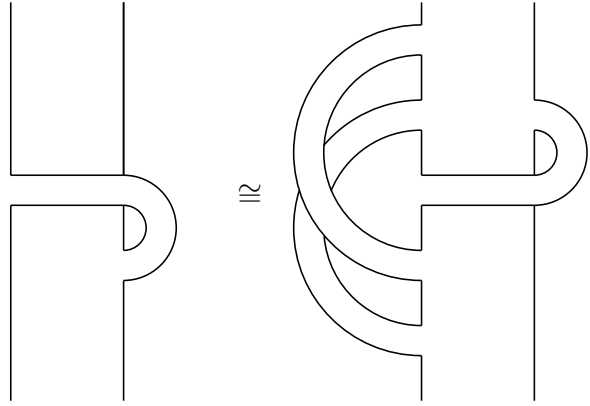
\begin{figure}
$$
\begin{pspicture}[shift=-2.8](-1,-3)(1.9,2.8)
\qline(.5,2.5)(.5,.2)
\qline(.5,2.5)(.5,.2)
\qline(.5,-.8)(.5,-.2) \qline(.5,-1.2)(.5,-2.8) \qline(-1,2.5)(-1,.2)
\qline(-1,-.2)(-1,-2.8) \psarc(.5,-.5){.3}{-90}{90}
\psarc(.5,-.5){.7}{-90}{90} \qline(-1,.2)(.5,.2) \qline(-1,-.2)(.5,-.2)
\end{pspicture} \cong
\begin{pspicture}[shift=-2.8](-3,-3)(1.5,2.8)
\qline(.5,2.5)(.5,1.2) \qline(.5,.8)(.5,.2) \qline(.5,-.2)(.5,-2.8)
\qline(-1,2.5)(-1,2.2) \qline(-1,1.8)(-1,1.2) \qline(-1,.8)(-1,.2)
\qline(-1,-.2)(-1,-.8) \qline(-1,-1.2)(-1,-1.8)
\qline(-1,-2.2)(-1,-2.8) \psarc(.5,.5){.3}{-90}{90}
\psarc(.5,.5){.7}{-90}{90} \qline(-1,.2)(.5,.2) \qline(-1,-.2)(.5,-.2)
\psarc(-1,.5){1.3}{90}{270} \psarc(-1,.5){1.7}{90}{270}
\psarc(-1,-.5){1.3}{90}{158} \psarc(-1,-.5){1.3}{185}{270}
\psarc(-1,-.5){1.7}{90}{140} \psarc(-1,-.5){1.7}{163}{270}
\end{pspicture}
$$
\caption{Changing the sign of a twisted band by two flat annuli plumbings.} \label{3annuli}
\end{figure}

The main idea of the theorem is as follows.
Let $D^2$ be a $2$-disc which is a union of discs bounded by
$n$ Seifert circles in the canonical Seifert surface of $L$
and half twisted bands presented by $\sigma_{n-1}\sigma_{n-2}\ldots\sigma_1$.
For each half twisted band presented by $\sigma_i^{-1}$ in $W$, it can be obtained
by an $A_0$ pluming along an arc in $D^2$.
For positive word $\sigma_i$ in $W$,
we first change $\sigma_i$ to $\sigma_i^{-1}$
by plumbing two flat plumbings as depicted in Figure~\ref{3annuli}.
Thus, the half twisted band presented by $\sigma_i$ can be obtained three
$A_0$ pluming along an arc in $D^2$.

To explain the idea of Corollary~\ref{falbkcor1},
let us deal with the knot $5_2$ in the following example.

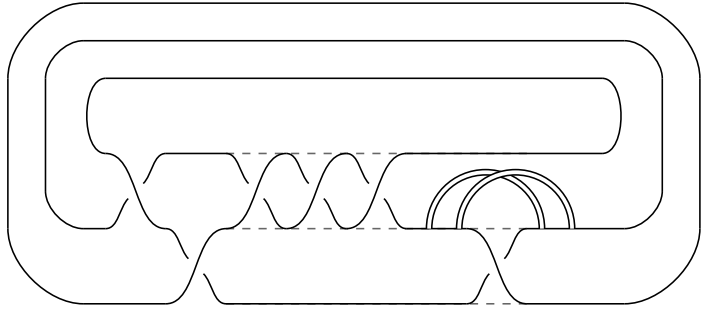
\begin{figure}
$$
\begin{pspicture}[shift=-1.2](0,-1.2)(9.5,3.1)
\psline[linecolor=pup, linestyle=dashed](2.9,1)(6.9,1)
\psline[linecolor=pup, linestyle=dashed](2.9,0)(6.9,0)
\psline[linecolor=pup, linestyle=dashed](2.9,-1)(6.9,-1)
\psline(0,0)(0,2)
\pccurve[angleA=90,angleB=180](0,2)(1,3)
\psline(1,3)(8.2,3)
\pccurve[angleA=0,angleB=90](8.2,3)(9.2,2)
\psline(9.2,2)(9.2,0)
\pccurve[angleA=-90,angleB=0](9.2,0)(8.2,-1)
\psline(8.2,-1)(6.9,-1)
\pccurve[angleA=180,angleB=0](6.9,-1)(6.1,0)
\psline(6.1,0)(5.3,0)
\pccurve[angleA=180,angleB=-45](5.3,0)(5,.4)
\pccurve[angleA=135,angleB=0](4.8,.6)(4.5,1)
\pccurve[angleA=180,angleB=0](4.5,1)(3.7,0)
\pccurve[angleA=180,angleB=-45](3.7,0)(3.4,.4)
\pccurve[angleA=135,angleB=0](3.2,.6)(2.9,1)
\psline(2.9,1)(2.1,1)
\pccurve[angleA=180,angleB=45](2.1,1)(1.8,.6)
\pccurve[angleA=-135,angleB=0](1.6,.4)(1.3,0)
\psline(1.3,0)(1,0)
\pccurve[angleA=180,angleB=-90](1,0)(.5,.5)
\psline(.5,.5)(.5,2)
\pccurve[angleA=90,angleB=180](.5,2)(1,2.5)
\psline(1,2.5)(8.2,2.5)
\pccurve[angleA=0,angleB=90](8.2,2.5)(8.7,2)
\psline(8.7,2)(8.7,.5)
\pccurve[angleA=-90,angleB=0](8.7,.5)(8.2,0)
\psline(8.2,0)(6.9,0)
\pccurve[angleA=180,angleB=45](6.9,0)(6.6,-.4)
\pccurve[angleA=-135,angleB=0](6.4,-.6)(6.1,-1)
\psline(6.1,-1)(2.9,-1)
\pccurve[angleA=180,angleB=-45](2.9,-1)(2.6,-.6)
\pccurve[angleA=135,angleB=0](2.4,-.4)(2.1,0)
\pccurve[angleA=180,angleB=0](2.1,0)(1.3,1)
\pccurve[angleA=180,angleB=180](1.3,1)(1.3,2)
\psline(1.3,2)(7.9,2)
\pccurve[angleA=0,angleB=0](7.9,2)(7.9,1)
\psline(7.9,1)(5.3,1)
\pccurve[angleA=180,angleB=0](5.3,1)(4.5,0)
\pccurve[angleA=180,angleB=-45](4.5,0)(4.2,.4)
\pccurve[angleA=135,angleB=0](4,.6)(3.7,1)
\pccurve[angleA=180,angleB=0](3.7,1)(2.9,0)
\pccurve[angleA=180,angleB=0](2.9,0)(2.1,-1)
\psline(2.1,-1)(1,-1)
\pccurve[angleA=180,angleB=-90](1,-1)(0,0)
\psarc[doubleline=true](6.35,0){.75}{0}{180}
\psarc[doubleline=true](6.75,0){.75}{0}{180}
\end{pspicture}$$
\caption{A flat plumbing basket surface of $5_2$.} \label{flatbasket52}
\end{figure}

\begin{exa}\label{fpbn52}
The flat plumbing basket number of the knot
$5_2$ is less than or equal to $6$.
\begin{proof}
A closed braid $\overline{\beta}$ form of the knot
$5_2$ is given in Figure~\ref{52complete} $(i)$ and
its braid word $\beta$ is $\sigma_2 \sigma_1^{-1} (\sigma_2)^{-3} (\sigma_1)^{-1} \in B_3$.
If we want to directly use Theorem~\ref{falbktheorem1},
we insert $\sigma_1 \sigma_1^{-1}$ to have
$\sigma_2 (\sigma_1 \sigma_1^{-1}) \sigma_1^{-1} (\sigma_2)^{-3} (\sigma_1)^{-1}$.
Beside $\sigma_2 \sigma_1$, all other letters are negative, $i.e.$, $p=0$.
Therefore, we find an upper bound $6$ by Theorem~\ref{falbktheorem1}.

However, we may choose a different $2$-disc
$D^2$ as depicted in Figure~\ref{52complete} $(i)$.
As we have seen in Example~\ref{bk52}, the half
twisted bands presented by $(\sigma_2)^{-3}$
can obtained by $A_0$ plumbing along the paths $\alpha_1$, $\alpha_2$ and $\alpha_3$
but the half twisted bands presented by the last word
$(\sigma_1)^{-1}$ was obtained by $A_1$ plumbing
along the path $\alpha_4$.
For the flat plumbing basket surface, we are only allowed
to use $A_0$ plumbing along the path in $D^2$.
By adding two $A_0$ plumbing, the sign of crossing
$(\sigma_1)^{-1}$ can be changed to $\sigma_1$ as in Figure~\ref{flatbasket52}.
Thus, we can build a flat plumbing basket surface of
the knot $5_2$ as depicted in Figure~\ref{flatbasket52}
by six $A_0$ plumbing along the paths
in the $2$-disc $D^2$.
\end{proof}
\end{exa}

\begin{figure}
$$
\begin{pspicture}[shift=-.8](-.7,-1.7)(6.2,2)
\psarc[doubleline=true](4,-.5){1.5}{-5}{185}
\psarc[doubleline=true](2.75,-.5){.75}{-5}{185}
\psarc[doubleline=true](2.25,-.5){.75}{-5}{185}
\psarc[doubleline=true](3,-.5){2}{-5}{185}
\psarc[doubleline=true](2.5,-.5){2}{-5}{185}
\psarc[doubleline=true](2,-.5){2}{-5}{185}
\psframe[linecolor=lightgray,fillstyle=solid,fillcolor=lightgray](-.5,-1.5)(6,-.5)
\psline(-.03,-.5)(-.5,-.5)(-.5,-1.5)(6,-1.5)(6,-.5)(5.53,-.5)
\psline(.03,-.5)(.47,-.5) \psline(.53,-.5)(.97,-.5) \psline(1.03,-.5)(1.47,-.5)
\psline(1.53,-.5)(1.97,-.5) \psline(2.03,-.5)(2.47,-.5) \psline(2.53,-.5)(2.97,-.5)
\psline(3.03,-.5)(3.47,-.5) \psline(3.53,-.5)(3.97,-.5) \psline(4.03,-.5)(4.47,-.5)
\psline(4.53,-.5)(4.97,-.5) \psline(5.03,-.5)(5.47,-.5)
\rput(2.75,-1){{$D^2$}}
\end{pspicture}
$$
\caption{A flat plumbing basket surface of $5_2$ in the trivial open book decomposition~\cite{CDK}.} \label{52fpbsopen}
\end{figure}
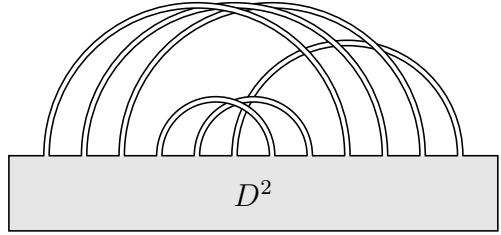

\begin{remark} \label{remark3}
Choi, Do and the author~\cite{CDK} recently proved that there
does not exist a nontrivial knot $K$ of the flat plumbing basket number $2$ and
a knot $K$ has the flat plumbing basket number $4$ if and only
of $K$ is either the trefoil knot or the figure eight knot.
By combining these results and Example~\ref{fpbn52}, we find
that the exact flat plumbing basket number of the knot $5_2$ is $6$.
A flat plumbing basket surface of $5_2$ in the trivial openbook decomposition
is give in Figure~\ref{52fpbsopen}~\cite{CDK}.
\end{remark}

As we have seen in Example~\ref{fpbn52},
as long as we pick only one of
$\sigma_i$, $\sigma_i^{-1}$ for $i=1, 2, \ldots, n-1$ to obtain the disc $D^2$,
Theorem~\ref{falbktheorem1}
works. So we obtained the following corollary.

\begin{cor} \label{falbkcor1}
Let $L$ be an oriented link which is a closed $n$-braid with a braid
word $\sigma_{n-1}^{\epsilon_{n-1}}\sigma_{n-2}^{\epsilon_{n-2}}
\ldots \sigma_1^{\epsilon_1}W$
where $\epsilon_i \in \{1, -1\}$ and the length of $W$ is $m$.
Let $s$ be the sum of the same sign letters;
$\sigma_{1}^{\epsilon_{1}}$, $\sigma_{2}^{\epsilon_{n-2}}$,
$\ldots$, $\sigma_{n-1}^{\epsilon_{n-1}}$ in $W$.
Then there exists a flat plumbing basket surface $S$ with $m + 2s$ bands such that
$\partial S$ is isotopic to $L$, $i.e.,$ $fpbk(L)\le m+2s$.
%
\end{cor}

The key ingredient of Theorem~\ref{falbktheorem1} is
that each crossing corresponding to the opposite sign letter can be obtained
by $A_0$ annulus plumbing as shown in Figure~\ref{52complete}. However,
crossings corresponding to the same sign letter can be obtained
by three $A_0$ annuli plumbings as shown in Figure~\ref{3annuli}.
However, depending on the $\beta \in B_n$ for the link $\overline{\beta}$,
an wise choice of $\epsilon_i$ can be made to obtain a minimum among the upper bounds
in Corollary~\ref{falbkcor1} for the
link $\overline{\beta}$ as described in Theorem~\ref{fpbktheoremext}.

\begin{theorem}
Let $L$ be an oriented link which is a closed $n$-braid with a braid
word $\beta$ whose length is $m$ and
let $ps(\sigma_i^{\pm 1})$ be the power sum
of $\sigma_i^{\pm 1}$ in $\beta$ for all $i=1, 2,
\ldots, n-1$.
Let $\gamma$ be the cardinality of the set
$$\Omega= \{ i | 1
\le i \le n-1, \sigma_i ~{\rm{and}}~ \sigma_i^{-1}~ {\rm{both}~\rm{appear}~ \rm{in}}~ \beta\}.$$
Let
$$\epsilon_i = \begin{cases} 1 ~~~& {\rm{if}}~ 1 \le
ps(\sigma_i^{1}) \le ps(\sigma_i^{- 1}) ~{\rm{or}}~ps(\sigma_i^{- 1})=0, \\
-1 ~~~&{\rm{if}}~ 1 \le ps(\sigma_i^{-1}) \le ps(\sigma_i^{1})
~{\rm{or}}~ps(\sigma_i^{1})=0. \end{cases}
$$
Then the flat plumbing basket number of $L$ is bounded by
$m + n-1 -4\gamma +2\sum_{i=1}^{n-1} ps(\sigma_i^{\epsilon_i})$, $i. e.,$
$$fpbk(L)\le m + n-1 -4\gamma +2\sum_{i=1}^{n-1} ps(\sigma_i^{\epsilon_i}).$$
\label{fpbktheoremext}
\begin{proof}
Since $L$ is nor splittable, for $i=1, 2, \ldots, n-1$,
its braid presentative $\beta\in B_n$ contains
at least one of $\sigma_i$ or $\sigma_i^{-1}$ by Remark~\ref{remark0} $(1)$.
For $1 \le i \le n-1$, we divide cases whether $i$ is in $\Omega$ or not.
For $i \in \Omega$, we chosen $\epsilon_i$ so that
$\sigma_i^{\epsilon_i}$ appears no more than $\sigma_i^{{-}\epsilon_i}$ does.
So it obviously contributes less or equal numbers for $s$ in Corollary~\ref{falbkcor1}.

For $i \not\in \Omega$, we chosen $\sigma_i^{\epsilon_i}$ to be the
one of $\sigma_i$ and $\sigma_i^{-1}$ which does not appears in $\beta$.
If we choose $\sigma_i^{-\epsilon_i}$ for the $2$-disc $D^2$,
then we will need $3\cdot (ps(\sigma_i^{-\epsilon_i})-1)$ times
flat plumbings as described in Figure~\ref{3annuli}.
On the other hand, since $L$ is prime, $\sigma_i^{-\epsilon_i}$
has to appear at least twice
in $\beta$ by Remark~\ref{remark0} $(2)$, $i.e.$, $ps(\sigma_i^{-\epsilon_i}) \ge 2$.
Instead of using $\sigma_i^{-\epsilon_i}$ for $D^2$,
we first perform a Reidemeister move II to insert $\sigma_i\sigma_i^{-1}$ in $\beta$.
Then we choose $\sigma_i^{\epsilon_i}$ for the $2$-disc $D^2$,
then we will need $ ps(\sigma_i^{-\epsilon_i})+1$ times
flat plumbings. It is easy to see that $ps(\sigma_i^{-\epsilon_i})+1
\le 3\cdot (ps(\sigma_i^{-\epsilon_i})-1)$ for all
$ps(\sigma_i^{-\epsilon_i}) \ge 2$.

Let $D^2$ be the $2$-disc the union of discs bounded by Seifert
circles of the closed braid $\overline{\beta}$
and connected by half twisted bands presented by
$\sigma_{1}^{\epsilon_{1}}$, $\sigma_{2}^{\epsilon_{n-2}}$,
$\ldots$, $\sigma_{n-1}^{\epsilon_{n-1}}$.
Then for the link diagram after inserting $\sigma_i\sigma_i^{-1}$
for all $i \not\in \Omega$,
its canonical Seifert surface is a flat plumbing basket surface
with
\begin{align*}
&\sum_{i\in \Omega} \left[ 3\cdot \left(ps(\sigma_i^{\epsilon_i})-1\right)
+ ps(\sigma_i^{-\epsilon_i}) \right]
+ \sum_{i\not\in \Omega} ps(\sigma_i^{-\epsilon_i})+1\\
&= \sum_{i=1}^{n-1} \left(ps(\sigma_i^{\epsilon_i}) +
ps(\sigma_i^{-\epsilon_i})\right) + 2\sum_{i\in \Omega} ps(\sigma_i^{\epsilon_i})
-3|\Omega| + (n-1-|\Omega|)\\
&=m+ (n-1)+2 \sum_{i=1}^{n-1} ps(\sigma_i^{\epsilon_i})-4\gamma
\end{align*}
many flat plumbings. It completes the proof of the theorem.
\end{proof}
\end{theorem}

Furthermore, the following example demonstrates that
the upper bound in Theorem~\ref{fpbktheoremext}
is shaper than one in Theorem~\ref{falbktheorem1}.

\begin{exa} \label{fpbsexa}
Let $\beta=\sigma_3 \sigma_2 \sigma_1 \sigma_3^{-1}\sigma_2  \sigma_1 \sigma_1
\sigma_2^{-1}\sigma_3 \sigma_2  \sigma_2 \sigma_1^{-1}\sigma_2  \sigma_3 \sigma_3 ~in ~B_4$.
Since $\beta$ is the exact shape as stated in Theorem~\ref{falbktheorem1},
we find an upper bound for the flat plumbing
basket number of $L = \bar{\beta}$
is $m+2s=12+2 \cdot 9=30$. To use Theorem~\ref{fpbktheoremext},
we can directly calculated
that for all $i=1,2,3$, $\epsilon_i =-1$ and $ps(\sigma_1)=3$,
$ps(\sigma_2)=5$, $ps(\sigma_3)=4$,
$ps(\sigma_1^{-1})=ps(\sigma_2^{-1})=ps(\sigma_3^{-1})=1$.
Thus, the upper bound in Theorem~\ref{fpbktheoremext} is
$15 +(4-1) +2 \cdot 3 - 4 \cdot 3 =12$.
\end{exa}

Now we want to find an upper bound for the flat plumbing basket
number of $L$ by using a canonical Seifert surface $S(L)$.
For a flat plumbing basket surface, the arc $\alpha$ we are plumbing along has
to be in the disc $D^2$ which was fixed from the beginning, and thus we
have to choose the disc $D^2$ carefully. The
Seifert graph $\Gamma(S)$ of the Seifert surface $S$ of a
closed braid is a path with multi-edges. Thus, there is
no ambiguity about the choice of a co-tree  edge alternating spanning tree for $\Gamma(S)$.
For an arbitrary canonical Seifert surface, it is completely different.

\begin{remark} \label{remark2}
For a canonical Seifert surface
of an arbitrary diagram of a link,
it is not trivial to find a co-tree  edge alternating spanning tree
in its Seifert graph and this is the
main reason that we have proven Theorem~\ref{maintheorem1}.
\end{remark}

\begin{figure}
$$
\begin{pspicture}[shift=-2](-.3,-.7)(2.3,3.5)
\psarc(.5,1){.5}{180}{0}
\psline(0,1)(1,3)
\pccurve[angleA=0,angleB=90](1,3)(2,1)
\pscircle[fillstyle=solid,fillcolor=darkgray](0,1){.075}
\pscircle[fillstyle=solid,fillcolor=darkgray](1,1){.075}
\pscircle[fillstyle=solid,fillcolor=darkgray](2,1){.075}
\pscircle[fillstyle=solid,fillcolor=darkgray](1,3){.075}
\rput[t](.5,3){$T$}
\rput[t](.45,2.5){$+$}\rput[t](.5,.8){$+$}
\rput[t](1.6,2){$+$}
\rput[t](-.3,1.25){$a$}\rput[t](1,1.45){$b$}
\rput[t](2.2,1.25){$c$}\rput[t](1,3.45){$d$}
\rput[t](1,-.2){$(i)$}
\end{pspicture}
\quad\quad
\begin{pspicture}[shift=-2](-.3,-.7)(2.3,3.5)
\psline[linewidth=1.5pt](0,1)(1,1)
\psline[linewidth=1.5pt](0,1)(1,3)
\psline[linewidth=1.5pt](1,3)(2,1)
\pscircle[fillstyle=solid,fillcolor=darkgray](0,1){.075}
\pscircle[fillstyle=solid,fillcolor=darkgray](1,1){.075}
\pscircle[fillstyle=solid,fillcolor=darkgray](2,1){.075}
\pscircle[fillstyle=solid,fillcolor=darkgray](1,3){.075}
\rput[t](0,3){${\overline T}$}
\rput[t](.45,2.5){$+$} \rput[t](.5,1.325){$-$} \rput[t](1.6,2.5){$-$}
\rput[t](-.3,1.25){$a$}\rput[t](1,1.45){$b$}
\rput[t](2.2,1.25){$c$}\rput[t](1,3.45){$d$}
\rput[t](1,-.2){$(ii)$}
\end{pspicture}
\quad\quad
\begin{pspicture}[shift=-2](-.3,-.7)(3.3,4.1)
\psline[linecolor=darkgray,linewidth=1pt](2.5,.5)(2.85,.5)(2.85,3.65)(1.4,3.65)(1.1,3.35)(.1,3.65)(.1,1.225)(1.6,1.475)(1.75,1.475)(1.75,.5)(2,.5)
\psline[linecolor=darkgray,linewidth=1pt](2.5,1.375)(2.65,1.5)(2.65,3.5)(.2,3.5)(.2,1.375)(2,1.375)
\psline[linecolor=darkgray,linewidth=1pt](1.75,3)(1.75,3.35)(1.4,3.35)(1.1,3.65)(.3,3.35)(.3,1.525)(1.6,1.275)(1.75,1.75)(1.75,2)
\psline(0,.25)(0,3.75)
\psarc(.25,3.75){0.25}{90}{180} \psline(.25,4)(.4,4)
\pccurve[angleA=0,angleB=120](.4,4)(.7,3.6)
\pccurve[angleA=-60,angleB=180](.8,3.4)(1,3)
\pccurve[angleA=0,angleB=-110](1,3)(1.227,3.35)
\pccurve[angleA=70,angleB=180](1.3,3.6)(1.5,4) \psline(1.5,4)(2.75,4)
\psarc(2.75,3.75){.25}{0}{90} \psline(3,3.75)(3,.25)
\psarc(2.75,.25){.25}{180}{0} \psline(2.5,.25)(2.5,2.75)
\psarc(2.25,2.75){.25}{0}{90} \psline(2.25,3)(1.5,3)
\pccurve[angleA=0,angleB=180](1,4)(1.5,3)
\pccurve[angleA=30,angleB=180](.5,3)(1,4)
\pccurve[angleA=90,angleB=-150](.75,1.75)(.5,3)
\pccurve[angleA=140,angleB=-90](1.125,1.375)(.75,1.75)
\pccurve[angleA=90,angleB=-40](1.5,1)(1.125,1.375)
\psline(1.5,.25)(1.5,1)
\psarc(1.75,.25){.25}{180}{0}
\psline(2,1.75)(2,.25)
\psarc(1.75,1.75){.25}{0}{180}
\pccurve[angleA=45,angleB=-90](1.175,1.425)(1.5,1.75)
\pccurve[angleA=90,angleB=-135](.75,1)(1.075,1.325)
\psline(.75,.25)(.75,1)
\psarc(.5,.25){.25}{270}{0}
\psline(.5,0)(.25,0)
\psarc(.25,.25){.25}{180}{270}
\rput[t](2.25,.6){$\alpha_1$} \rput[t](2.25,1.475){$\alpha_2$} \rput[t](1.75,2.9){$\alpha_3$}
\rput[t](1.5,-.2){$(iii)$} \rput[t](.4,.7){$D^2$}
\end{pspicture}
$$
$$
\begin{pspicture}[shift=-2](-.3,-.7)(2.3,4.1)
\psline[linewidth=1.5pt](0,1)(1,1)
\psline(1,1)(2,1)
\psline[linewidth=1.5pt](1,3)(2,1)
\psarc(1.5,1){.5}{180}{0}
\psline[linewidth=1.5pt](0,1)(1,3)
\pccurve[angleA=70,angleB=110](1,1)(2,1)
\pscircle[fillstyle=solid,fillcolor=darkgray](0,1){.075}
\pscircle[fillstyle=solid,fillcolor=darkgray](1,1){.075}
\pscircle[fillstyle=solid,fillcolor=darkgray](2,1){.075}
\pscircle[fillstyle=solid,fillcolor=darkgray](1,3){.075}
\rput[t](1.5,1.55){$+$} \rput[t](1.5,.45){$+$} \rput[t](1.5,.9){$+$}
\rput[t](.35,2.3){$+$} \rput[t](1.6,2.5){$-$}
\rput[t](.5,1.325){$-$} \rput[t](-.3,1.25){$a$}
\rput[t](1,1.45){$b$} \rput[t](2.2,1.25){$c$}
\rput[t](1,3.45){$d$} \rput[t](1,-.2){$(iv)$}
\rput[t](2,3.3){$\Gamma(S_3)$}
\end{pspicture}
\quad\quad
\begin{pspicture}[shift=-2](-.3,-.7)(3.3,4.1)
\psline[linecolor=darkgray,linewidth=1pt](.75,.4)(.2,.4)(.2,1.525)(1.65,1.275)(1.65,.4)(1.5,.4)
\psline[linecolor=darkgray,linewidth=1pt](.75,.8)(.6,.8)(.6,1.275)(1.85,1.475)(1.85,.8)(1.5,.8)
\psline[linecolor=darkgray,linewidth=1pt](1.17,3.85)(.9,3.85)(.9,3.5)(1.6,3.5)(1.6,3.85)(1.37,3.85)
\psline(0,.25)(0,3.75)
\psarc(.25,3.75){0.25}{90}{180} \psline(.25,4)(.4,4)
\pccurve[angleA=0,angleB=120](.4,4)(.7,3.6)
\pccurve[angleA=-60,angleB=180](.8,3.4)(1,3)
\pccurve[angleA=0,angleB=-110](1,3)(1.227,3.35)
\pccurve[angleA=70,angleB=180](1.3,3.6)(1.5,4) \psline(1.5,4)(2.75,4)
\psarc(2.75,3.75){.25}{0}{90} \psline(3,3.75)(3,.25)
\psarc(2.75,.25){.25}{180}{0}
\pccurve[angleA=90,angleB=-40](2.5,.25)(2.3,.45)
\pccurve[angleA=140,angleB=-90](2.2,.55)(2,.75)
\pccurve[angleA=90,angleB=-90](2,.25)(2.5,.75)
\psline(2.5,.75)(2.5,1.25) \psline(2,.75)(2,1.25)
\pccurve[angleA=90,angleB=-40](2.5,1.25)(2.3,1.45)
\pccurve[angleA=140,angleB=-90](2.2,1.55)(2,1.75)
\pccurve[angleA=90,angleB=-90](2,1.25)(2.5,1.75)
\psline(2.5,1.75)(2.5,2.75)
\psarc(2.25,2.75){.25}{0}{90} \psline(2.25,3)(2,3)
\pccurve[angleA=0,angleB=180](1,4)(1.5,3)
\pccurve[angleA=30,angleB=180](.5,3)(1,4)
\pccurve[angleA=90,angleB=-150](.75,1.75)(.5,3)
\pccurve[angleA=140,angleB=-90](1.125,1.375)(.75,1.75)
\pccurve[angleA=90,angleB=-40](1.5,1)(1.125,1.375)
\psline(1.5,.25)(1.5,1)
\psarc(1.75,.25){.25}{180}{0}
\pccurve[angleA=90,angleB=-40](2,1.75)(1.75,2.375)
\pccurve[angleA=140,angleB=0](1.75,2.375)(1.5,3)
\pccurve[angleA=90,angleB=-120](1.5,1.75)(1.7,2.25)
\pccurve[angleA=60,angleB=180](1.8,2.5)(2,3)
\pccurve[angleA=45,angleB=-90](1.175,1.425)(1.5,1.75)
\pccurve[angleA=90,angleB=-135](.75,1)(1.075,1.325)
\psline(.75,.25)(.75,1)
\psarc(.5,.25){.25}{270}{0}
\psline(.5,0)(.25,0)
\psarc(.25,.25){.25}{180}{270} \rput[t](1.8,3.4){$\alpha_6$}
\rput[t](1.125,.5){$\alpha_4$} \rput[t](1.125,.9){$\alpha_5$} 
\rput[t](1.5,-.2){$(v)$} \rput[t](2.3,3.7){$S_3$}
\end{pspicture} \quad\quad
\begin{pspicture}[shift=-2](-.3,-.7)(2.3,4.1)
\psline[linewidth=1.5pt](0,1)(1,1)
\psline(1,1)(2,1) \pccurve[angleA=30,angleB=30](1,3)(2,1)
\psline[linewidth=1.5pt](1,3)(2,1)
\psarc(.5,1){.5}{180}{0}
\psarc(1.5,1){.5}{180}{0}
\psline[linewidth=1.5pt](0,1)(1,3)
\pccurve[angleA=70,angleB=110](1,1)(2,1)
\pccurve[angleA=-90,angleB=135](0,1)(.2,.2)
\pccurve[angleA=-45,angleB=-135](.2,.2)(.8,.2)
\pccurve[angleA=45,angleB=-90](.8,.2)(1,1)
\pscircle[fillstyle=solid,fillcolor=darkgray](0,1){.075}
\pscircle[fillstyle=solid,fillcolor=darkgray](1,1){.075}
\pscircle[fillstyle=solid,fillcolor=darkgray](2,1){.075}
\pscircle[fillstyle=solid,fillcolor=darkgray](1,3){.075}
\rput[t](.5,.45){$+$} \rput[t](.5,.05){$+$}
\rput[t](1.5,1.55){$+$} \rput[t](1.5,.45){$+$} \rput[t](1.5,.9){$+$}
\rput[t](.35,2.3){$+$} \rput[t](1.25,2.2){$-$} \rput[t](1.6,2.6){$+$}
\rput[t](.5,1.325){$-$} \rput[t](-.3,1.25){$a$}\rput[t](1,1.45){$b$}
\rput[t](2.2,1.25){$c$}\rput[t](1,3.45){$d$}
\rput[t](1,-.2){$(vi)$} \rput[t](2,3.3){$\Gamma(S_6)$}
\end{pspicture}
$$
$$
\begin{pspicture}[shift=-2](-.3,-.7)(3.3,4.1)
\psline[linecolor=darkgray,linewidth=1pt](1,3.5)(.85,3.5)(.95,3.13)(1.55,3.13)(1.65,3.5)(1.5,3.5)
\psline[linecolor=darkgray,linewidth=1pt](.75,1.24)(.5,1.24)(.5,1.5)(1.75,1.5)(1.75,1.24)(1.5,1.24)
\psline(0,.25)(0,3.75)
\psarc(.25,3.75){0.25}{90}{180}\psline(.25,4)(.4,4)
\pccurve[angleA=0,angleB=120](.4,4)(.7,3.6)
\pccurve[angleA=-60,angleB=180](.8,3.4)(1,3)
\pccurve[angleA=0,angleB=-150](1,3)(1.2,3.1)
\pccurve[angleA=30,angleB=-90](1.3,3.15)(1.5,3.3)
\psline(1.5,3.3)(1.5,3.7)
\pccurve[angleA=90,angleB=180](1,3.7)(1.5,4)
\pccurve[angleA=0,angleB=150](1,4)(1.15,3.95)
\pccurve[angleA=-30,angleB=90](1.25,3.9)(1.5,3.7)
\psline(1,3.7)(1,3.3)
\pccurve[angleA=-90,angleB=150](1,3.3)(1.25,3.13)
\pccurve[angleA=-30,angleB=180](1.25,3.13)(1.5,3)
\psline(1.5,4)(2.75,4)
\pccurve[angleA=30,angleB=180](.5,3)(1,4)
\pccurve[angleA=90,angleB=-150](.75,1.75)(.5,3)
\psarc(2.75,3.75){.25}{0}{90} \psline(3,3.75)(3,.25)
\psarc(2.75,.25){.25}{180}{0}
\pccurve[angleA=90,angleB=-40](2.5,.25)(2.3,.45)
\pccurve[angleA=140,angleB=-90](2.2,.55)(2,.75)
\pccurve[angleA=90,angleB=-90](2,.25)(2.5,.75)
\psline(2.5,.75)(2.5,1.25) \psline(2,.75)(2,1.25)
\pccurve[angleA=90,angleB=-40](2.5,1.25)(2.3,1.45)
\pccurve[angleA=140,angleB=-90](2.2,1.55)(2,1.75)
\pccurve[angleA=90,angleB=-90](2,1.25)(2.5,1.75)
\psline(2.5,1.75)(2.5,2.75)
\psarc(2.25,2.75){.25}{0}{90} \psline(2.25,3)(2,3)
\psarc(1.75,.25){.25}{180}{0}
\pccurve[angleA=90,angleB=-40](2,1.75)(1.75,2.375)
\pccurve[angleA=140,angleB=0](1.75,2.375)(1.5,3)
\pccurve[angleA=90,angleB=-120](1.5,1.75)(1.7,2.25)
\pccurve[angleA=60,angleB=180](1.8,2.5)(2,3)
\pccurve[angleA=90,angleB=-150](.75,.25)(1.125,.5)
\pccurve[angleA=30,angleB=-90](1.125,.5)(1.5,.75)
\pccurve[angleA=90,angleB=-30](1.5,.25)(1.175,.47)
\pccurve[angleA=150,angleB=-90](1.075,.53)(.75,.75)
\pccurve[angleA=90,angleB=-150](.75,.75)(1.125,1)
\pccurve[angleA=30,angleB=-90](1.125,1)(1.5,1.25)
\pccurve[angleA=90,angleB=-30](1.5,.75)(1.175,.97)
\pccurve[angleA=150,angleB=-90](1.075,1.03)(.75,1.25)
\pccurve[angleA=90,angleB=-150](.75,1.25)(1.075,1.47)
\pccurve[angleA=30,angleB=-90](1.175,1.53)(1.5,1.75)
\pccurve[angleA=90,angleB=-30](1.5,1.25)(1.125,1.5)
\pccurve[angleA=150,angleB=-90](1.125,1.5)(.75,1.75)
\psarc(.5,.25){.25}{270}{0}
\psline(.5,0)(.25,0)
\psarc(.25,.25){.25}{180}{270}
\rput[t](.25,1.4){$\alpha_7$} \rput[t](1.9,3.4){$\alpha_8$} 
\rput[t](1.5,-.4){$(vii)$} \rput[t](2.3,3.7){$S_6$}
\end{pspicture}
\quad\quad
\begin{pspicture}[shift=-2](-.3,-.7)(2.3,4.1)
\psline[linewidth=1.5pt](0,1)(1,1)
\pccurve[angleA=-70,angleB=-110, linewidth=1.5pt, linestyle=dashed](0,1)(1,1)
\psline(1,1)(2,1)
\psline[linewidth=1.5pt](0,1)(1,3)
\psline[linewidth=1.5pt](1,3)(2,1)
\psarc(.5,1){.5}{180}{0}
\psarc(1.5,1){.5}{180}{0}
\pccurve[angleA=-30,angleB=90, linewidth=1.5pt, linestyle=dashed](1,3)(2,1)
\pccurve[angleA=30,angleB=30](1,3)(2,1)
\pccurve[angleA=70,angleB=110](1,1)(2,1)
\pccurve[angleA=-90,angleB=135](0,1)(.2,.2)
\pccurve[angleA=-45,angleB=-135](.2,.2)(.8,.2)
\pccurve[angleA=45,angleB=-90](.8,.2)(1,1)
\pscircle[fillstyle=solid,fillcolor=darkgray](0,1){.075}
\pscircle[fillstyle=solid,fillcolor=darkgray](1,1){.075}
\pscircle[fillstyle=solid,fillcolor=darkgray](2,1){.075}
\pscircle[fillstyle=solid,fillcolor=darkgray](1,3){.075}
\rput[t](.5,.45){$+$} \rput[t](.5,.05){$+$} \rput[t](.5,1.35){$-$}
\rput[t](1.5,.45){$+$}\rput[t](1.5,.9){$+$}\rput[t](1.5,1.55){$+$}
\rput[t](.48,2.5){$+$} \rput[t](2.1,2.65){$+$}
\rput[t](1.25,2.2){$-$} \rput[t](-.3,1.25){$a$}\rput[t](1,1.45){$b$}
\rput[t](2.2,.85){$c$}\rput[t](1,3.45){$d$}
\rput[t](1,-.4){$(vii)$} \rput[t](0,3){$\Gamma(S_8)$}
\end{pspicture}
\quad\quad
\begin{pspicture}[shift=-2](-.3,-.7)(3.3,4.1)
\psline(0,.25)(0,3.75)
\psarc(.25,3.75){0.25}{90}{180}\psline(.25,4)(.4,4)
\pccurve[angleA=0,angleB=120](.4,4)(.7,3.6)
\pccurve[angleA=-60,angleB=180](.8,3.4)(1,3)
\pccurve[angleA=0,angleB=-150](1,3)(1.2,3.1)
\pccurve[angleA=30,angleB=-90](1.3,3.15)(1.5,3.3)
\pccurve[angleA=90,angleB=-30](1.5,3.3)(1.3,3.47)
\pccurve[angleA=150,angleB=-90](1.2,3.53)(1,3.7)
\pccurve[angleA=90,angleB=180](1,3.7)(1.5,4)
\pccurve[angleA=0,angleB=150](1,4)(1.15,3.95)
\pccurve[angleA=-30,angleB=90](1.25,3.9)(1.5,3.7)
\pccurve[angleA=-90,angleB=30](1.5,3.7)(1.25,3.5)
\pccurve[angleA=-150,angleB=90](1.25,3.5)(1,3.3)
\pccurve[angleA=-90,angleB=150](1,3.3)(1.25,3.13)
\pccurve[angleA=-30,angleB=180](1.25,3.13)(1.5,3)
\psline(1.5,4)(2.75,4)
\pccurve[angleA=30,angleB=180](.5,3)(1,4)
\pccurve[angleA=90,angleB=-150](.75,1.85)(.5,3)
\psarc(2.75,3.75){.25}{0}{90} \psline(3,3.75)(3,.25)
\psarc(2.75,.25){.25}{180}{0}
\pccurve[angleA=90,angleB=-40](2.5,.25)(2.3,.45)
\pccurve[angleA=140,angleB=-90](2.2,.55)(2,.75)
\pccurve[angleA=90,angleB=-90](2,.25)(2.5,.75)
\psline(2.5,.75)(2.5,1.25) \psline(2,.75)(2,1.25)
\pccurve[angleA=90,angleB=-40](2.5,1.25)(2.3,1.45)
\pccurve[angleA=140,angleB=-90](2.2,1.55)(2,1.75)
\pccurve[angleA=90,angleB=-90](2,1.25)(2.5,1.75)
\psline(2.5,1.75)(2.5,2.75)
\psarc(2.25,2.75){.25}{0}{90} \psline(2.25,3)(2,3)
\psarc(1.75,.25){.25}{180}{0}
\pccurve[angleA=90,angleB=-40](2,1.75)(1.75,2.375)
\pccurve[angleA=140,angleB=0](1.75,2.375)(1.5,3)
\pccurve[angleA=90,angleB=-120](1.5,1.85)(1.7,2.25)
\pccurve[angleA=60,angleB=180](1.8,2.5)(2,3)
\pccurve[angleA=90,angleB=-150](.75,.25)(1.125,.45)
\pccurve[angleA=30,angleB=-90](1.125,.45)(1.5,.65)
\pccurve[angleA=90,angleB=-30](1.5,.25)(1.175,.425)
\pccurve[angleA=150,angleB=-90](1.075,.475)(.75,.65)
\pccurve[angleA=90,angleB=-150](.75,.65)(1.125,.85)
\pccurve[angleA=30,angleB=-90](1.125,.85)(1.5,1.05)
\pccurve[angleA=90,angleB=-30](1.5,.65)(1.175,.825)
\pccurve[angleA=150,angleB=-90](1.075,.875)(.75,1.05)
\pccurve[angleA=90,angleB=-150](.75,1.05)(1.125,1.25)
\pccurve[angleA=30,angleB=-90](1.125,1.25)(1.5,1.45)
\pccurve[angleA=90,angleB=-30](1.5,1.05)(1.175,1.225)
\pccurve[angleA=150,angleB=-90](1.075,1.275)(.75,1.45)
\pccurve[angleA=90,angleB=-150](.75,1.45)(1.075,1.625)
\pccurve[angleA=30,angleB=-90](1.175,1.675)(1.5,1.85)
\pccurve[angleA=90,angleB=-30](1.5,1.45)(1.125,1.65)
\pccurve[angleA=150,angleB=-90](1.125,1.65)(.75,1.85)
\psarc(.5,.25){.25}{270}{0}
\psline(.5,0)(.25,0)
\psarc(.25,.25){.25}{180}{270}
\rput[t](1.5,-.4){$(ix)$} \rput[t](2.3,3.7){$\overline{S}(7_5)$}
\end{pspicture}
$$
\caption{The process to obtain a flat plumbing basket surface of the knot $7_5$.} \label{graph2}
\end{figure}
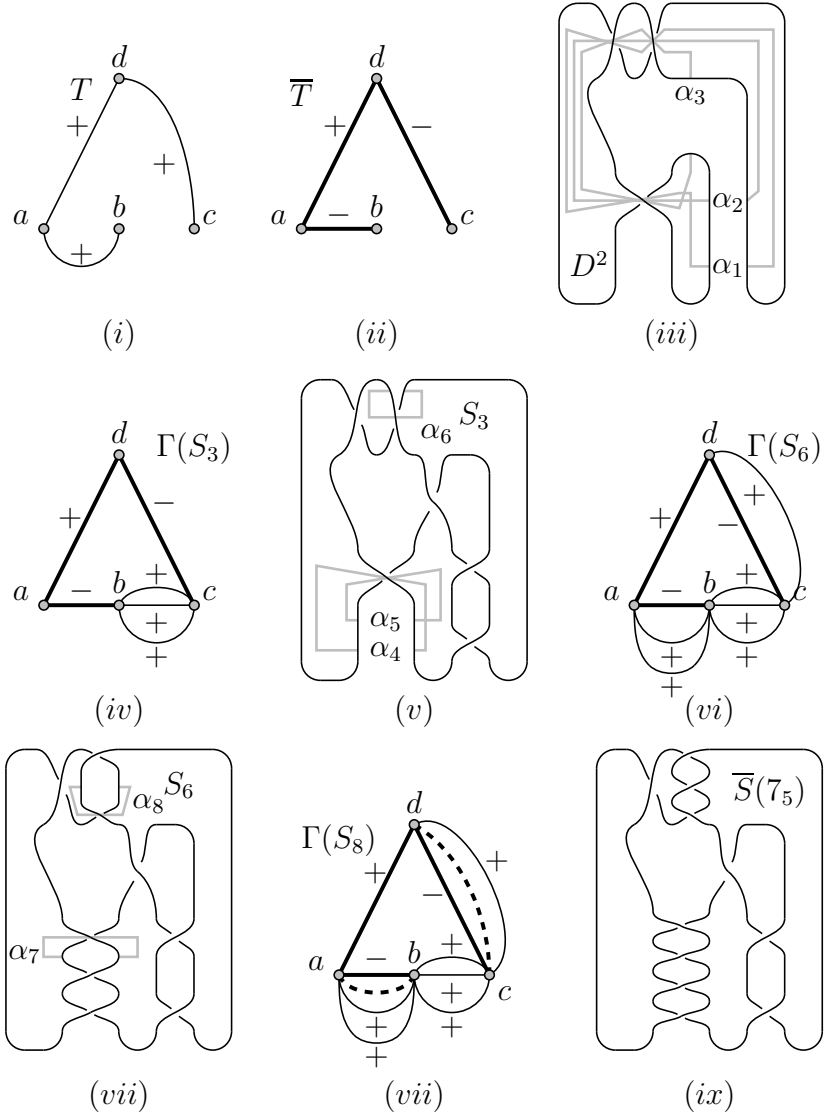

For example, let us consider the knot $7_5$ as depicted in Figure~\ref{SeifertS}.
Every spanning tree of the Seifert graph of $\Gamma(7_5)$ in Figure~\ref{graph} $(ii)$
is not a co-tree  edge alternating spanning tree because there
exist two consecutive edges of the same sign which will produce a
full twist that is prohibited in flat plumbing for every spanning tree.
Thus, none of these spanning trees
can be used directly to construct a flat plumbing basket surface in
neither Theorem~\ref{falbktheorem1} nor Theorem~\ref{fpbktheoremext}.
Before we proceed to find any upper bound from canonical Seifert surfaces,
let us deal with the knot $7_5$ in Example~\ref{52flat} to explain
the Remark~\ref{remark-1} and the key idea of Theorem~\ref{flatbktheorem3}.

\begin{figure}
$$
\begin{pspicture}[shift=-.8](-.7,-2.2)(8.2,2.8)
\psarc[doubleline=true](3.5,-.5){3}{-5}{185}
\psarc[doubleline=true](3,-.5){3}{-5}{185}
\psarc[doubleline=true](5.5,-.5){2}{-5}{185}
\psarc[doubleline=true](5,-.5){2}{-5}{185}
\psarc[doubleline=true](4,-.5){1.5}{-5}{185}
\psarc[doubleline=true](3.5,-.5){1.5}{-5}{185}
\psarc[doubleline=true](3,-.5){1.5}{-5}{185}
\psarc[doubleline=true](2.5,-.5){1.5}{-5}{185}
\psframe[linecolor=lightgray,fillstyle=solid,fillcolor=lightgray](-.5,-1.5)(8,-.5)
\psline(-.03,-.5)(-.5,-.5)(-.5,-1.5)(8,-1.5)(8,-.5)(7.53,-.5)
\psline(.03,-.5)(.47,-.5) \psline(.53,-.5)(.97,-.5) \psline(1.03,-.5)(1.47,-.5)
\psline(1.53,-.5)(1.97,-.5) \psline(2.03,-.5)(2.47,-.5) \psline(2.53,-.5)(2.97,-.5)
\psline(3.03,-.5)(3.47,-.5) \psline(3.53,-.5)(3.97,-.5) \psline(4.03,-.5)(4.47,-.5)
\psline(4.53,-.5)(4.97,-.5) \psline(5.03,-.5)(5.47,-.5)
\psline(5.53,-.5)(5.97,-.5) \psline(6.03,-.5)(6.47,-.5)
\psline(6.53,-.5)(6.97,-.5) \psline(7.03,-.5)(7.47,-.5)
\rput(3.75,-1){{$\mathcal{D}$}}
\rput(3.75,-1.9){{$(e)$}}
\end{pspicture}
$$
\caption{A flat plumbing basket surface of $7_5$
in the trivial open book decomposition~\cite{CDK}.} \label{75open}
\end{figure}
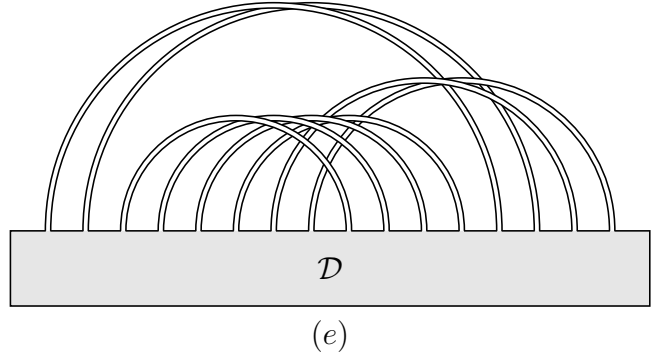

\begin{exa} \label{52flat}
A reduced braid presentation $\beta$ of $7_5$
is $\sigma_1 \sigma_2 \sigma_1^{-1} \sigma_2^2 \sigma_1^3 \in B_3$, thus, we
find that the upper bound of the flat plumbing basket number for $7_5$ by applying
Theorem~\ref{falbktheorem1} is $19$ and $10$ by Theorem~\ref{fpbktheoremext}.
The canonical Seifert surface of the diagram
in Figure~\ref{graph2} $(ix)$ of the knot $7_5$ is a flat plumbing basket surface.
Thus, we find
the upper bound of the flat plumbing basket number for the knot $7_5$ is $8$.
In fact, using the result which find
all knots of the flat plumbing basket number $6$ in~\cite{CDK},
the exact flat plumbing basket number of the knot $7_5$ is $8$.
A flat plumbing basket surface of $7_5$ in the trivial open book decomposition is
illustrated in Figure~\ref{75open}.
\begin{proof}
For all spanning trees of $\Gamma(S(7_5))$ in Figure~\ref{graph} $(ii)$,
there exist two consecutive edges of the same sign which will produce a
full twist which obstructs the existence of a
$3$-ball $B_{\alpha}$ in a flat plumbing.
Therefore, Theorem~\ref{fptheorem1} cannot be directly applied to construct
a flat plumbing basket surface of $7_5$ as described in Remark~\ref{remark2} (1).

We first consider $G(\Gamma)$ which is a simple graph
obtained from $\Gamma$ by identifying edges in the same parallel class
which was defined as the induced graph of $\Gamma$.
For this case,  $G(\Gamma)$ is the cycle graph $C_4$ on four vertices.
Instead of a given spanning tree $T$ in Figure~\ref{graph2} $(i)$,
we choose a co-tree edge alternating spanning tree $\overline{T}$
of $G(\Gamma)$ given by Theorem~\ref{maintheorem1}
as illustrated in Figure~\ref{graph2} $(ii)$.
We choose the disc $D^2$ which is union of four Seifert discs
and three half twisted bands corresponding to $\overline{T}$ as depicted in Figure~\ref{graph2} $(iii)$.

Now, we perform flat plumbings in the following order.
We first do three flat plumbings along arcs in $D^2$
as given in Figure~\ref{graph2} $(iii)$ to obtain a flat plumbing basket surface $S_3$
as depicted in Figure~\ref{graph2} $(v)$. Let us remark that to guarantee
the existence of $3$-ball as described in Figure\ref{fpbkprooffig},
we produce a half twisted band presented by
the three edges $\{b, c\}$ by flat plumbing along the arcs $\alpha_1$,
$\alpha_2$ and $\alpha_3$ (exactly in this order).
Then all remaining edges in $E(\Gamma(S(7_5)))\setminus E(\overline{T})$ can be obtained by
flat plumbings along the cars $\alpha_4$, $\alpha_5$ and $\alpha_6$ as illustrated in Figure~\ref{graph2} $(v)$.

Since two edges of `$-$' sign in the co-tree edge alternating spanning tree $\overline T$
were not in $\Gamma(S(7_5))$,
if we rebuild a Seifert surface from the Seifert graph $\Gamma(S_6)$ as depicted in Figure~\ref{graph2} $(vii)$,
then it is a flat plumbing basket surface but its boundary is not isotopic to the original link $7_5$.
However, adding two adjacent edges of different signs
is, in fact, a Reidemeister type II move which guarantees the resulting link is isotopic to
the original one as illustrated in Figure~\ref{Reide2}.
So at last, we add three flat plumbings along the
arcs $\alpha_7$ and $\alpha_8$ as presented in Figure~\ref{graph2} $(vii)$
to obtain the flat plumbing basket surface $S_8$ which is the desired surface $\overline{S}(7_5)$.
The canonical Seifert surface of the knot diagram in Figure~\ref{graph2} $(ix)$
is indeed a flat plumbing basket surface of the knot $7_5$ and it does require $8$ flat plumbings.
\end{proof}
\end{exa}

Let us deal with general cases for the flat plumbing
basket number of $L$ by using a canonical Seifert surface $S(L)$.
Let us recall some notations and definitions before we state theorem.
For a Seifert graph $\Gamma$, the induced graph $G(\Gamma)$ of $\Gamma$ is a simple graph
obtained from $\Gamma$ by identifying edges in the same parallel class.
Since the link $L$ is not splittable, its Seifert graph $\Gamma$ is connected.
Since Seifert surfaces are orientable, the induced
graph $G(\Gamma)$ is bipartite. By applying Theorem~\ref{maintheorem1}, there exists a
co-tree  edge alternating spanning tree of $\Gamma$ and alternating labeling $\mu$ on $T$.
For an edge $e$ or a subgraph $H$ of $G$, we denoted
$\Gamma(e)= \{ \overline{e} \in \Gamma | \overline{e}$ reduced to $e \in G\}$
and  $\Gamma(H)= \{ \overline{f} \in \Gamma | \overline{f}$ reduced to $f$ where $f$ is en edge in $H\}$.

\begin{theorem}
Let $\Gamma$ be an Seifert graph of canonical Seifert surface $S$ of a link $L$ with
$|V(\Gamma)|=n$, $|E(\Gamma)|=m$ and the sign labeling $\phi$.
Let $G(\Gamma)$ be the induced graph of $\Gamma$.
Let $T$ be a co-tree  edge alternating spanning tree of $\Gamma$ and $\mu$ a
labeling on $T$ chosen by Theorem~\ref{maintheorem1}.
Let $\delta(T)$ be the cardinality of the set
$$\Psi(T) =\{e\in E(T)|~\mu(e)\neq\phi(\overline{e}) ~{\rm{for~all}}~ \overline{e} \in \Gamma(e) \},$$
and let $\zeta(T)$ be the cardinality of the set
$$\Upsilon(T) =\{\overline{e}\in E(\Gamma(T))~|~\mu(e)=\phi(\overline{e}),~\overline{e}\in \Gamma(e),~e \in E(T)-\Psi(T) \}.$$
and let $\eta(T)$ be the cardinality of the set
$$\Phi(T) = \{ \overline{e}\in E(\Gamma)-E(\Gamma(T)) ~|~\mu(\overline{e}) = \nu(e) \}$$
where $\nu(e) =+$($-$, resp.) if there is one extra positive(negative, respectively) sign
in the path $P_e$ joining end vertices of the edge $e$ in $T$.
Then the
flat plumbing basket number of $L$ is bounded by $m - 3(n-1) + 2 ( 2\delta(T)+\zeta(T) + \eta(T))$, $i. e.$,
$$ fpbk(L)\le m - 3(n-1) + 2 ( 2\delta(T)+\zeta(T) + \eta(T)).$$ \label{flatbktheorem3}
\begin{proof}
The proof consists two parts. First we examine the existence
of $3$-ball $B_{\alpha}$ for the validity of flat plumbings.
Second, we enumerate the total number of
flat plumbing used for $\overline{S}$.

We claim that there exists a co-tree  edge $e \in E(\Gamma)\setminus E(T)$
such that the half twisted band presented by the edge $e$ can be deplumbed.
We induct on the number of co-tree  edges in $E(\Gamma) \setminus E(T)$ and
the number of Seifert circles in lexicographic order.
If there is no co-tree  edge in $E(\Gamma)\setminus E(T)$ or there is only
one Seifert circle, then the Seifert surface $S$ is a disc $D_2$
so it is a flat plumbing basket surface.

Suppose there exist at least two Seifert circles and at least
one co-tree  edge in $E(\Gamma)\setminus E(T)$. We divide cases depending on
the existence of two concentric Seifert circles.
If there exist at least two Seifert circles which are concentric,
there exists at least one edge in $T$ and one edge $e$ in
$E(\Gamma)\setminus E(T)$ between two vertices which
correspond to two adjacent Seifert circles. Then the half twisted
band presented by the edge $e$ can be deplumbed
as we have proven for Theorem~\ref{fpbktheoremext}.

Suppose that there exist no Seifert circles which are concentric,
then the Seifert graph $\Gamma$ is planar.
Since $\Gamma$ is finite, there exists a co-tree  edge $e \in E(\Gamma)\setminus E(T)$
such that the disc $D$ bounded by $P_e$, the path joining
the both end vertices of the edge $e$ in $T$ and the edge $e$ does
not contain any other co-tree edges of $E(\Gamma)\setminus E(T)$.
Our next claim is that the interior of the disc $D$ bounded by $P_e$
and the edge $e$ does not contain any other edges of $\Gamma$.
Since the interior of the disc bounded by $P_e$ and the edge $e$
does not contain any co-tree edges of $E(\Gamma)\setminus E(T)$, then there might exist
some edges of the spanning tree $T$ in the interior of the disc bounded 
by $P_e$ and the edge $e$. Since there does not exist any co-tree edge
in the disc $D$, then there exists a vertex of valency $1$ in the disc $D$ but
it can be removed by Reidemeister move I, which reduces the number of Seifert circles.
By the induction hypothesis, we can see that the second claim is true.
Finally, we can see that the half twisted band presented by the edge $e$ can be deplumbed when
the $3$-ball $D_{\alpha}$ can be chosen along $P_e$ as depicted in Figure~\ref{fpbkprooffig}.

Next, we enumerate the total number of
flat plumbing used for $\overline{S}$. We divide $E(\Gamma)$ into
five subsets, $\Gamma(\Psi(T))$, $\Gamma(E(T)\setminus\Psi(T)) \cap \Upsilon(T)$,
$\Gamma(E(T)\setminus\Psi(T)) \setminus \Upsilon(T)$, $\Phi(T)$ and $E(\Gamma)- (E(\Gamma(T))\cup\Phi(T))$.

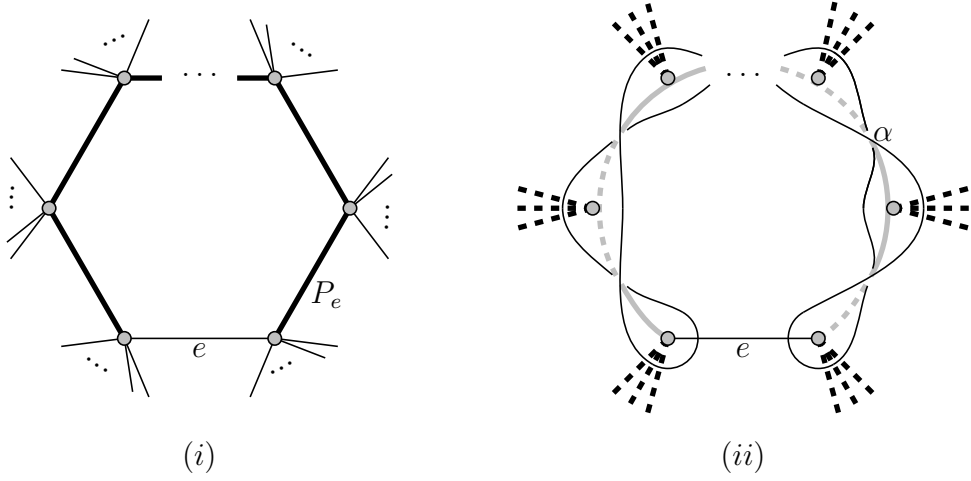
\begin{figure}
$$
\begin{pspicture}[shift=-3](-3.2,-3.5)(3.2,3)
\psline[linewidth=2pt](2;-60)(2;0)(2;60)(.5,1.732)
\psline[linewidth=2pt](-.5,1.732)(2;120)(2;180)(2;240)
\psline(2;-60)(2;-120) \psline(2;-60)(2.6;-45)
\psline(2;-60)(2.6;-50) \psline(2;-60)(2.6;-75)
\rput(2.4;-60){$\cdot$} \rput(2.4;-63){$\cdot$} \rput(2.4;-66){$\cdot$}
\psline(2;0)(2.6;15) \psline(2;0)(2.6;10) \psline(2;0)(2.6;-15)
\rput(2.5;0){$\cdot$} \rput(2.5;-3){$\cdot$} \rput(2.5;-6){$\cdot$}
\psline(2;60)(2.6;45) \psline(2;60)(2.6;70) \psline(2;60)(2.6;75)
\rput(2.5;60){$\cdot$} \rput(2.5;57){$\cdot$} \rput(2.5;54){$\cdot$}
\psline(2;120)(2.6;105) \psline(2;120)(2.6;130) \psline(2;120)(2.6;135)
\rput(2.5;120){$\cdot$} \rput(2.5;117){$\cdot$} \rput(2.5;114){$\cdot$}
\psline(2;180)(2.6;165) \psline(2;180)(2.6;190) \psline(2;180)(2.6;195)
\rput(2.5;180){$\cdot$}  \rput(2.5;177){$\cdot$} \rput(2.5;174){$\cdot$}
\psline(2;240)(2.6;225) \psline(2;240)(2.6;250) \psline(2;240)(2.6;255)
\rput(2.5;240){$\cdot$} \rput(2.5;237){$\cdot$} \rput(2.5;234){$\cdot$}
\rput[t](1.95;-30){$P_e$} \rput[t](1.8;270){$e$} \
\rput[t](1.8;90){$\ldots$} \rput[t](3.1;270){$(i)$}
\pscircle[fillstyle=solid,fillcolor=darkgray](2;-60){.1}
\pscircle[fillstyle=solid,fillcolor=darkgray](2;-120){.1}
\pscircle[fillstyle=solid,fillcolor=darkgray](2;60){.1}
\pscircle[fillstyle=solid,fillcolor=darkgray](2;120){.1}
\pscircle[fillstyle=solid,fillcolor=darkgray](2;0){.1}
\pscircle[fillstyle=solid,fillcolor=darkgray](2;-180){.1}
\end{pspicture}
\quad\quad
\begin{pspicture}[shift=-3](-3.2,-3.5)(3.2,3)
\pccurve[angleA=30,angleB=-120,linecolor=darkgray,linestyle=dashed,linewidth=2pt](2;-60)(1.925;-32)
\pccurve[angleA=60,angleB=-60,linecolor=darkgray,linewidth=2pt](1.925;-28)(1.925;28)
\pccurve[angleA=120,angleB=-15,linecolor=darkgray,linestyle=dashed,linewidth=2pt](1.925;32)(1.925;75)
\pccurve[angleA=150,angleB=-60,linecolor=darkgray,linewidth=2pt](2;-120)(1.91;-148)
\pccurve[angleA=120,angleB=-120,linecolor=darkgray,linestyle=dashed,linewidth=2pt](1.91;-152)(1.91;152)
\pccurve[angleA=60,angleB=-165,linecolor=darkgray,linewidth=2pt](1.91;148)(1.925;105)
\psline(2;-60)(2;-120)
\pccurve[angleA=30,angleB=-90](1.6;-60)(2.4;0)
\pccurve[angleA=30,angleB=-110](2.4;-60)(1.95;-32)
\pccurve[angleA=70,angleB=-90](1.9;-25)(1.6;0)
\pccurve[angleA=90,angleB=-30](2.4;0)(1.6;60)
\pccurve[angleA=-30,angleB=110](2.4;60)(1.95;32)
\pccurve[angleA=-70,angleB=90](1.9;25)(1.6;0)
\pccurve[angleA=150,angleB=20](2.4;60)(2.1;75)
\pccurve[angleA=150,angleB=-40](1.6;60)(1.7;75)
\pccurve[angleA=-30,angleB=110](2.4;60)(1.95;32)
\pccurve[angleA=-70,angleB=90](1.9;25)(1.6;0)
\pccurve[angleA=150,angleB=-90](2.4;-120)(1.6;180)
\pccurve[angleA=150,angleB=-40](1.6;-120)(1.86;-146)
\pccurve[angleA=140,angleB=-90](1.95;-153)(2.4;180)
\pccurve[angleA=90,angleB=-150](1.6;180)(2.4;120)
\pccurve[angleA=-150,angleB=40](1.6;120)(1.86;146)
\pccurve[angleA=-140,angleB=90](1.95;153)(2.4;180)
\pccurve[angleA=30,angleB=160](2.4;120)(2.1;105)
\pccurve[angleA=30,angleB=-140](1.6;120)(1.7;105)
\psarc(2;-60){.4}{120}{300}
\psarc(2;-120){.4}{240}{60}
\psline[linestyle=dashed,linewidth=2pt](2;-60)(3;-60)
\psline[linestyle=dashed,linewidth=2pt](2;-60)(3;-55)
\psline[linestyle=dashed,linewidth=2pt](2;-60)(3;-65)
\psline[linestyle=dashed,linewidth=2pt](2;0)(3;0)
\psline[linestyle=dashed,linewidth=2pt](2;0)(3;5)
\psline[linestyle=dashed,linewidth=2pt](2;0)(3;-5)
\psline[linestyle=dashed,linewidth=2pt](2;60)(3;60)
\psline[linestyle=dashed,linewidth=2pt](2;60)(3;55)
\psline[linestyle=dashed,linewidth=2pt](2;60)(3;66)
\psline[linestyle=dashed,linewidth=2pt](2;-120)(3;-120)
\psline[linestyle=dashed,linewidth=2pt](2;-120)(3;-115)
\psline[linestyle=dashed,linewidth=2pt](2;-120)(3;-125)
\psline[linestyle=dashed,linewidth=2pt](2;120)(3;120)
\psline[linestyle=dashed,linewidth=2pt](2;120)(3;115)
\psline[linestyle=dashed,linewidth=2pt](2;120)(3;125)
\psline[linestyle=dashed,linewidth=2pt](2;180)(3;180)
\psline[linestyle=dashed,linewidth=2pt](2;180)(3;175)
\psline[linestyle=dashed,linewidth=2pt](2;180)(3;185)
\pscircle[fillstyle=solid,fillcolor=darkgray](2;-60){.1}
\pscircle[fillstyle=solid,fillcolor=darkgray](2;-120){.1}
\pscircle[fillstyle=solid,fillcolor=darkgray](2;0){.1}
\pscircle[fillstyle=solid,fillcolor=darkgray](2;60){.1}
\pscircle[fillstyle=solid,fillcolor=darkgray](2;120){.1}
\pscircle[fillstyle=solid,fillcolor=darkgray](2;180){.1}
\rput[t](1.8;270){$e$}\rput[t](1.8;90){$\ldots$}
\rput[t](3.1;270){$(ii)$} \rput[t](2.15;30){$\alpha$}
\end{pspicture}
$$
\caption{ $(i)$ A part of the graph containing the path $P_e$ and the edge $e$ in $\Gamma$ where the interior of the disc bounded by $P_e$ and the edge $e$
does not contain any other edges of $\Gamma$ $(ii)$ A 3-ball $D_{\alpha}$ which is a 3-dimensional tube
 along $\alpha$ attached from the side of $S$ which has the solid gray line where the
thick black dashed lines present some bands towards the outside and
thick gray line correspond to the arc $\alpha$ of deplumbing.} \label{fpbkprooffig}
\end{figure}

As we have seen in Example~\ref{52flat}, for each edge $e \in \Psi(T)$,
we have to use a Reidemeister move II for flat plumbings.
Thus, for an edge $e\in \Psi(T)$, we need $|\Gamma(e)|+1$ flat plumbings.
For an edge $e\in E(T)- \Psi(T)$, one edge $\overline{e}\in \Gamma(e)\cap \Upsilon(T)$
will be used in for the co-tree edge spanning tree $T$. For rest of edges
in $\Gamma(e)\cap \Upsilon(T)$ contributes three flat plumbings while
one flat plumbing for all edges in $\Gamma(e)\cap (E(\Psi(T)) \cup\Upsilon(T))^c$.

As we have seen in Example~\ref{52flat}, for each co-tree edge $f$ in $E(\Gamma)-\Gamma(T)$,
if $f\in E(\Gamma)- (E(\Gamma(T))\cup\Phi(T))$, $i.e.$, $\nu(f)$ is different from $\phi(f)$,
then the half twisted band presented by $f$ can be obtained by one flat plumbing.
If $f\in \Phi(T)$, then it requires three flat plumbings where two of three flat plumbings
are required to change the sign of $\phi(f)$ to the opposite sign
as illustrated in Figure~\ref{3annuli}.
Thus, for an edge $e\in \Phi(T)$, we need three flat plumbings while
we need one flat plumbing for an edge $\overline{e}\in E(\Gamma)- (E(\Gamma(T))\cup\Phi(T))$.

By summing all numbers of flat plumbings with respect to
a partition of the edge set $E(\Gamma)$ in five subsets,
we have

\begin{align*}
fpbk(L)&\le\sum_{e \in \Psi(T)} (|\Gamma(e)|+1) + \sum_{e \in E(T)-\Psi(T)}
\left( \sum_{\overline{e} \in \Gamma(e)\cap\Upsilon(T)} 3 + \sum_{\overline{e} \in \Gamma(e) \cap (\Gamma(T)\setminus \Upsilon(T))} 1  \right) \\
&-3(n-1-\delta(T)) + \sum_{\overline{e} \in \Psi(T)} 3 + \sum_{e \in (E(\Gamma)- (E(\Gamma(T))\cup\Phi(T))} 1\\
&= \left(\sum_{e \in \Psi(T)} |\Gamma(e)|\right) +\delta(T) + \left(\sum_{e \in E(T)\setminus \Psi(T)} |\Gamma(e)|\right) + 2 \zeta(T) \\
& -3(n-1-\delta(T))+ \left(\sum_{e\in E(G)\setminus E(T)} |\Gamma(e)| \right) + 2\eta(T) \\
&= m - 3(n-1) + 2 ( 2\delta(T)+\zeta(T) +\eta(T)).
\end{align*}

This completes the proof of the theorem.
\end{proof}
\end{theorem}

For the Seifert graph $\Gamma(S(7_5))$ of $7_5$ in Figure~\ref{graph} $(ii)$ and the
co-tree alternating spanning tree $\overline{T}$ in Figure~\ref{75open} $(ii)$,
let us check how Theorem~\ref{flatbktheorem3} can be applied.
First of all, the number of vertices and the number of edges in $\Gamma(S(7_5))$
are $4$ and $7$ respectively.
Two edges $\{c,d\}$ and $\{a,b\}$ belong to $\Psi(T)$ so $\delta(T)=2$.
The edge $\{a,d\}$ belongs to $\Upsilon(T)$ and so $\zeta(T)=1$.
None of three edges between vertex $b$ and $c$ belong to $\Phi(T)$ thus
$\eta(T)=0$. The upper bound obtained in Theorem~\ref{flatbktheorem3}
is $7-3(4-1)+2(2\cdot 2+1+0)=8$ as we have seen in Example~\ref{52flat}.

Using Theorem~\ref{flatbktheorem3}, we obtain an upper bound
for the basket number of a link using the canonical Seifert surface in the
following corollary.

\begin{cor}
Let $\Gamma$ be an Seifert graph of canonical Seifert surface $S$ of a link $L$ with
$|V(\Gamma)|=n$, $|E(\Gamma)|=m$ and the sign labeling $\phi$.
Let $G(\Gamma)$ be the induced graph of $\Gamma$.
Let $T$ be a co-tree  edge alternating spanning tree of $\Gamma$ and $\mu$ a
labeling on $T$ chosen by Theorem~\ref{maintheorem1}.
Let $\delta(T)$ be the cardinality of the set
$$\Psi(T) =\{e\in E(T)|~\mu(e)\neq\phi(\overline{e}) ~{\rm{for~all}}~ \overline{e} \in \Gamma(e) \}.$$
Then the basket number of $L$ is bounded by $m -n + 1 + 2\delta(T)$, $i. e.$,
$$ bk(L)\le m -n + 1 + 2\delta(T).$$ \label{bktheorem3}
\begin{proof}
In Theorem~\ref{flatbktheorem3}, whenever we perform three flat plumbings
for the edges in $\Upsilon(T)$ or $\Phi(T)$, it can be done by a single $A_{\pm 1}$ plumbing
for a basket surface. Thus, we have

\begin{align*}
bk(L)&\le\sum_{e \in \Psi(T)} (|\Gamma(e)|+1) + \sum_{e \in E(T)-\Psi(T)} (|\Gamma(e)|-1)
+ \sum_{e \in E(G(\Gamma))-E(T)} |\Gamma (e)|  \\
&=\left(\sum_{e \in \Psi(T)} |\Gamma(e)|\right) +\delta(T) + \left(\sum_{e \in E(T)-\Psi(T)} |\Gamma(e)|\right) -(n-1-\delta(T))\\
&+ \left(\sum_{e \in E(G(\Gamma))-E(T)} |\Gamma (e)|\right) \\
&= m - n +1 + 2\delta(T).
\end{align*}
\end{proof}
\end{cor}

\begin{figure}
$$
\begin{pspicture}[shift=-2](-.2,-2.8)(3.7,2.2)
\psline(.75,2)(3.25,2)  \psline(.75,1.2)(1,1.2)
\psline(1.5,1.2)(2,1.2) \psline(2.5,1.2)(3,1.2)
\psline(.25,.4)(.5,.4) \psline(0,0)(0,-.4)
\psline(1.5,.4)(1.5,-.4) \psline(2.5,-.4)(3,-.4)
\psline(1.5,-1.2)(1.75,-1.2) \psline(2,-.8)(2,-.4)
\psline(3.5,-.4)(3.5,-1.2) \psline(1,-2)(3,-2)
\psline(3.5,1.2)(3.5,1.6)
\pccurve[angleA=180,angleB=90](.25,2)(0,1.6)
\pccurve[angleA=180,angleB=-90](.25,1.2)(0,1.6)
\pccurve[angleA=0,angleB=180](.25,2)(.75,1.2)
\pccurve[angleA=0,angleB=-120](.25,1.2)(.45,1.5)
\pccurve[angleA=60,angleB=180](.55,1.7)(.75,2)
\pccurve[angleA=0,angleB=120](1,1.2)(1.25,.8)
\pccurve[angleA=-60,angleB=90](1.25,.8)(1.5,.4)
\pccurve[angleA=180,angleB=60](1.5,1.2)(1.3,.85)
\pccurve[angleA=-120,angleB=60](1.2,.75)(.3,-.75)
\pccurve[angleA=-120,angleB=90](.2,-.85)(0,-1.2)
\pccurve[angleA=-90,angleB=180](0,-1.2)(.5,-2)
\pccurve[angleA=0,angleB=-120](.5,-2)(1,-1.2)
\pccurve[angleA=60,angleB=-120](1,-1.2)(1.2,-.85)
\pccurve[angleA=60,angleB=-90](1.3,-.75)(1.5,-.4)
\pccurve[angleA=0,angleB=120](.5,.4)(.7,.1)
\pccurve[angleA=-60,angleB=120](.8,-.1)(1.25,-.8)
\pccurve[angleA=-60,angleB=180](1.25,-.8)(1.5,-1.2)
\pccurve[angleA=0,angleB=-90](1.75,-1.2)(2,-.8)
\pccurve[angleA=180,angleB=90](.25,.4)(0,0)
\pccurve[angleA=-90,angleB=120](0,-.4)(.25,-.8)
\pccurve[angleA=-60,angleB=120](.25,-.8)(.7,-1.7)
\pccurve[angleA=-60,angleB=180](.8,-1.8)(1,-2)
\pccurve[angleA=0,angleB=120](2,1.2)(2.25,.8)
\pccurve[angleA=-60,angleB=90](2.25,.8)(2.5,.4)
\pccurve[angleA=-90,angleB=60](2.5,.4)(2.3,.05)
\pccurve[angleA=-120,angleB=90](2.2,-.05)(2,-.4)
\pccurve[angleA=180,angleB=60](2.5,1.2)(2.3,.85)
\pccurve[angleA=-120,angleB=90](2.2,.75)(2,.4)
\pccurve[angleA=-90,angleB=120](2,.4)(2.25,0)
\pccurve[angleA=-60,angleB=180](2.25,0)(2.5,-.4)
\pccurve[angleA=0,angleB=90](3.25,2)(3.5,1.6)
\pccurve[angleA=0,angleB=-90](3,-2)(3.5,-1.2)
\pccurve[angleA=0,angleB=120](3,1.2)(3.25,.8)
\pccurve[angleA=-60,angleB=90](3.25,.8)(3.5,.4)
\pccurve[angleA=-90,angleB=60](3.5,.4)(3.3,.05)
\pccurve[angleA=-120,angleB=0](3.2,-.05)(3,-.4)
\pccurve[angleA=-90,angleB=60](3.5,1.2)(3.3,.85)
\pccurve[angleA=-120,angleB=90](3.2,.75)(3,.4)
\pccurve[angleA=-90,angleB=120](3,.4)(3.25,0)
\pccurve[angleA=-60,angleB=90](3.25,0)(3.5,-.4)
\rput[t](1.75,-2.3){$(i)$}
\end{pspicture} \quad
\begin{pspicture}[shift=-2](-.8,-2.8)(1.7,2.2)
\psline(0,1)(-.5,0)(-.5,-1)(.5,-2)
\psline(.5,0)(0,1)(1.5,0)(.5,-1)(.5,0)(0,1)(-.5,0)(.5,-1)(.5,-2)
\rput(0,1){$\bullet$}
\rput(-.5,0){$\bullet$} \rput(.5,0){$\bullet$}
\rput(1.5,0){$\bullet$} \rput(-.5,-1){$\bullet$}
\rput(.5,-1){$\bullet$} \rput(.5,-2){$\bullet$}
\rput(-.5,.4){$+$}
\rput(1,.6){$+$} \rput(.1,.4){$+$}
\rput(1.3,-.5){$+$} \rput(.7,-.5){$+$}
\rput(-.7,-.5){$+$} \rput(-.25,-.5){$+$}
\rput(.8,-1.5){$+$} \rput(-.3,-1.5){$+$}
\rput(-.8,0){$v_1$}
\rput(.8,0){$v_3$} \rput(1.7,-.3){$v_4$}
\rput(-.8,-1.1){$v_5$} \rput(.8,-1.1){$v_6$}
\rput(.8,-2.1){$v_7$} \rput(0,1.3){$v_1$}
\rput[t](1.5,1){$\Gamma$} \rput[t](.5,-2.3){$(ii)$}
\end{pspicture}
\quad
\begin{pspicture}[shift=-2](-.8,-2.8)(1.7,2.2)
 \psline(0,1)(-.5,0)
\psline(-.5,0)(-.5,-1) \psline(-.5,-1)(.5,-2)
\psline[linestyle=dashed](.5,0)(0,1)
\psline[linestyle=dashed](0,1)(1.5,0)
\psline(1.5,0)(.5,-1) \psline(.5,-1)(.5,0)
\psline(0,1)(-.5,0) \psline[linestyle=dashed](-.5,0)(.5,-1)
\psline(.5,-1)(.5,-2) \rput(0,1){$\bullet$}
 \rput(-.5,0){$\bullet$}
\rput(.5,0){$\bullet$} \rput(1.5,0){$\bullet$}
\rput(-.5,-1){$\bullet$} \rput(.5,-1){$\bullet$}
\rput(.5,-2){$\bullet$}
\rput(-.5,.4){$+$} \rput(1.3,-.5){$+$}
\rput(-.3,-1.5){$+$} \rput(-.7,-.5){$-$}
\rput(.3,-.5){$+$} \rput(.3,-1.5){$-$}
\rput(0,1.3){$v_1$} \rput[t](1.5,1){$\overline T_1$}
\rput[t](.5,-2.3){$(iii)$}
\end{pspicture}
\quad
\begin{pspicture}[shift=-2](-.8,-1.8)(1.7,2.5)
\psline(0,2)(-.5,1)
\psline[linestyle=dashed](-.5,1)(-.5,0) \psline(-.5,0)(.5,-1)
\psline[linestyle=dashed](.5,1)(0,2)
\psline[linestyle=dashed](0,2)(1.5,1)
\psline(1.5,1)(.5,0) \psline(.5,0)(.5,1)
\psline(0,2)(-.5,1) \psline(-.5,1)(.5,0)
\psline(.5,0)(.5,-1) \rput(0,2){$\bullet$}
 \rput(-.5,1){$\bullet$}
\rput(.5,1){$\bullet$} \rput(1.5,1){$\bullet$}
\rput(-.5,0){$\bullet$} \rput(.5,0){$\bullet$}
\rput(.5,-1){$\bullet$}
\rput(-.5,1.4){$+$} \rput(1.3,.5){$+$}
\rput(-.2,-.6){$-$} \rput(-.2,.4){$-$}
\rput(.7,.5){$+$} \rput(.3,-.5){$+$}
\rput(0,2.3){$v_1$} \rput[t](1.5,2){$\overline T_2$}
\rput[t](.5,-1.3){$(iv)$}
\end{pspicture}
$$
$$
\begin{pspicture}[shift=-2](-.8,-1.8)(1.8,2.5)
 \psline[linestyle=dashed](-.5,1)(0,2)
\psline(.5,1)(0,2) \psline[linestyle=dashed](1.5,1)(0,2)
\psline(.5,0)(.5,1) \psline(-.5,0)(-.5,1)
\psline(.5,0)(1.5,1) \psline(.5,0)(-.5,1)
\psline[linestyle=dashed](.5,-1)(.5,0)
\psline(.5,-1)(-.5,0) \rput(0,2){$\bullet$}
 \rput(-.5,1){$\bullet$}
\rput(.5,1){$\bullet$} \rput(1.5,1){$\bullet$}
\rput(-.5,0){$\bullet$} \rput(.5,0){$\bullet$}
\rput(.5,-1){$\bullet$}
\rput(-.2,.4){$+$} \rput(1.2,.4){$+$}
\rput(.1,1.4){$+$} \rput(-.7,.5){$-$}
\rput(.3,.5){$-$} \rput(-.3,-.6){$+$}
\rput(0,2.3){$v_1$} \rput[t](1.5,2){$\overline T_3$}
\rput[t](.5,-1.5){$(v)$}
\end{pspicture}
\quad
\begin{pspicture}[shift=-2](-.8,-1.8)(1.8,2.5)
 \psline[linestyle=dashed](-.5,1)(0,2)
\psline(.5,1)(0,2) \psline[linestyle=dashed](1.5,1)(0,2)
\psline(.5,0)(.5,1) \psline(-.5,0)(-.5,1)
\psline(.5,0)(1.5,1) \psline[linestyle=dashed](.5,0)(-.5,1)
\psline(.5,-1)(.5,0) \psline(.5,-1)(-.5,0)
\rput(0,2){$\bullet$}
\rput(-.5,1){$\bullet$} \rput(.5,1){$\bullet$}
\rput(1.5,1){$\bullet$} \rput(-.5,0){$\bullet$}
\rput(.5,0){$\bullet$} \rput(.5,-1){$\bullet$}
 \rput(.7,-.4){$+$}
\rput(1.2,.4){$+$} \rput(.1,1.4){$+$}
\rput(-.7,.5){$+$} \rput(.3,.5){$-$}
\rput(-.3,-.6){$-$} \rput(0,2.3){$v_1$}
\rput[t](1.5,2){$\overline T_4$}
\rput[t](.5,-1.5){$(vi)$}
\end{pspicture}
\quad
\begin{pspicture}[shift=-2](-.8,-1.8)(1.8,2.5)
\psline[linestyle=dashed](-.5,1)(0,2)
\psline[linestyle=dashed](.5,1)(0,2)
\psline(1.5,1)(0,2) \psline(.5,0)(.5,1)
\psline(-.5,0)(-.5,1) \psline(.5,0)(1.5,1)
\psline(.5,0)(-.5,1) \psline[linestyle=dashed](.5,-1)(.5,0)
\psline(.5,-1)(-.5,0) \rput(0,2){$\bullet$}
 \rput(-.5,1){$\bullet$}
\rput(.5,1){$\bullet$} \rput(1.5,1){$\bullet$}
\rput(-.5,0){$\bullet$} \rput(.5,0){$\bullet$}
\rput(.5,-1){$\bullet$}
\rput(-.2,.4){$+$} \rput(1.2,.4){$-$}
\rput(1,1.6){$+$} \rput(-.7,.5){$-$}
\rput(.3,.5){$+$} \rput(-.3,-.6){$+$}
\rput(0,2.3){$v_1$} \rput[t](1.5,2){$\overline T_5$}
\rput[t](.5,-1.5){$(vii)$}
\end{pspicture}
\quad
\begin{pspicture}[shift=-2](-.8,-1.8)(1.8,2.5)
 \psline[linestyle=dashed](-.5,1)(0,2)
\psline[linestyle=dashed](.5,1)(0,2)
\psline(1.5,1)(0,2) \psline(.5,0)(.5,1)
\psline(-.5,0)(-.5,1) \psline(.5,0)(1.5,1)
\psline[linestyle=dashed](.5,0)(-.5,1)
\psline(.5,-1)(.5,0) \psline(.5,-1)(-.5,0)
\rput(0,2){$\bullet$} \rput(-1.5,1){$\bullet$}
\rput(-.5,1){$\bullet$} \rput(.5,1){$\bullet$}
\rput(1.5,1){$\bullet$} \rput(-.5,0){$\bullet$}
\rput(.5,0){$\bullet$} \rput(.5,-1){$\bullet$}
 \rput(.7,-.4){$+$}
\rput(1.2,.4){$-$} \rput(1,1.6){$+$}
\rput(-.7,.5){$+$} \rput(.3,.5){$+$}
\rput(-.3,-.6){$-$} \rput(0,2.3){$v_1$}
\rput[t](1.5,2){$\overline T_6$}
\rput[t](.5,-1.5){$(viii)$}
\end{pspicture}
$$
\caption{Figures described in Example~\ref{fpbsexa2}.}  \label{graphtree1}
\end{figure}
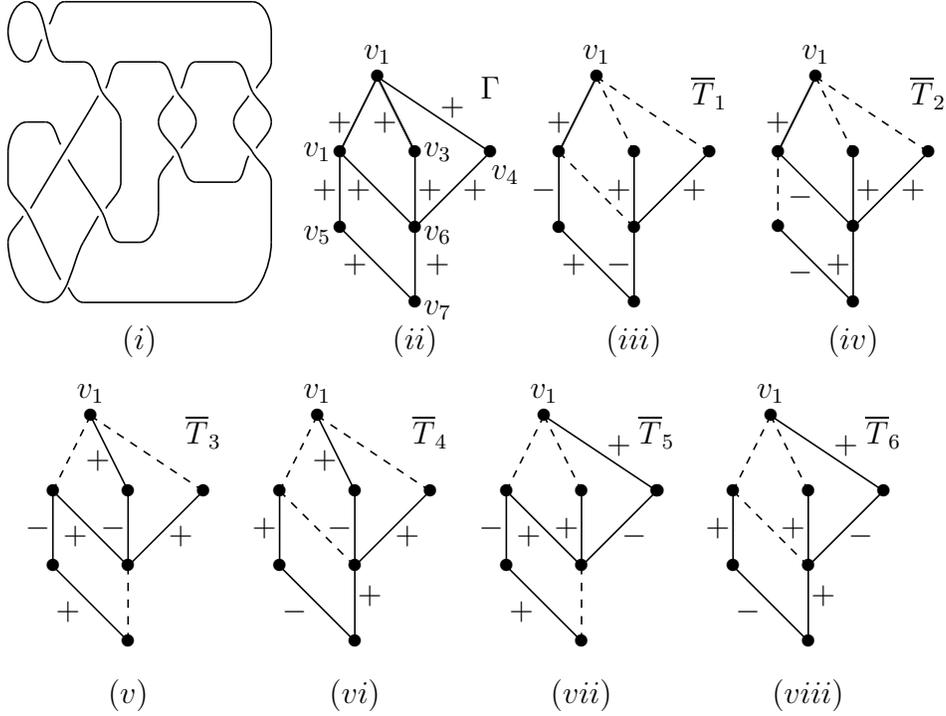

Since $m, n$ are given, the minimum can be attained 
by the minimum of $2\delta(T)+\zeta(T) +\eta(T)$
among all co-tree  edge alternating spanning trees of $G$. Thus, we may obtain a
better upper bound for Theorem~\ref{flatbktheorem3} if we
use the method in the proof of Theorem~\ref{fpbktheoremext} by
considering all possible pairs of (co-tree edge alternating 
spanning tree $T$, an alternating labeling $\mu$ on $T$)
obtained from Theorem~\ref{maintheorem1} as demonstrated in the following example.
Let us remark that we ignore the vertex of valency $1$ 
since it can simply be removed by Reidemeister move I.

\begin{exa} \label{fpbsexa2}
For the link depicted in Figure~\ref{graphtree1} $(i)$ for which the
Seifert graph is given in Figure~\ref{graph1} $(i)$, if we fix the
vertex $v_1$ as a root and the depth labeling with the root $v_1$,
the upper bound attained in Theorem~\ref{flatbktheorem3} is $10$
which can be obtained from co-tree edge alternating spanning tree $\overline{T}_1$
and $\overline{T}_2$ as depicted in Figure~\ref{graphtree1}.
\begin{proof}
If we fix the vertex $v$ and the depth labeling with a root $v$,
there are six possible spanning trees for which Theorem~\ref{maintheorem1}
holds as illustrated in Figure~\ref{graphtree1} $(iii-viii)$.
From spanning tree $\overline{T}_i$, for $i=1, 2$, we have
$2(2\delta(\overline{T}_i)+\zeta(\overline{T}_i)+\eta(T_i))=20$.
For the others $3\le j\le 6$, we have $2(2\delta(\overline{T}_j)+\zeta(\overline{T}_j)+\eta(T_j))=22$.
Therefore, the upper bound for the flat plumbing basket number attained in
Theorem~\ref{flatbktheorem3} for the fixed vertex $v_1$ is $10$.
By considering all other vertex in $\Gamma$ as a root,
we find that the smallest upper bound for the flat
plumbing basket number attained in Theorem~\ref{flatbktheorem3} is $10$.
\end{proof}
\end{exa}

\subsection{Flat plumbing number} \label{flats}

For the flat plumbing number,
we have much more flexibility that
the arc $\alpha$ we are plumbing along
can pass an annulus we have made by previous plumbings.

Hayashi and Wada~\cite{HW:plumbing} first showed the existence of
a flat plumbing surface of a given link $L$ as follow. First they showed that
the link $L$ can be presented by closed braid $\overline{\beta}$, and we can
assume $\beta \in B_n$ contains all generators and their inverse at least once by
applying the Reidemeister move $II$.
Then they choose $D^2$ to be the disc bounded by Seifert circles
and connected by half twisted bands presented by all generators exactly once.
Then, all inverse of generators can be obtained by $A_0$ annulus plumbing
as described in Figure~\ref{52complete} where all paths are contained in $D^2$. Then, for
the other words in $\beta$ can be constructed by $A_0$ annulus plumbing along a path
in the flat plumbing surface of the previous step.
By careful managing unnecessary $A_0$ annulus plumbings,
we obtain a new upper bound for the flat plumbing number.

\begin{theorem}
Let $L$ be an oriented link which is a closed $n$-braid $\overline{\beta}$ where
$\beta\in B_n$ is reduced and the length of $\beta$ is $m$.
Let $\gamma$ be the cardinality of the set
$$\Omega= \{ i | 1
\le i \le n-1, \sigma_i ~{\rm{and}}~ \sigma_i^{-1}~ {\rm{both}~\rm{appear}~ \rm{in}}~ \beta\}.$$
Then there exists a flat
plumbing surface $S$ obtained by at most $m + n -1-2\gamma$
flat plumbings such that $\partial S$ is
isotopic to $L$, $i.e.,$ $fp(L)\le m+n-1-2\gamma$. \label{fptheorem1}
\end{theorem}
\begin{proof}
Since $L$ is nor splittable, for $i=1, 2, \ldots, n-1$,
its braid presentative $\beta\in B_n$ contains
at least one of $\sigma_i$ or $\sigma_i^{-1}$ by Remark~\ref{remark0} $(1)$.
For $i \in \Omega$, we choose $\sigma_i^{\epsilon_i}$
if $\sigma_i^{\epsilon_i}$ appears no more than $\sigma_i^{{-}\epsilon_i}$ does.
For $i \not\in \Omega$, we choose $\sigma_i^{\epsilon_i}$
to be the one of $\sigma_i$ and $\sigma_i^{-1}$ which appears in $W$.
Since $L$ is prime and reduced, $\sigma_i^{\epsilon_i}$ has to
appear at least twice in $W$ by Remark~\ref{remark0} $(2)$.
Let $D^2$ be the $2$-disc bounded by Seifert circles
and connected by half twisted bands presented by the
above $\sigma_i^{\epsilon_i}$ for $i=1, 2, \ldots, n-1$.
To construct a flat plumbing surface $S$,
we divide cases depend on $i\in \Omega$ or not.

If $i\not\in \Omega$ and $\sigma_i^{\epsilon_i}$ appear more than once,
we perform a Reidemeister move II between $i$-th and
$i+1$-th string of the closed braids, and the effect of this move
inserts $\sigma_i\sigma_i^{-1}$ to the word $W$.
First, we perform $A_0$ plumbing along
an arc passing on the disc $D^2$ to generate the half twisted band
presented by $\sigma_i^{{-}\epsilon_i}$ which was produced by Reidemeister move II.
All the other half twisted bands can be obtained by $A_0$ plumbing along
an arc passing on the disc $D^2$ and the half twisted band presented by $\sigma_i^{-\epsilon_i}$.

If $i\in \Omega$, we first perform $A_0$ plumbing along
an arc passing on the disc $D^2$ to generate all half twisted bands
presented by $\sigma_i^{-\epsilon_i}$. For the other
half twisted bands can be obtained by $A_0$ plumbing along
an arc passing on the disc $D^2$ and the half twisted band presented by $\sigma_i^{-\epsilon_i}$.

The canonical Seifert surface of a closed braid
$\overline{W}$ which is obtained from $W$ by $n-\gamma-1$ times
insertions of $\sigma_i\sigma_i^{-1}$ for $i\not\in \Omega$,
is a flat plumbing surface whose boundary is
isotopic to $L$. The total number of $A_0$ plumbing
is $m+n-1-2\gamma$, we obtained the desired upper bound in the theorem.
\end{proof}

Let us remark that if $W\in B_n$ contains both of
$\sigma_i$, $\sigma_i^{-1}$ for all $i=1,2, \ldots, n-1$, then
the upper bound obtained by Theorem~\ref{fptheorem1}
is $m-n+1$ which is the best upper bound possible.
If $\beta$ contains only one of $\sigma_i$,
$\sigma_i^{-1}$ for all $i=1,2, \ldots, n-1$ which include a positive (all letters in the
braid word $\beta$ are $\sigma_i$) or negative (all letters in the
braid word $\beta$ are $\sigma_i^{-1}$), then the
upper bound obtained by Theorem~\ref{fptheorem1} is $m+n-1$.

Now we want to find an upper bound for the flat
plumbing number of $L$ by using a canonical Seifert surface $S(L)$.
Because the Seifert graph of the canonical Seifert
surface of a closed braid is a path with multi edges,
Theorem~\ref{maintheorem1} can be directly obtained
without changing anything.
For a canonical Seifert surface of an arbitrary
diagram of a link, we will use a co-tree  edge alternating spanning tree of
the Seifert graph of the given link diagram, but
the number of flat plumbings are reduced compare to
the case for the flat plumbing basket surfaces.
By using the similar idea of the proof of Theorem~\ref{fptheorem1},
we obtain the following theorem.

\begin{theorem}
Let $S$ be a canonical Seifert surface of a link $L$. Let $\Gamma$ be the
Seifert graph of the Seifert surface $S$ with
edge labeling $\phi : E(\Gamma) \rightarrow \{+,-\}$ and
$G(\Gamma)$ be the induced graph of $\Gamma$.
Let $|E(\Gamma)|=m$ and $|V(\Gamma)|=n$.
By Theorem~\ref{maintheorem1}, there exists a co-tree  edge alternating
spanning tree $T$ of $G(\Gamma)$ with respect to a labeling $\mu$.
Let $\delta$ to be the cardinality of the set
$$\Psi(T) =\{e\in E(T)|~\mu(e)\neq\phi(\overline{e}) ~{\rm{for~all}}~ \overline{e} \in \Gamma(e) \}.$$
Let $\xi$ to be the cardinality of the set
$$\overline{\Phi}(T) =\{e\in E(G)-E(T)| \phi(\overline{e}) =\phi(\overline{f}), ~
\nu(e) = \phi(\overline{e})
~{\rm{for~all}}~ \overline{e}, \overline{f} \in \Gamma(e) \}$$
where $\nu(e) =+$($-$, resp.) if there is one extra positive(negative, respectively) sign
in the path $P_e$ joining end vertices of the edge $e$ in $T$.
Then, the flat plumbing number of $L$ is bounded by
$m- n+1 + 2( \delta(T) + \xi(T))$, $i. e.$,
$$ fp(L)\le m- n+1 + 2( \delta(T) + \xi(T)). $$ \label{flattheorem2}
\begin{proof}
Since the link $L$ is not splittable, its Seifert graph $\Gamma$ is connected.
Since Seifert surfaces are orientable, the induced
graph $G(\Gamma)$ of $\Gamma$ is bipartite.
For each edge $e \in \Psi(T)$,
we have to use a Reidemeister type II move for flat plumbings.
Thus, for an edge $e\in \Psi(T)$, we need $|\Gamma(e)| +1$ many flat plumbings while
we need $|\Gamma(e)| -1$ many flat plumbings for an edge $e\in E(T)\setminus  \Psi(T)$.

In Example~\ref{52flat}, for each co-tree  edge $\overline{f}$ in $E(\Gamma) \setminus  \Gamma(E(T))$,
$\nu(f)$ is different from $\phi(f)$ where $f$ is the parallel class of $\overline{f}$ in $G(\Gamma)$.
Then, it can be obtained by a flat plumbing.
If $\nu(f)$ is same to $\phi(f)$, $i.e.$,
$f \in\overline{\Phi}(T)$,
then to obtain the half twisted band presented by $\overline{f}$,
we require three flat plumbings, two of three flat plumbings
are required to change the sign of $\overline{f}$ to the opposite sign
as illustrated in Figure~\ref{3annuli}. For flat plumbing basket
surfaces, this phenomena is inevitable. However, for
flat plumbing surfaces, the arc of which we are plumbing can pass
the previously added annuli. So, adding three
flat plumbings are necessary just once for $f\in \Phi(T)$.

For $e\in E(G(\Gamma))- (E(T)\cup\Phi(T))$, $\Gamma(e)$ contains
two edges of the opposite signs. Thus, one of these edges has different $\nu$ and $\phi$
which can be realized by a flat plumbing along an arc in $D^2$.
Then all edges in $\Gamma(e)$ of the opposite sign can be realized by a flat plumbing along an arc in
the annulus the previously plumbed. All remaining edges $\Gamma(e)$ can
be realized by a flat plumbing along an arc in
the annulus obtained for the edges in $\Gamma(e)$ of the opposite sign.
Therefore, for an edge $e\in \Phi(T)$, we need $|\Gamma(e)| +2$ many flat plumbings while
we need $|\Gamma(e)|$ many flat plumbings for an edge $e\in E(G(\Gamma))- (E(T)\cup\Phi(T))$.

By summing all numbers of flat plumbings with respect to
a partition of the edge set $E(\Gamma)$ in four subsets,
we have

\begin{align*}
fp(L)&\le\sum_{e \in \Psi(T)} (|\Gamma(e)| +1) + \sum_{e \in E(T)\setminus \Psi(T)} (|\Gamma(e)| -1)
+ \sum_{e \in \Phi(T)} (|\Gamma(e)| +2) + \sum_{e \in H} |\Gamma(e)|\\
&= \left(\sum_{e \in E(T)} (|\Gamma(e)|-1)\right) + 2|\Psi(T)|
+ \left(\sum_{e \in E(G)\setminus E(T)} |\Gamma(e)| \right) + 2|\Phi(T)|\\
&= \left(\sum_{e \in E(G)} |\Gamma(e)| \right) - |E(T)| + 2( \delta(T) + \xi(T)) \\
&= m-(n-1)+2( \delta(T) + \xi(T))
\end{align*}
where $H=E(G(\Gamma))- (E(T)\cup\Phi(T))$.
\end{proof}
\end{theorem}

By considering all possible co-tree  edge alternating spanning tree $T$s of $G(\Gamma)$
and since $\delta(T)+\xi(T) \le |E(T)|+|E(G)|-|E(T)|=|E(G)|$,
we obtain the following corollary.

\begin{cor}\label{fpcor1}
Let $S$ be a canonical Seifert surface of a link $L$ with $s(S)$
Seifert circles and $c(S)$ half twisted bands. Let $\Gamma$ be the
Seifert graph of the Seifert surface $S$ with
edge labeling $\phi : E(\Gamma) \rightarrow \{+,-\}$.
Let $\lambda(T)$ be the the minimum of $\delta(T) + \xi(T)$
over all co-tree  edge alternating spanning tree $T$s of  $G(\Gamma)$.
Then,
$$fp(L) \le c(S)-s(S)+1 + 2 \lambda(T) \le c(S)- s(S)+1 + 2|E(G)|.$$
\end{cor}

The upper bound in Corollary~\ref{fpcor1} can be written in a different form
if we have information about the genus of the surface $S$,
denoted by $g$ and the number of components
of the link $L$ , denoted by $\ell(L)$. In fact, the number of
components of $L$ is the number of faces in the Seifert
graph $\Gamma(S(L))$ on the surface $S$. By Euler's characteristic formula; $s(S)-c(S)+\ell(L)=2-2g$,
we get a different expression of the upper bound in Corollary~\ref{fpcor1}
$$fp(L)\le 2g+\ell(L)-1 + 2\lambda(T).$$

\subsection{Relations between plumbing numbers and classical link invariants}
\label{relation}

An easy observation for these plumbing numbers is the relations between three plumbing numbers.
The fundamental inequalities regarding these three plumbing numbers are

$$ bk(L) \le fpbk(L)~{\rm{and}}~ fp(L) \le fpbk(L).$$

For the second inequality, we consider the link $S(4,1)$ discussed
in Proposition 4.1. in~\cite{FHK:openbook}. It admits a flat plumbing
surface with the flat plumbing number $3$ as shown in Figure 4.2
in~\cite{FHK:openbook}. In fact, they claimed the flat plumbing number
of the link $S(4,1)$ is not less than nor equal to $3$. But, as we have mentioned in
Remark~\ref{remark2} (1), its flat plumbing number has to be nonzero because
it is a nontrivial link. By Remark~\ref{remark2} (2), its flat plumbing number has to
either $1$ or $3$ because the number of components of the link $S(4,1)$ is $2$.
It is easy to see that the link of the flat plumbing number $1$ must be a trivial
link of two components. They also showed the flat plumbing basket number of
the link $S(4,1)$ is bigger than or equal to $5$.
Therefore, the second inequality is proper for the link $S(4,1)$.

For the first inequality, let consider the trefoil knot. Example~\ref{basketexa}
shows its basket number is $2$ while its flat plumbing basket number is $4$
by Remark~\ref{remark3}. Then naturally we can ask the following question.

\begin{que} \label{que1}
Are there links for which the difference of plumbing numbers is
arbitrarily large ?
\end{que}

We conjecture that the link $L_{2n}$, the closure of $(\sigma_1)^{2n}\in B_2$
is one link for Question~\ref{que1} for the basket number and the flat plumbing basket number.
It is easy to see that its basket number is $1$. Showing its flat plumbing basket number
is bigger than $2n$ could be done using some classical invariants of links.

At least, this example $L_{2n}$ demonstrates that there exists a link of
different basket number and flat plumbing number. It is
fairly easy to see that any link whose flat plumbing number is less
than $3$ is trivial; either the unknot, trivial links of two or
three components.

We first thought the inequality $ bk(L) \le fp(L)$ is true; all theorems
in the article seems to work but these are upper bounds. Therefore, we can
not prove the inequality and there are not much classification theorems
about links of a fixed basket number or flat plumbing number while
there are some progress on links of a fixed flat plumbing basket number~\cite{CDK}.

For the last part of subsection, we compare these
plumbing numbers with the genus and canonical genus of a link.
Let us recall the definitions first. The \emph{genus} of
a link $L$ is the minimal genus among all Seifert
surfaces of $L$, denoted by $g(L)$. A Seifert surface $S$ of $L$
with the minimal genus $g(L)$ is called a \emph{minimal genus
Seifert surface} of $L$. A Seifert surface of $L$ is said to be
\emph{canonical} if it is obtained from a diagram of $L$ by applying
Seifert's algorithm. Then the minimal genus among all canonical
Seifert surfaces of $L$ is called the \emph{canonical genus} of $L$,
denoted by $g_c(L)$. A Seifert surface $S$ of $L$ is said to be
\emph{free} if the fundamental group $\pi_1(\mathbb{S}^3 \setminus  S)$
of the complement of $S$ is a free group. Then the minimal
genus among all free Seifert surfaces of $L$ is called the
\emph{free genus} for $L$, denoted by $g_f (L)$. Since any canonical
Seifert surface is free, we have the following inequalities.
$$g(L)\le g_f(L) \le g_c(L).$$
There are many interesting results about the above
inequalities; there are links of arbitrary large differences
for each inequalities~\cite{Brittenham:free, crowell:genera, KK, Moriah,
Nakamura, sakuma:minimal}.

From Corollary~\ref{fpcor1}, we obtain the following
corollary.

\begin{cor}
Let $S$ be a minimal genus canonical Seifert surface of a link $L$ with $s(S)$
Seifert circles and $c(S)$ half twisted bands and let $\ell(L)$ be the
number of components of $L$. Let $\Gamma$ be the
Seifert graph of a Seifert surface $S$. Let $T$ be a co-tree  edge alternating spanning tree
with the minimum $\lambda(T)$ as described in Corollary~\ref{fpcor1}.
Then,
$$2g(L)+\ell(L)-1 \le fp(L) \le 2 g_c (L)+ 2 \lambda(T) + \ell(L)-1.$$ \label{coroflat}
\begin{proof}
Let $v$, $e$ and $f$ be the numbers of vertices, edges and faces,
respectively in a $2$-cell embedding of $\Gamma$ onto $S$~\cite{GT1}.  From
Corollary~\ref{fpcor1}, we find
\begin{align*}fp(S)&\le c(S)-s(S)+1 + 2\lambda(T)=e-v+1+ 2\lambda(T)\\
&=-(v-e+f)+f+1 + 2\lambda(T)=2g_c(L)-2+ f+1+ 2\lambda(T)\\
&=2g_c(L)+ 2 \lambda(T)+ \ell(L)-1.\end{align*}
Since a flat plumbing surface is
also a Seifert surface of $L$ and for a connected minimal genus seifert surface
can be considered as a boundary of a bouquet of $n$-circles (recalled that $1-n+\ell(L)= 2-2g(L)$),
then $n\le fb(L)$ and thus the first inequality follows.
\end{proof}
\end{cor}

By applying Theorem~\ref{flatbktheorem3}, we obtain the following
corollary using a similar argument to that used in Corollary~\ref{coroflat}.

\begin{cor}
Let $\Gamma$ be an Seifert graph of a minimal genus
canonical Seifert surface $S$ of a link $L$ with
vertex set $V(\Gamma)$ and edge set $E(\Gamma)$, and
let $\ell(L)$ be the number of components of $L$.
Let $T$ be a co-tree edge alternating spanning tree of the induced graph $G(\Gamma)$
given in Theorem~\ref{maintheorem1}.
Let $\delta(T)$, $\zeta(T)$ and $\eta(T)$ be the numbers
defined in Theorem~\ref{flatbktheorem3}. Then,
$$2g(L)+\ell-1 \le fpbk(L) \le 2 g_c (L)+ 2(
 2\delta(T) + \zeta(T)+\eta(T) -|V(\Gamma)|) + \ell(L)-1.$$ \label{coroflatbk}
\begin{proof}
The first inequality follows by the same reason as explained in Corollary~\ref{fpcor1}.
From Theorem~\ref{flatbktheorem3}, we find
\begin{align*}
fpbk(L)&\le |E(\Gamma)|-3(|V(\Gamma)|-1) + 2 ( 2\delta(T)+\zeta(T) +\eta(T)) \\
&=(|E(\Gamma)|-|V(\Gamma)|+\ell(L)) -2|V(\Gamma)| + 3 -\ell(L) + 2 ( 2\delta(T)+\zeta(T) +\eta(T))\\
&=2g_c(L)-2|V(\Gamma)|+ 2( 2\delta(T)+\zeta(T) +\eta(T)-|V(\Gamma)|)+ \ell(L)+1.
\end{align*}
\end{proof}
\end{cor}

Let us remark that Hirose and Nakashima recently found a better lower bound for the
flat plumbing basket number of a knot using the genus and 
the Alexander polynomial of the knot~\cite{HN}.

\begin{theorem} (\cite{HN}) \label{thm:lowerbound}
Let $K$ be a non-trivial knot, $\Delta_K(t)$ be the Alexander polynomial of $K$,
$\deg \Delta_K(t) =$ (the maximal degree of $\Delta_K(t)$) $-$
(the minimal degree of $\Delta_K(t)$),
$a$ be the leading coefficient of $\Delta_K(t)$,
and $g(K)$ be the minimal genus of the Seifert surface (i.e. three genus)
of $K$. Then $fpbk(K)$ is evaluated as follows:
\begin{enumerate}
\item If $a = \pm 1$ then $fpbk(K) \geq
\max \{ 2g(K)+2, \deg \Delta_K(t) +2 \}$,
\item If $a \not= \pm 1$ then $fpbk(K) \geq
\max\{ 2g(K)+2, \deg \Delta_K(t) +4 \}$.
\end{enumerate}
\end{theorem}

\section*{Acknowledgments}
The referee's keen and thoughtful observations, which
led the article to take this current form, are thankfully acknowledged.
He also like to thank Hunki Baek, Younghae Do, Tsuyoshi Kobayashi,
Young Soo Kwon, Jaeun Lee and Myoungsoo Seo for
helpful discussion and Lee Rudolph for his comments. A critical error in the first version
was pointed out by Tsuyoshi Kobayashi. The author has learnt critical
facts about the Seifert graphs from Jaeun Lee and Young Soo Kwon.
The \TeX\, macro package
PSTricks~\cite{PSTricks} was essential for typesetting the equations
and figures.

\end{document}